\documentclass[12pt, leqno]{amsart}

\usepackage{graphicx}
\usepackage{amssymb,amsmath,amsthm,amscd,accents}
\usepackage{stackengine,scalerel}
\usepackage{mathrsfs}
\usepackage{lscape}
\usepackage{enumerate}
\usepackage{upgreek}
\usepackage{bbm}
\usepackage{stmaryrd}
\usepackage[usenames,dvipsnames]{color}
\usepackage[colorlinks=true, pdfstartview=FitV,
 linkcolor=blue,citecolor=blue,urlcolor=blue]{hyperref}
\usepackage{tikz}
\usetikzlibrary{calc,tikzmark,decorations.pathmorphing}
\tikzset{dynkdot/.style={circle,draw,scale=.38}}
\usetikzlibrary{arrows.meta,arrows}
\usepackage{bm}
\usepackage{marginnote}
\usepackage{verbatim}
\usepackage[normalem]{ulem}  
\usepackage{xspace}

\usepackage[all]{xy}

\allowdisplaybreaks[3]

\newcommand{\nc}{\newcommand}
\numberwithin{equation}{section}

\newenvironment{red}{\relax\color{red}}{\relax}
\newenvironment{blue}{\relax\color{blue}}{\hspace*{.5ex}\relax}
\newenvironment{purple}{\relax\color{blue}}{\hspace*{.5ex}\relax}

\newenvironment{magenta}{\relax\color{magenta}}{\hspace*{.5ex}\relax}
\newenvironment{jaune}{\relax\marginnote{\ber \scalebox{.6}{\sc{To be deleted}} \er}
 \color{Orchid}}{\hspace*{.5ex}\relax}

\newcommand{\bep}{\begin{purple}}
\newcommand{\eep}{\end{purple}}

\newcommand{\ber}{\begin{red}}
\newcommand{\er}{\end{red}}
\nc{\bj}{\begin{jaune}}
\nc{\ej}{\end{jaune}}
\newcommand{\beb}{\begin{blue}}
\newcommand{\eb}{\end{blue}}
\newcommand{\bgr}{\begin{jaune}}
\newcommand{\egr}{\end{jaune}}
\newcommand{\bem}{\begin{magenta}}
\newcommand{\eem}{\end{magenta}}

\nc{\eq}{\begin{myequation}} \nc{\eneq}{\end{myequation}}
\nc{\eqn}{\begin{myequationn}} \nc{\eneqn}{\end{myequationn}}

\newcommand{\bero}[1]{\begin{red}{}\marginnote{\fbox{\scshape\lowercase{O}}}%
#1}

\nc{\ms}{\mspace}

\nc{\qR}[1]{\ttq_{\mspace{-2mu}\raisebox{-.8ex}{${\scriptstyle{#1}}$}}}

\theoremstyle{plain}
\newtheorem{lemma}{Lemma}[section]
\newtheorem{proposition}[lemma]{Proposition}
\newtheorem{theorem}[lemma]{Theorem}

\newtheorem{corollary}[lemma]{Corollary}
\newtheorem{conjecture}[lemma]{Conjecture}

\theoremstyle{definition}
\newtheorem{remark}[lemma]{Remark}
\newtheorem{example}[lemma]{Example}

\newtheorem{definition}[lemma]{Definition}

\nc{\Def}{\begin{definition}}
\nc{\edf}{\end{definition}}
\nc{\mTh}{\begin{mtheorem}}
\nc{\enmth}{\end{mtheorem}}
\nc{\Th}{\begin{theorem}}
\nc{\enth}{\end{theorem}}
\nc{\Prop}{\begin{proposition}}
\nc{\enprop}{\end{proposition}}
\nc{\Lemma}{\begin{lemma}}
\nc{\enlemma}{\end{lemma}}
\nc{\Cor}{\begin{corollary}}
\nc{\encor}{\end{corollary}}
\nc{\Rem}{\begin{remark}}
\nc{\enrem}{\end{remark}}
\nc{\Conj}{\begin{conjecture}}
\nc{\enconj}{\end{conjecture}}

\newcommand{\Vkm}[1]{%
  \ttV\bigl(\hspace{0.2ex}\scalebox{.8}{\(\boxed{#1}\)}\hspace{0.2ex}\bigr)%
}

\nc{\ledot}{\mathrel{\le\ms{-11mu}\raisebox{.2ex}{$\cdot$}}}
\nc{\predot}{\mathrel{\preceq\ms{-9mu}\raisebox{.35ex}{$\centerdot$}}}

\newcommand{\rmat}[1]{{\mathbf r}_{\mspace{-2mu}\raisebox{-.5ex}{${\scriptstyle{#1}}$}}}

\renewcommand{\le}{\leqslant}
\renewcommand{\ge}{\geqslant}
\renewcommand{\preceq}{\preccurlyeq}

\newcommand{\head}{{\operatorname{hd}}}
\newcommand{\soc}{{\operatorname{soc}}}

\newcommand{\seteq}{\mathbin{:=}}
\newcommand{\conv}{\mathop{\mathbin{\mbox{\large $\circ$}}}}
\newcommand{\soplus}{\mathop{\mbox{\normalsize$\bigoplus$}}\limits}

\newcommand{\tens}{\mathop\otimes}

\newcommand{\ex}{\mathrm{ex}}
\newcommand{\fr}{\mathrm{fr}}

\newcommand{\g}{\mathfrak{g}}

\newcommand{\n}{\mathfrak{n}}

\newcommand{\C}{\mathbb{C}}
\newcommand{\Q}{\mathbb{Q}}
\newcommand{\Z}{\mathbb{Z}\ms{1mu}}

\newcommand{\be}{{\ms{1mu}\beta}}

\newcommand{\ve}{\boldsymbol{\varepsilon}}

\newcommand{\wt}{{\rm wt}}
\newcommand{\bx}{{\rm box}}

\newcommand{\codeg}{{\rm codeg}}
\newcommand{\trop}{{\rm trop}}

\newcommand{\Fold}{{\rm Fold}}
\newcommand{\re}{{\rm re}}

\newcommand{\ch}{{\rm ch}}

\newcommand{\de}{\mathfrak{d}}

\newcommand{\Dynkin}{\triangle}  

\newcommand{\im}{\imath}
\newcommand{\jm}{\jmath}

\newcommand{\ii}{ \textbf{\textit{i}}}
\newcommand{\jj}{ \textbf{\textit{j}}}

\newcommand{\Fib}{\operatorname{Fib}}

\newcommand{\HOM}{\mathrm{H{\scriptstyle OM}}}

\newcommand{\tC}{\widetilde{C}}

\newcommand{\tc}{\widetilde{c}}

\newcommand{\tx}{\widetilde{x}}
\newcommand{\ty}{\widetilde{y}}

\newcommand{\hC}{\widehat{C}}

\newcommand{\hF}{\widehat{F}}

\newcommand{\hI}{\widehat{I}}

\newcommand{\hL}{\widehat{L}}
\newcommand{\hM}{\widehat{M}}
\newcommand{\hN}{\widehat{N}}

\newcommand{\hQ}{\widehat{Q}}

\newcommand{\ha}{\widehat{a}}
\newcommand{\hb}{\widehat{b}}
\newcommand{\hc}{\widehat{c}}

\newcommand{\hg}{{\widehat{\g}}}

\newcommand{\hi}{\widehat{i}}

\newcommand{\hl}{\widehat{l}}

\newcommand{\hs}{\widehat{s}}

\newcommand{\hw}{\widehat{w}}
\newcommand{\hx}{\widehat{x}}
\newcommand{\hy}{\widehat{y}}
\newcommand{\hz}{\widehat{z}}

\newcommand{\hbe}{\widehat{\be}}

\newcommand{\hii}{{\widehat{\ii}}}
\newcommand{\hjj}{{\widehat{\jj}}}

\newcommand{\frakB}{\mathfrak{B}}
\newcommand{\frakC}{\mathfrak{C}}
\newcommand{\frakD}{\mathfrak{D}}
\newcommand{\frakE}{\mathfrak{E}}
\newcommand{\frakF}{\mathfrak{F}}

\newcommand{\frakT}{\mathfrak{T}}

\newcommand{\frakc}{\mathfrak{c}}

\newcommand{\frake}{\mathfrak{e}}

\newcommand{\frako}{\mathfrak{o}}

\newcommand{\tfrakc}{\widetilde{\frakc}}

\newcommand{\hfrakC}{\widehat{\frakC}}
\newcommand{\hfrakc}{\widehat{\frakc}}
\newcommand{\hfrake}{\widehat{\frake}}

\newcommand{\hfrakD}{\widehat{\frakD}}

\newcommand{\sfM}{\mathsf{M}}

\newcommand{\sfP}{\mathsf{P}}
\newcommand{\sfQ}{\mathsf{Q}}

\newcommand{\sfT}{\mathsf{T}}

\newcommand{\sfr}{\mathsf{r}}

\newcommand{\hsfQ}{\widehat{\sfQ}}

\newcommand{\hsfr}{\widehat{\sfr}}

\newcommand{\bbD}{\mathbb{D}}

\newcommand{\bbF}{\mathbb{F}}

\newcommand{\bbT}{\mathbb{T}}

\newcommand{\brho}{{\boldsymbol{\rho}}}

\newcommand{\hbfx}{\widehat{\bfx}}

\newcommand{\dbfy}{\dot{\bfy}}
\newcommand{\dy}{\dot{y}}
\newcommand{\hdbfy}{\widehat{\dbfy}}
\newcommand{\hdy}{\widehat{\dy}}

\newcommand{\hbfg}{\widehat{\bfg}}

\newcommand{\bfa}{\mathbf{a}}
\newcommand{\bfb}{\mathbf{b}}

\newcommand{\bfe}{\mathbf{e}}
\newcommand{\bff}{\mathbf{f}}
\newcommand{\bfg}{\mathbf{g}}
\newcommand{\bfh}{\mathbf{h}}
\newcommand{\bfi}{\mathbf{i}}

\newcommand{\bfr}{\mathbf{r}}

\newcommand{\bfv}{\mathbf{v}}
\newcommand{\bfw}{\mathbf{w}}
\newcommand{\bfx}{\mathbf{x}}
\newcommand{\bfy}{\mathbf{y}}

\newcommand{\calA}{\mathcal{A}}

\newcommand{\calE}{\mathcal{E}}
\newcommand{\calF}{\mathcal{F}}
\newcommand{\calG}{\mathcal{G}}

\newcommand{\calL}{\mathcal{L}}
\newcommand{\calM}{\mathcal{M}}

\newcommand{\calQ}{\mathcal{Q}}
\newcommand{\calR}{\mathcal{R}}

\newcommand{\calT}{\mathcal{T}}

\newcommand{\hcalF}{{\widehat{\calF}}}
\newcommand{\hcalG}{{\widehat{\calG}}}

\newcommand{\hcalM}{\widehat{\calM}}

\newcommand{\scrC}{\mathscr{C}}

\newcommand{\ttB}{\mathtt{B}}

\newcommand{\ttR}{\mathtt{R}}

\newcommand{\ttV}{\mathtt{V}}

\newcommand{\ttb}{\mathtt{b}}

\newcommand{\ttq}{\mathtt{q}}

\nc{\bg}{\sigma}

\newcommand{\httb}{\widehat{\ttb}}

\newlength{\mylength}
\newcommand*{\para}{%
  \rlap{\rotatebox{-30}{\rule[.05ex]{.4pt}{.77em}}}%
  \kern.04em%
  \rlap{\kern.36em\raisebox{0.649519052835em}{\rule{.6em}{.4pt}}}%
  \rule{.6em}{.4pt}\kern-.04em%
  \rotatebox{-30}{\rule[.05ex]{.4pt}{.77em}}}

\newcommand{\hDynkin}{{\widehat{\Dynkin}}}

\newcommand{\tDynkin}{{\widetilde{\Dynkin}}}

\newcommand{\diag}{\mathrm{diag}}

\newcommand{\qtq}[1][{and}]{\quad\text{{#1}}\quad}
\newcommand{\qtqh}[1][{and hence}]{\quad\text{{#1}}\quad}
\newcommand{\qtoq}[1][{or}]{\quad\text{{#1}}\quad}

\newcommand{\ee}{\end{enumerate}}
\newcommand{\bitem}{\begin{itemize}}
\newcommand{\eitem}{\end{itemize}}
\newcommand{\ben}{\begin{enumerate}[{\rm (1)}]}
\newcommand{\bnum}{\begin{enumerate}[{\rm (i)}]}
\newcommand{\bnump}{\begin{enumerate}[{\rm (i)$'$}]}
\newcommand{\bna}{\begin{enumerate}[{\rm (a)}]}
\newcommand{\bnA}{\begin{enumerate}[{\rm (A)}]}
\newcommand{\bc}{\begin{cases}}
\newcommand{\ec}{\end{cases}}
\newcommand{\ba}{\begin{array}}
\newcommand{\ea}{\end{array}}

\newenvironment{myequation}
{\relax\setlength{\arraycolsep}{1pt}\begin{eqnarray}}
{\end{eqnarray}}
\newenvironment{myequationn}
{\relax\setlength{\arraycolsep}{1pt}\begin{eqnarray*}}
{\end{eqnarray*}}

\nc{\eqs}[1]{\underset{\raisebox{.4ex}[.7ex][0ex]{$\scriptstyle{#1}$}}{=}}

\nc{\seq}[1]{ \boldsymbol{(} {#1} \boldsymbol{)}   }
\nc{\Bseq}[1]{ \boldsymbol{\left(} {#1} \boldsymbol{\right)}   }

\newcommand{\Dcan}{{\bbD_\can}}

\newcommand{\longepito}[1][]{\xymatrix@C=4ex{{}\ar@{->>}[r]^{#1}&{}}}

\newcommand{\uii}{{\underline{\boldsymbol{\im}}}}

\nc{\ake}[1][2ex]{\rule[-.5ex]{0ex}{#1}}

\newcommand{\ujj}{{\underline{\boldsymbol{\jm}}}}

\newcommand{\Red}{\mathrm{Red}}
\newcommand{\Seq}{\mathrm{Seq}}

\nc{\col}{\colon}
\nc{\ord}{\mathrm{ord}}
\nc{\catQ}{\scrC_\calQ}
\nc{\catD}{\scrC_\bbD}
\nc{\vs}{\vspace*}
\nc{\D}{\mathscr{D}}
\nc{\Proof}{\begin{proof}}
  \nc{\QED}{\end{proof}}

\nc{\tLa}{\widetilde{\Lambda}}
\nc{\hLa}{\widehat{\Lambda}}
\nc{\hSigma}{\widehat{\Sigma}}
\nc{\ro}{{\rm(}}
\nc{\rf}{{\rm)}\xspace}
\nc{\Aut}{\mathrm{Aut}}
\nc{\can}{\mathrm{can}}
\nc{\Dc}{{\bbD_{\can}}}
\nc{\Cgo}{\scrC_{\mathfrak{g}}^{0}}
\nc{\bb}{\mathtt{b}}
\nc{\catC}{\scrC}
\nc{\braid}{\ttB}
\nc{\Ass}[1][\bbD]{\calE_{#1}}
\nc{\Ci}{C^\uii}
\nc{\Cj}{C^\ujj}
\nc{\Di}[1][{\uii}]{{\bbD,{#1}}}
\nc{\Vi}[1][\uii]{V^{#1}}
\nc{\Pii}[1][\uii]{P^{#1}}
\nc{\Vdi}[1][{\Di}]{V^{#1}}
\nc{\Pdi}[1][{\Di}]{P^{#1}} 
\nc{\PBix}[1][{[l,r]}]{\Z_{\ge0}^{\oplus{#1}}}
\nc{\mcf}[1][{$[a,b]$}]{maximal commuting family of $\bfi$-boxes in {#1}\xspace}
\nc{\nn}{\nonumber}
\nc{\bwr}{\mbox{\large$\wr$}}
\newcommand{\tch}{\widetilde{\ch}}
\nc{\tchDcan}{\ms{2mu}\tch_{\mspace{.1mu}\raisebox{-.4ex}{${\scriptstyle{\Dcan}}$}}}
\nc{\monoto}[1][]{\xymatrix@C=2ex{\ar@{>->}[r]^-{{#1}}&}\ms{-8mu}}

\usepackage[colorinlistoftodos]{todonotes}

\title[Folding of \(\de\)-invariants] {Folding of \(\de\)-invariants via \(F\)-invariants\\
with applications to quiver Hecke KR-modules}

\author[S.-j. Oh]{Se-jin Oh}
\thanks{ The research of S.-j.\ Oh was supported by the National Research Foundation of
  Korea (NRF) Grant funded by the Korea government (MSIT) (NRF-2022R1A2C1004045).}
\address[S.-j. Oh]{ Department of Mathematics, Sungkyunkwan University, Suwon, South Korea}
\email[S.-j. Oh]{sejin092@gmail.com}

\date{\today}

\begin{document}

\begin{abstract}
The \(\de\)-invariant measures the interaction between two simple
modules and controls the behavior of their convolution product.
Using categorical interpretations of Cao's tropical and
\(F\)-invariants, we establish folding formulas for the
\(\Lambda\)-invariants between determinantial modules in arbitrary
symmetrizable type and derive the corresponding \(\de\)-formula
by symmetrization.
 We apply the \(\de\)-formula to quiver Hecke Kirillov--Reshetikhin modules 
 and show that, in every finite type, their
\(\de\)-invariants are completely determined by the coefficients
of the symmetrized inverse \(t\)-quantized Cartan matrix. 
We obtain explicit formulas in all finite types, including closed
formulas in the classical types \(A\), \(B\), \(C\), and \(D\).
\end{abstract}

\maketitle
\tableofcontents

\section*{Introduction}
\label{sec: introduction}

Cluster algebras and their quantum analogues were introduced as
classes of algebras equipped with distinguished generators and
mutation relations \cite{FZ02,BZ05}.  They have since appeared
naturally in representation theory, where cluster monomials are
often expected to correspond to distinguished basis elements.
Monoidal categorification provides a categorical form of this
relationship: cluster variables are represented by real prime simple
modules, cluster monomials by tensor products of pairwise
strongly commuting real prime simple modules, and exchange relations
by short exact sequences \cite{HL10}.

\smallskip

Quiver Hecke algebras were introduced independently in
\cite{KL1,KL2,R08} to categorify the negative half of a quantum
group. Their renormalized \(R\)-matrix calculus for
simple modules gives rise to the numerical invariants \(\Lambda\)
and \(\de\).  These
invariants have played a central role in proving monoidal
categorification conjectures \cite{KKKO18}.  In particular, \(\Lambda(M,N)\) 
records the \emph{degree} of the renormalized \(R\)-matrix from \(M\conv N\) to \(N\conv M\), whereas the symmetric invariant
\[
\de(M,N) = \frac{\Lambda(M,N)+\Lambda(N,M)}{2}
\]
measures the interaction between \(M\) and \(N\).  Under the usual realness assumptions, it detects strong commutation and controls
the structure of their convolution product (\cite{KKKO15,KKKO18}).  Thus formulas for \(\de\)-invariants provide concrete information on
the multiplicative structure of these monoidal categories.

\smallskip

Despite the simplicity of their definition, \(\de\)-invariants
are generally quite difficult to compute, since explicit
descriptions of simple modules, and often even their dimensions,
are unavailable.
Generalized quantum affine Schur--Weyl duality provides a way
to compute these invariants in certain cases.
Our approach is to compute them through tropical cluster
invariants and folding.

\smallskip

Let \(w\) be a Weyl group element and choose a reduced 
sequence \(\ii\in\Red(w)\).  In the category \(\scrC_w\) of 
quiver Hecke modules associated with \(w\), every simple 
object is obtained, uniquely up to a grading shift, as the 
head of an ordered convolution product of the \emph{cuspidal} 
modules associated with \(\ii\)~\cite{Kato14,McNa15,TW16}.  The \(\de\)-invariants 
between these cuspidal modules govern
their commutation relations, while the \(\Lambda\)-invariants
play a central role in determining the \emph{normality}
of the corresponding sequences. Moreover, \emph{determinantial modules} provide categorical counterparts
of \emph{quantum minors} \cite{KKOP18}. The determinantial module associated with an \(\bfi\)-box, also called a box module, is the unique simple head of the ordered convolution product of the cuspidal modules indexed by the occurrences of its color in that box.  The classes of these modules correspond to cluster variables in box seeds for the quantum unipotent coordinate algebra \(A_q(\n(w))\) \cite{KKOP24A,KK24}, whose quantum cluster algebra structure was established in \cite{GLS13,GY17,GY21}.

\smallskip

A parallel picture appears in the representation theory of quantum
affine algebras $U_q'(\g^{(1)})$ \cite{CP95A,Kas02,KKOP23P}.  
Fundamental representations play the role of PBW generators on the quantum affine side, 
and generalized quantum affine Schur--Weyl duality identifies the corresponding quiver
Hecke cuspidal modules with them \cite{KKK18}.  A ubiquitous
combinatorial object in this theory is  the inverse
\(C(q)^{-1}\) of the quantum Cartan matrix \(C(q)\): its
coefficients occur in the study of \(q\)-characters, quantum
Grothendieck rings, the commutative part of universal
\(R\)-matrices, and denominator formulas for normalized
\(R\)-matrices \cite{FR99,Nak04,HL15,FO21,FM21}.  In this sense, the
inverse quantum Cartan matrix provides a common key to  several
structures in the representation theory of $U_q'(\g^{(1)})$; see
\cite[Introduction]{FM21} for this point of view.

\smallskip

An especially important family in the  representation
theory of quantum affine algebras is formed by the Kirillov--Reshetikhin modules 
(see~\cite{Nak03KR,Her06KR,CH10Beyond} and the references therein).  
In simply-laced finite type, the Schur--Weyl duality associated with a
Dynkin quiver $Q$ realizes the corresponding quiver Hecke
modules as determinantial modules.  Their \(\de\)-invariants
may therefore be studied through both quiver Hecke and quantum
affine representation theory (\cite{KKOP23P}).

\smallskip

On the cluster-algebraic side, Cao introduced seed-independent
tropical and \(F\)-invariants defined in terms of pointed Laurent
expansions \cite{Cao23}.  The \(F\)-invariant is the
symmetrization of the tropical invariant.
When a cluster algebra admits a suitable \emph{monoidal}
categorification, the \emph{Laurent maximum property} gives
categorical interpretations of these invariants:
\begin{align}\label{eq: interpretation}
\bigl\langle\Phi([M]),\Phi([N])\bigr\rangle_{\mathrm{trop}}
=\Lambda(N,M),\qquad 
\bigl(\Phi([M])\Vert\Phi([N])\bigr)_F
=2\de(M,N).
\end{align}
This provides a direct bridge between tropical cluster data and
renormalized \(R\)-matrices.   For a general categorification, however,
these interpretations of the tropical and \(F\)-invariants are
\emph{not} automatic.  We
obtain such an interpretation~\eqref{eq: interpretation} by applying
the maximum formula for Laurent families in
\cite[Theorem~2.13]{KKOP26} to categorical box seeds
(Theorem~\ref{thm: categorical interpretation of F-invariant}). We then use cluster-algebra folding \cite{JO26}
to derive the corresponding identities in the folded setting
from those for the upstairs box modules. 

\medskip 

We briefly fix the folding setting used below.  Let \(\hg\) be a
simply-laced Kac--Moody algebra with vertex set \(\hI\), equipped
with an admissible action of a finite group \(\Gamma\), and let \(\g\)
be the symmetrizable Kac--Moody algebra associated with the folded
Cartan datum.  Thus
\begin{align} \label{eq: frame}
\begin{array}{c}
(\hg,\hI)\qquad \ \ \ \text{simply-laced upstairs data}
\\[2pt]
\Big\downarrow\;/\Gamma
\\[2pt]
(\g,I=\hI/\Gamma)\qquad\text{symmetrizable downstairs data}.
\end{array}
\end{align}
Set
\[
\pi:\hI\longrightarrow I,
\qquad
I_i\seteq\pi^{-1}(i).
\]
We refer to the data associated with \(\hg\) and \(\g\) as the
\emph{upstairs} and \emph{downstairs} data, respectively, and
usually denote upstairs objects by hatted symbols.  We normalize
the symmetrizer
\[
D^\g=\diag(d_i\mid i\in I)
\quad
\text{ so that }
\quad
d_i|I_i|=m \ \  \text{is independent of \(i\).}
\]
  Unless explicitly stated otherwise, no
finite-type assumption is imposed.

\smallskip

The framework in~\eqref{eq: frame} is
provided by  \(\bfi\)-boxes and \(\bfi\)-box seeds.  
For a reduced sequence \(\ii\in\Red(w)\), each \(\bfi\)-box indexes a determinantial module in
\(\scrC_w\), or equivalently a unipotent quantum minor.  Maximal commuting families of \(\bfi\)-boxes give
cluster seeds, while box moves describe the mutations between them
\cite{KKOP24A,KK24} (see
Lemma~\ref{lem: two types of box moves},
Proposition~\ref{prop: box moves are categorical mutations}, and
Theorem~\ref{thm: categorical realization of box seeds}).  Thus
\(\bfi\)-box seeds provide a common
language for determinantial modules, cluster variables, and their
Laurent expansions. 

\smallskip

This box realization is particularly well suited to \emph{folding}.
The cluster-algebra folding developed in \cite{JO26} organizes the
upstairs \(\bfi\)-boxes into fibers and folds the associated box seeds
and cluster variables to their downstairs counterparts. Using this
folding, we derive the corresponding \(F\)-polynomial identity
(Proposition~\ref{prop: folding of pointed F-parts}) and recover the
tropical and \(F\)-invariants between downstairs \(\bfi\)-box cluster
variables from their simply-laced counterparts. At the
cluster-algebraic level, these identities hold in the general
symmetrizable setting without a finite Cartan type assumption.

\smallskip

The categorical folding formulas do not require a monoidal
functor between the upstairs and downstairs categories.

\smallskip

The same folding philosophy also appears in the theory of
\(t\)-quantized Cartan matrices $C(t)$.   In finite Cartan type, the works
\cite{KO23,KO24,JLO1,JLO2} show that \(t\)-quantized Cartan
matrices and their symmetrized inverses encode important representation-theoretic information and
behave naturally under folding.  This suggests that, for symmetrizable Cartan data not necessarily of
finite type, the symmetrized inverse \(t\)-quantized Cartan matrix
should play a fundamental role in the representation theory of
quiver Hecke algebras, analogous to that played by the inverse
quantum Cartan matrix in the representation theory of $U_q'(\g^{(1)})$.

\smallskip

We now state our main results more precisely.  Our first result is purely cluster-algebraic and requires neither
reducedness nor finite Cartan type.
For a positive braid \(\ttb\), let \(\ii,\jj\in\Seq(\ttb)\), and let \(\frakC\) and
\(\frakD\) be admissible chains associated with \(\ii\) and
\(\jj\), respectively. Let
\[
\frakc\in\calF(\frakC),
\qquad
\frakc'\in\calF(\frakD).
\]
For any
\(\hfrakc'\in\Fib_{\frakD}(\frakc')\), we prove that
\[
\begin{aligned}
\langle x_\frakc,x_{\frakc'}\rangle_{\trop}
=d_{i_\frakc}
\langle\hx_{\Fib_{\frakC}(\frakc)},\hx_{\hfrakc'}\rangle_{\trop},\qquad
\langle x_{\frakc'},x_\frakc\rangle_{\trop}
=d_{i_\frakc}
\langle\hx_{\hfrakc'},\hx_{\Fib_{\frakC}(\frakc)}\rangle_{\trop},
\end{aligned}
\]
where
\(
\hx_{\Fib_{\frakC}(\frakc)}
\seteq \prod_{\hfrakc\in\Fib_{\frakC}(\frakc)}\hx_{\hfrakc}
\).
Adding these two ordered identities gives
\[
\bigl(x_\frakc\Vert x_{\frakc'}\bigr)_F
=  d_{i_\frakc} \bigl( \hx_{\Fib_{\frakC}(\frakc)} \Vert \hx_{\hfrakc'}
\bigr)_F,
\]
(Theorem~\ref{thm: folding of F-invariants between box variables};
see also
\eqref{eq: folding of first ordered tropical pairing}--\eqref{eq: folding of F-invariants between box variables}).
The right-hand sides are independent of the choice of
\(\hfrakc'\).  Thus both ordered tropical pairings, and hence their
sum defining the \(F\)-invariant, are recovered from the lifts.

\smallskip

Let \(\ii,\jj\in\Red(w)\), and let \(M^\ii[\frakc]\) and
\(M^\jj[\frakc']\) be the determinantial modules associated with
two downstairs \(\bfi\)-boxes.  Choose
\(\hfrakc'\in\Fib_{\frakD}(\frakc')\).
Our categorical interpretation of the tropical invariant gives
folding formulas for both orders of the \(\Lambda\)-invariant
(Theorem~\ref{thm: categorical folding of directional Lambda invariants}).
For these box modules, which may be attached to different reduced
words, the formulas are
\[
\begin{aligned}
\Lambda\bigl(M^\ii[\frakc],M^\jj[\frakc']\bigr)
&=d_{i_\frakc}\sum_{\hfrakc\in\Fib_{\frakC}(\frakc)}
\Lambda\bigl(M^\hii[\hfrakc],M^\hjj[\hfrakc']\bigr),\\
\Lambda\bigl(M^\jj[\frakc'],M^\ii[\frakc]\bigr)
&=d_{i_\frakc}\sum_{\hfrakc\in\Fib_{\frakC}(\frakc)}
\Lambda\bigl(M^\hjj[\hfrakc'],M^\hii[\hfrakc]\bigr),
\end{aligned}
\]
as stated in Corollary~\ref{cor: folding of d-values between box modules}.
Thus, symmetrizing the two identities gives the corresponding folding formula for
\(\de\)
(Theorem~\ref{thm: categorical folding of d-invariants}):
\[
\de\bigl(M^\ii[\frakc],M^\jj[\frakc']\bigr)= 
d_{i_\frakc} \sum_{\hfrakc\in\Fib_{\frakC}(\frakc)}
\de\bigl(M^\hii[\hfrakc],M^\hjj[\hfrakc']\bigr).
\]
All three expressions on the right are independent of the
chosen lift \(\hfrakc'\). Hence the
interaction between two determinantial modules downstairs is
\emph{completely} determined by the interactions among their simply-laced
lifts.  While the cluster-algebraic formula requires neither reducedness nor
finite Cartan type, the categorical formula requires reducedness but remains
valid for arbitrary symmetrizable Cartan data 
(Corollary~\ref{cor: folding of d-values between box modules} and~\eqref{eq: folded d-value between box modules}).

\smallskip

We next establish, without assuming finite Cartan type, a parallel
folding formula for the symmetrized inverse \(t\)-quantized Cartan
matrix.  Write
\[
\tC^\g(t)\seteq D^\g C^\g(t)^{-1},
\]
and define \(\tC^\hg(t)\) similarly.  For \(i,j\in I\), choose
\(\im\in I_i\) and \(\jm\in I_j\).  We prove that
\[
\tC_{i,j}^{\g}(t)
=
d_i\sum_{\im'\in I_i} \tC_{\im',\jm}^{\hg}(t)
= \frac{m}{|\Gamma|} \sum_{g\in \Gamma} \tC_{g(\im),\jm}^{\hg}(t).
\]
The expressions on the right are independent of the chosen lifts
(Theorem~\ref{thm: folding of inverse t-Cartan matrices} and~\eqref{eq: folding inverse t-Cartan first fiber}).
For the standard finite foldings
\[
A_{2n-1}\longrightarrow C_n,\qquad
D_{n+1}\longrightarrow B_n,\qquad
E_6\longrightarrow F_4,\qquad
D_4\longrightarrow G_2,
\]
one has \(m=|\Gamma|\), and hence the final expression is simply the
group sum.

\smallskip

In finite type, \cite[Theorem~8.2]{KO24} identify the coefficients
of $\tC(t)$ with \(\de\)-invariants between \(Q\)-cuspidal modules.  Under this
interpretation, the specialization of our categorical folding
formula to cuspidal modules  gives the corresponding coefficient identities in~\eqref{eq: folding inverse t-Cartan first fiber}.  In particular, the same factor \(d_i\) appears in the categorical
folding formula and in the folding formula for the
symmetrized inverse \(t\)-Cartan matrix (Theorem~\ref{thm: finite-type cuspidal d-values}).

\smallskip

This agreement leads to our principal finite-type application.  Let
\(\g\) be of finite type, and let
\(\Vkm{i^k}_p\) and \(\Vkm{j^\ell}_s\) be (\(Q\)-adapted) \emph{quiver
Hecke Kirillov--Reshetikhin modules}~\eqref{eq: KR-module over quiver Hecke}.  Assume that all the vertices
defining these two modules lie in the Auslander-Reiten quiver \(\Upsilon_Q\).  Set
\[
K\seteq\max(k,\ell), \qquad L\seteq\min(k,\ell), \qquad u\seteq|p-s|.
\]

Let \(h\) be the Coxeter number of $\g$, and set
\[
\tC_{i,j}^{\g,[h]}(t) \seteq \sum_{v=1}^{h}\tc_{i,j}^{\g}(v)t^v,
\qquad
[L]_t \seteq t^{L-1}+t^{L-3}+\cdots+t^{1-L}.
\]
Writing \([t^u]f(t)\) for the coefficient of \(t^u\) in \(f(t)\),
we have
\begin{align} \label{eq: closed in intro}
\de\bigl(\Vkm{i^k}_p,\Vkm{j^\ell}_s\bigr)
= [t^u]\!\left(t^K[L]_t\,\tC_{i,j}^{\g,[h]}(t)\right)    
\end{align}
(Theorems~\ref{thm: ADE quantum affine KR d-values}
and~\ref{thm: BCFG quiver Hecke KR d-values}; see
\eqref{eq: ADE quantum affine KR d-value}
and~\eqref{eq: BCFG quiver Hecke KR d-value}). 
Thus the \(\de\)-invariants between quiver Hecke
KR-modules can be read \emph{completely} from $\tC(t)$ in every finite type.  Our results
therefore \emph{support} the expectation that this matrix plays a
fundamental role in the representation theory of quiver Hecke
algebras.

\smallskip

The closed formulas~\eqref{eq: closed in intro} are particularly transparent for the
folding \(A_{2n-1}\to C_n\), where
\(\Gamma\simeq\mathbb Z/2\mathbb Z\) acts by
\(i\leftrightarrow 2n-i\).  For \(r,L\in\Z_{>0}\) and \(v\in\Q\), set
\[
\sfT^{(1)}_{r,L}(v) = \delta(v \in \Z)\big[\min(v,r,L,r+L-v)\big]_+, \quad
\kappa\seteq|k-\ell|, \quad u\seteq|p-s|,
\]
where \([a]_+=\max(a,0)\) and \(\delta(P)\) denotes the indicator
of a statement \(P\). 

In type \(A_n\), set
\(\brho^{(n)}_{i,j}\seteq\min(i,j,n+1-i,n+1-j)\).
Then we have
\[
\de\bigl(\Vkm{i^k}_p,\Vkm{j^\ell}_s\bigr)
= \sfT^{(1)}_{\brho^{(n)}_{i,j},L} \left( \frac{u-\kappa-|i-j|}{2} \right)
\]
(Corollary~\ref{cor: type A closed d-value formula}; see
also \eqref{eq: type A closed d-value formula}).

The folding \(A_{2n-1}\longrightarrow C_n\) makes the corresponding
type \(C_n\) formula equally explicit. 
If \(1\leq i,j<n\), then
\[ 
\de\bigl(\Vkm{i^k}_p,\Vkm{j^\ell}_s\bigr)
=\sfT^{(1)}_{\brho^{(2n-1)}_{i,j},L}
\left( \frac{u-\kappa-|i-j|}{2} \right)
+\sfT^{(1)}_{\brho^{(2n-1)}_{2n-i,j},L}
\left( \frac{u-\kappa-(2n-i-j)}{2} \right).
\]
If \(i=n\) or \(j=n\), where \(n\) is the long node, then
\[
\de\bigl(\Vkm{i^k}_p,\Vkm{j^\ell}_s\bigr)
= 2\sfT^{(1)}_{\brho^{(2n-1)}_{i,j},L}
\left( \frac{u-\kappa-|i-j|}{2}
\right).
\]
The two summands in the first formula arise from the upstairs pairs
of colors \((i,j)\) and \((2n-i,j)\), respectively, while the factor
\(2\) in the second is the symmetrizer attached to the fixed vertex
\(n\).  Here
\(\brho^{(2n-1)}_{i,j}=\brho^{(2n-1)}_{2n-i,j}=\min(i,j)\)
for \(1\leq i,j<n\).  Thus the closed type
\(C_n\) formula visibly retains the folding from type \(A_{2n-1}\)
by Corollary~\ref{cor: type C quiver Hecke KR closed formula}.

Analogous closed formulas are obtained in types \(B\) and \(D\) 
(Corollaries~\ref{cor: type B quiver Hecke KR closed formula}
and~\ref{cor: type D closed d-value formulas}).
The uniform coefficient formula also applies to the exceptional
types \(F\) and \(G\), as illustrated by explicit examples 
(Examples~\ref{ex: type F4 folded KR d-value}
and~\ref{ex: type G2 folded KR d-value}).  In particular, we obtain
new closed formulas for quiver Hecke KR-modules in types \(B\) and
\(C\).
Together with the formula in
Proposition~\ref{prop: KKOP directional box Lambda formula}, these
\(\de\)-formulas allow one to compute both orders of the
\(\Lambda\)-invariant for KR-modules defined by boxes in the same
\(Q\)-adapted reduced sequence (see Remark~\ref{rmk: Lambda computable} and Example~\ref{rmk: Lambda computable}).

\smallskip

 Although \eqref{eq: closed in intro} is proved here only in
finite type, the coefficient formula using the full \(t\)-adic matrix
also holds in three special cases for arbitrary symmetrizable
Cartan data.   For boxes in a \(Q\)-adapted reduced sequence, the same coefficient
expression using the full \(t\)-adic symmetrized inverse \(t\)-Cartan
matrix gives \(d_i\) both for an \(\bfi\)-box and its adjacent singleton
\(\bfi\)-box and for the exchange pair of boxes in a \(T\)-system. 
It vanishes for quiver Hecke KR-modules corresponding to
commuting \(\bfi\)-boxes.  In type \(B_3^{(1)}\), we illustrate this vanishing with a 
pair of non-nested commuting \(\bfi\)-boxes of different colors.  Thus \(\tC(t)\)
detects these categorical interactions beyond finite
type; see Proposition~\ref{prop: local coefficient formula beyond finite type}
and Example~\ref{ex: affine B3 commuting boxes of different colors}. 

\smallskip

The folding \(A_5\longrightarrow C_3\) illustrates the
\(F\)-invariant, categorical, and inverse \(t\)-Cartan formulas
throughout the paper.

\subsection*{Organization}
The paper is organized as follows.
Section~\ref{sec: F-invariant} recalls Cao's tropical and
\(F\)-invariants.
Section~\ref{sec: folding cluster algebras} studies the folding
of \(\bfi\)-box seeds, \(F\)-polynomials, and tropical invariants.
Section~\ref{sec: categorical folding} establishes categorical interpretations of Cao's invariants and the \(\Lambda\)- and \(\de\)-folding
formulas.
Finally, Section~\ref{sec: KR modules and inverse} introduces
quiver Hecke KR-modules, proves the folding formula for
symmetrized inverse \(t\)-Cartan matrices, derives finite-type
\(\de\)-formulas, and studies three special cases 
of the coefficient formula beyond finite type.

\subsection*{Convention}
Throughout this paper, we use the following conventions.
\bitem
\item For a statement $P$, we set $\delta(P)$ to be $1$ or $0$ depending on whether $P$ is true or not. In particular, $\delta_{i,j}=\delta(i=j)$
is the Kronecker delta.
\item For integers $a\le b$, we write $[a,b]=\{a,a+1,\ldots,b\}$.
We understand $[a,b] = \emptyset$ if $b<a$ and set $[a]\seteq [a,a]$.
\item For a power series $f(t) = \sum_{i=0}^\infty a_i t^i$ over a commutative ring $A$, we set 
$$
[t^n] f(t) \seteq a_n \quad \text{ for $n\in \Z_{\ge 0}$.}
$$
\eitem

\section{Tropical and \(F\)-invariants}
\label{sec: F-invariant}

In this section, we recall Cao's tropical and \(F\)-invariants
\cite{Cao23} for compatibly pointed elements of a
\(\Uplambda\)-cluster algebra.
Unless otherwise stated, neither finite type nor symmetry of the
underlying Cartan datum is assumed.
Their categorical interpretations in terms of the
\(\Lambda\)- and \(\de\)-invariants for box seeds will be
established in Section~\ref{sec: categorical folding}.

We follow \cite{FZ02,BZ05,FWZ1} for compatible
pairs and \(\Lambda\)-cluster patterns, and
\cite{Qin17,Cao23,KK19,KKOP23F} for dominance orders and pointed elements.

Let \(J=J_{\ex}\sqcup J_{\fr}\) be a finite index set, where
\(J_{\ex}\) and \(J_{\fr}\) are the exchangeable and frozen indices,
respectively. Let 
\bitem
\item $B=(b_{i,j})_{i\in J,\,j\in J_{\ex}}$ \footnote{An extended
exchange matrix is commonly denoted by \( \widetilde B \).  We write
\(B\) for the extended exchange matrix and \(B^{\mathrm{prin}}\) for its principal part.}
be an integer matrix whose principal part
\[
B^{\mathrm{prin}} \seteq
(b_{i,j})_{i,\,j\in J_{\ex}} \ \ \text{is skew-symmetrizable,}
\]

\item  $L=(l_{i,j})_{i,j\in J}$ be a skew-symmetric integer matrix.  
\eitem 
The pair \((L,B)\) is
\emph{compatible of type \(S=\diag(s_j\mid j\in J_{\ex})\)} if
\begin{equation}\label{eq: compatible pair}
B^{T}L=(S\mid0) \qquad   (s_j\in2\Z_{>0}).
\end{equation}
Here the block decomposition on the right-hand side is taken with
respect to \(J=J_{\ex}\sqcup J_{\fr}\).

Set $\bbF\seteq\C(X_j\mid j\in J)$, where the \(X_j\) are algebraically independent.  A \emph{\(\Uplambda\)-seed of type \(S\)} is a triple
\[
\Sigma=(\bfx,L,B), \qquad \bfx=(x_j)_{j\in J},
\]
where \(\bfx\) is algebraically independent over \(\C\),
\(\bbF=\C(x_j\mid j\in J)\), and \((L,B)\) is compatible of type \(S\). 

For \(k\in J_{\ex}\), its mutation
\[
\mu_k(\Sigma) = \bigl(\mu_k(\bfx),\mu_k(L),\mu_k(B)\bigr)
\]
is defined by the usual seed and compatible-pair mutations.  The
mutated cluster \(\mu_k(\bfx)=(x'_j)_{j\in J}\) is given by
\[
x'_j=x_j\quad(j\ne k), \qquad
x_kx'_k = \prod_{i\in J}x_i^{[b_{i,k}]_+} + \prod_{i\in J}x_i^{[-b_{i,k}]_+}.
\]
In particular, \(\mu_k(L)\) and \(\mu_k(B)\) are again compatible
of type \(S\). We call \(x'_k\) the \emph{cluster variable} obtained by mutation in
direction \(k\).

A \emph{\(\Uplambda\)-cluster pattern of type \(S\)} is a family of
\(\Uplambda\)-seeds
\[
\Sigma_t=(\bfx_t,L_t,B_t) \qquad(t\in\bbT),
\]
indexed by the vertices of the \(J_{\ex}\)-regular tree \(\bbT\),
such that adjacent seeds are related by the mutation labeling the
corresponding edge.  For each \(t\in\bbT\), let
\[
\calL(\Sigma_t) \seteq \C[x_{j;t}^{\pm1}\mid j\in J] \subset\bbF
\]
be the Laurent polynomial ring associated with \(\Sigma_t\).
The associated \emph{\(\Uplambda\)-cluster algebra} is the
\(\C\)-subalgebra\footnote{Cluster algebras are usually defined over
\(\Z\); here we work over \(\C\) for compatibility with the
complexified Grothendieck rings used below.}
\[
\calA \seteq \C[x_{j;t}\mid j\in J,\ t\in\bbT] \subset\bbF.
\]
By the Laurent phenomenon~\cite{FZ02}, one has
\begin{align} \label{eq: Laurent phe}
\calA\subset\calL(\Sigma_t) \qquad(t\in\bbT) 
\qtqh 
\calA \subset  \bigcap_{t\in\bbT}\calL(\Sigma_t) \subset\bbF.
\end{align}

An element of the form
\[
\bfx_t^{\bfa}\seteq\prod_{j\in J}x_{j;t}^{a_j}
\qquad
\bigl(\bfa=(a_j)_{j\in J}\in\Z_{\geq0}^{\oplus J}\bigr)
\]
is called a \emph{cluster monomial}.

For \(t\in\bbT\) and \(j\in J_{\ex}\), set \footnote{These Laurent monomials are often denoted by
\(\widehat y_{j;t}\); we use the dot notation since hats are reserved
for upstairs objects throughout this paper.}
\[
\dy_{j;t}\seteq\prod_{i\in J}x_{i;t}^{b_{i,j;t}},
\qquad
\dbfy_t\seteq(\dy_{j;t})_{j\in J_{\ex}}.
\]
For \(\bfg,\bfg'\in\Z^{\oplus J}\), we write
\[
\bfg'\preceq_t\bfg \qquad \text{if $\bfg'=\bfg+B_t\bfv$ for some $\bfv\in\Z_{\geq0}^{\oplus J_{\ex}}$.}
\]
Equivalently,
\[
\bfx_t^{\bfg'}=\bfx_t^{\bfg}\dbfy_t^{\bfv},
\quad \text{where} \quad
\bfx_t^{\bfg}\seteq\prod_{i\in J}x_{i;t}^{g_i}
\qtq 
\dbfy_t^{\bfv}\seteq\prod_{j\in J_{\ex}}\dy_{j;t}^{v_j}.
\]
The compatibility condition~\eqref{eq: compatible pair} implies that the columns of \(B_t\)
are linearly independent, and hence \(\preceq_t\) is a partial order.

By the Laurent phenomenon~\eqref{eq: Laurent phe}, every \(u\in\calA\) has a finite
expansion
\begin{align} \label{eq: u expression}
u=\sum_{\bfh\in\Z^{\oplus J}}c_{\bfh}\bfx_t^{\bfh}, \qquad c_{\bfh}\in\C.    
\end{align}
Following the terminology of pointed and copointed elements, we say
that \(u\) is \emph{pointed with respect to \(t\)} if there exists a
unique vector \(\deg_t(u)\in\Z^{\oplus J}\) such that
\[
c_{\deg_t(u)}=1, \qquad c_{\bfh}\ne0\Longrightarrow\bfh\preceq_t\deg_t(u).
\]
Dually, \(u\) is \emph{copointed with respect to \(t\)} if there
exists a unique vector \(\codeg_t(u)\in\Z^{\oplus J}\) such that
\[
c_{\codeg_t(u)}=1, \qquad c_{\bfh}\ne0\Longrightarrow\codeg_t(u)\preceq_t\bfh.
\]
We call \(\deg_t(u)\) and \(\codeg_t(u)\) the \emph{degree} and
\emph{codegree} of \(u\) with respect to \(t\), respectively.  An
element is \emph{bipointed with respect to \(t\)} if it is both
pointed and copointed with respect to \(t\).

Since the columns of \(B_t\) are linearly independent, every
pointed element \(u\) has a unique expansion
\begin{equation}\label{eq: pointed F expansion}
u=\bfx_t^{\deg_t(u)}F_{u;t}(\dbfy_t).
\end{equation}
 Here 
\[
F_{u;t}(\bfy)\in\C[y_j\mid j\in J_{\ex}], \qquad F_{u;t}(0)=1.
\]
We call \(F_{u;t}(\bfy)\) the \emph{\(F\)-polynomial of \(u\) with
respect to \(t\)}.

For \(j\in J_{\ex}\), let \(f_{j;u;t}\) be the maximal exponent of
\(y_j\) among the monomials occurring in \(F_{u;t}(\bfy)\), and set
\[
\bff_{u;t}\seteq(f_{j;u;t})_{j\in J_{\ex}} \in \Z_{\ge0}^{J_\ex}.
\]
We call \(\bff_{u;t}\) the \emph{\(f\)-vector of \(u\) with respect
to \(t\)}.
Then \(u\) is bipointed with respect to \(t\) if and only if
\(\bfy^{\bff_{u;t}}\) occurs in \(F_{u;t}(\bfy)\) with coefficient
\(1\).  In this case,
\[
\codeg_t(u) = \deg_t(u)+B_t\bff_{u;t}.
\]

Following Qin's notion of compatibly pointed elements
\cite{Qin24}, as recalled in \cite[Definition~3.6]{Cao23}, we say
that \(u\in\calA\) is \emph{compatibly pointed} if it is pointed
with respect to every \(t\in\bbT\) and its degree family forms a
tropical point of the cluster pattern. In this case, the
family
\[
[\bfg_u]\seteq\bigl(\deg_t(u)\bigr)_{t\in\bbT}
\]
is called the \emph{tropical point} of \(u\).  A compatibly pointed
element is called \emph{bipointed} if it is bipointed with respect
to every \(t\in\bbT\).

\begin{proposition}[{\cite{FZ07,GHKK18} (see also \cite[Proposition~3.7]{Cao23})}]
\label{prop: cluster monomials are bipointed}
Every cluster monomial in \(\calA\) is compatibly bipointed.
\end{proposition}

A compatibly pointed element \(u\in\calA\) is called \emph{good} if
it is \emph{universally positive}; that is, its Laurent expansion
with respect to every seed belongs to
\[
\Z_{\geq0}[x_{j;t}^{\pm1}\mid j\in J].
\]
In particular, every cluster monomial is good (see~\cite[Proposition~3.7]{Cao23}). 

For a nonzero polynomial $F(\bfy)=\sum_{\bfv}c_{\bfv}\bfy^{\bfv}$ and \(\bfr\in\Z^{\oplus J_{\ex}}\), set
\[
F[\bfr]\seteq \max\{\bfv^T\bfr\mid c_{\bfv}\ne0\}.
\]
If \(F(0)=1\), then \(F[\bfr]\geq0\).

Let \(u,u'\in\calA\) be good.  For \(t\in\bbT\), set
\[
\bfg_{u;t}\seteq\deg_t(u), \qquad S\bfg_{u;t} \seteq (s_jg_{j;u;t})_{j\in J_{\ex}}.
\]

We regard \(\calA\) as a subalgebra of the corresponding upper
cluster algebra, in which Cao's tropical and \(F\)-invariants are
defined.

\begin{definition}\label{def: tropical and F-invariants}
 {\cite[Definition~4.6]{Cao23}} For $t \in \bbT$, we define
\[
\langle u,u'\rangle_{\trop,t} \seteq
\bfg_{u;t}^TL_t\bfg_{u';t} + F_{u;t}[S\bfg_{u';t}]
\]
and
\[
(u\Vert u')_{F,t} \seteq
\langle u,u'\rangle_{\trop,t} + \langle u',u\rangle_{\trop,t}.
\]
\end{definition}

\begin{theorem}[{\cite[Theorem~4.10]{Cao23}}]
\label{thm: Cao F-invariant}
The quantities in {\rm Definition~\ref{def: tropical and F-invariants}} are independent of
\(t\).  Moreover,
\[
(u\Vert u')_F = F_{u;t}[S\bfg_{u';t}] + F_{u';t}[S\bfg_{u;t}] \in2\Z_{\geq0}.
\]
\end{theorem}

The evenness follows from \(s_j\in2\Z_{>0}\).  Following~\cite{CFM26}, we call
\begin{align} \label{eq: frako}
F_{u;t}[S\bfg_{u';t}]
\end{align}
the \emph{partial \(F\)-invariant} of \((u,u')\) with respect to
\(t\). Unlike \((u\Vert u')_F\), it may depend on \(t\) in general.

If \(k\in J_{\ex}\), then
\[
\bfg_{x_{k;t};t}=\mathbf e_k \qtq F_{x_{k;t};t}=1.
\]
Hence, for every good \(u\in\calA\),
\begin{equation}\label{eq: F-invariant with initial variable}
(x_{k;t}\Vert u)_F = F_{u;t}[S\mathbf e_k] = s_kf_{k;u;t}.
\end{equation}

\section{Folding of cluster algebras and tropical invariants}
\label{sec: folding cluster algebras}

In this section, we establish folding formulas for Cao's tropical
and \(F\)-invariants in cluster algebras associated with positive
braids.
Starting from an admissible folding of Cartan data, we compare
the compatible pairs and cluster variables attached to
corresponding \(\bfi\)-box seeds.
We then show that \(F\)-polynomials fold under the substitution
induced by the box fibers and derive folding formulas for both
orders of the tropical invariant between \(\bfi\)-box cluster
variables.  Their symmetrization gives the corresponding formula
for the \(F\)-invariant.
Neither the Cartan datum nor its simply-laced cover is assumed
to be of finite type.

\subsection{Admissible foldings and lifted sequences}
\label{subsec: admissible foldings and lifted sequences}

We follow the folding conventions of 
\cite[Sections~2.1--2.3]{JO26}, to which we refer for further details.

Let \(\hg\) be a simply-laced Kac--Moody algebra with index set
\(\hI\) and Cartan matrix
\[
\hC=(\hc_{\im,\jm})_{\im,\jm\in\hI}.
\]
Let \((\widehat\alpha_\im)_{\im\in\hI}\) be its simple roots, and set
\[
\hsfQ\seteq\soplus_{\im\in\hI}\Z\widehat\alpha_\im,
\qquad
\hsfQ_+\seteq \sum_{\im\in\hI}\Z_{\geq0}\widehat\alpha_\im,
\qquad
\hsfQ_-\seteq-\hsfQ_+.
\]

Let \(W_\hg\) and \(\ttB_\hg^+\) be the Weyl group and the positive
braid monoid of \(\hg\), with standard generators \(s_\im\) and
\(\sfr_\im\) $(\im \in \hI)$, respectively.

Let \(\widehat\Delta\) be the root system of \(\hg\), and set
\(\widehat\Delta_+\seteq\widehat\Delta\cap\hsfQ_+\). Let 
\[
\widehat\Delta_{\re} \seteq \bigl\{ \hw(\widehat\alpha_\im) \mid \hw\in W_\hg,\ \im\in\hI \bigr\}
\]
be the set of real roots, and set
\[
\widehat\Delta_{\re,+} \seteq \widehat\Delta_{\re}\cap\widehat\Delta_+, \qquad \widehat\Delta_{\re,-} \seteq -\widehat\Delta_{\re,+}.
\]

Let \(\Gamma\) be a finite group acting on \(\hI\) by diagram
automorphisms.  We assume that the action is \emph{admissible}, that
is,
\[
\hc_{\im,g(\im)}=0
\qquad
\bigl(\im\in\hI,\ g\in \Gamma,\ g(\im)\ne\im\bigr).
\]
Set
\[
I\seteq\hI/\Gamma, \qquad \pi:\hI\longrightarrow I, \qquad I_i\seteq\pi^{-1}(i), \qquad e_i\seteq|I_i|.
\]
Thus the vertices belonging to each \(I_i\) are pairwise
nonadjacent. To distinguish the two index sets, we generally use
\(\im,\jm\in\hI\) for the simply-laced datum and \(i,j\in I\) for
the quotient datum.

\smallskip

The quotient Cartan matrix \(C=(c_{i,j})_{i,j\in I}\) is defined by
\[
c_{i,j} \seteq \sum_{\im\in I_i}\hc_{\im,\jm}, \qquad \jm\in I_j.
\]
The right-hand side is independent of the choice of
\(\jm\in I_j\).  Let
\[
D=\diag(d_i\mid i\in I)
\]
be the normalized symmetrizer of \(C\).  Then there exists
\(m\in\Z_{>0}\) such that 
\begin{equation}\label{eq: symmetrizer and fiber cardinality}
d_ie_i=m \qquad(i\in I).
\end{equation}
We denote by \(\g\) the symmetrizable Kac--Moody algebra associated
with \(C\).
Let \((\alpha_i)_{i\in I}\) be its simple roots, and set
\[
\sfQ\seteq\bigoplus_{i\in I}\Z\alpha_i, \qquad
\sfQ_+\seteq\sum_{i\in I}\Z_{\geq0}\alpha_i, \qquad
\sfQ_-\seteq-\sfQ_+.
\]
We use the same symbol \(\pi\) for the induced map of root lattices
\[
\pi:\hsfQ\longrightarrow\sfQ, \qquad
\widehat\alpha_\im\longmapsto\alpha_{\pi(\im)}.
\]

When \(\hg\) is of finite type, the nontrivial admissible foldings
used below are summarized as follows:
\begin{align} \label{eq: finite folding}
\begin{array}{c|c|c}
\hg & \Gamma & \g=\hg/\Gamma\\
\hline
A_{2n-1} & \Z/2\Z & C_n\\
D_{n+1}  & \Z/2\Z & B_n\\
E_6      & \Z/2\Z & F_4\\
D_4      & \Z/3\Z & G_2
\end{array}    
\end{align}
Here \(\Gamma\) acts by the standard diagram automorphisms.

We equip \(\sfQ\) with the symmetric bilinear form
\[
(\alpha_i,\alpha_j)_\g
\seteq
d_i c_{i,j}
\qquad(i,j\in I).
\]

Let \(W_\g\) and \(\ttB_\g^+\) be the Weyl group and the positive
braid monoid of \(\g\), with standard generators \(s_i\) and
\(\sfr_i\), respectively.
For \(i\in I\), set
\[
\hs_i \seteq \prod_{\im\in I_i}s_\im \in W_\hg,
\qquad
\hsfr_i \seteq \prod_{\im\in I_i}\sfr_\im \in \ttB^+_\hg.
\]
The products are independent of the ordering of \(I_i\) by
admissibility.  The assignments
\[
s_i\longmapsto\hs_i,
\qquad
\sfr_i\longmapsto\hsfr_i
\]
induce embeddings
\[
W_\g\lhook\joinrel\longrightarrow W_\hg,
\qquad
\ttB_\g^+\lhook\joinrel\longrightarrow\ttB_\hg^+.
\]
See~\cite[Theorem~28]{Castella08} and~\cite[Corollary~2.5]{JO26}.

Now let $\ii=(i_1,\ldots,i_l)\in I^l$ be an arbitrary finite sequence.  
Fix an ordering
\[
I_i=(\im_{i,1},\ldots,\im_{i,e_i})
\qquad(i\in I),
\]
and define the blockwise lift of \(\ii\) by
\begin{align}\label{eq: seq lift}
\hii \seteq I_{i_1}\star\cdots\star I_{i_l} = (\im_{i_1,1},\ldots,\im_{i_1,e_{i_1}},  \ldots,  \im_{i_l,1},\ldots,\im_{i_l,e_{i_l}}).    
\end{align}
Here \(\star\) denotes concatenation of the ordered blocks.
Set
\[
N_0\seteq0,
\qquad
N_k\seteq\sum_{a=1}^ke_{i_a},
\qquad
\kappa(k,r)\seteq N_{k-1}+r
\]
for \(1\leq k\leq l\) and \(1\leq r\leq e_{i_k}\).   The integer \(\kappa(k,r)\) is the position in \(\hii\) of the
\(r\)-th entry of the \(k\)-th block \(I_{i_k}\).

Writing
$\sfr_\ii \seteq \sfr_{i_1}\cdots\sfr_{i_l}$, the lifted sequence represents
$\sfr_\hii = \hsfr_{i_1}\cdots \hsfr_{i_l}.$ 
Consequently, if $\ttb=\sfr_\ii\in\ttB_\g^+,$ we denote its canonical lift by $\httb \seteq \sfr_\hii\in\ttB_\hg^+.$
Moreover, \(\ii\) is reduced if and only if \(\hii\) is reduced.
For \(\ttb\in\ttB_\g^+\), we write \(\Seq(\ttb)\) for the set of
finite sequences representing \(\ttb\).

For later use, define the root sequence of \(\ii\) by
\begin{equation}\label{eq: root sequence}
\be_k^\ii \seteq s_{i_1}\cdots s_{i_{k-1}}(\alpha_{i_k}) \qquad(1\leq k\leq l).
\end{equation}
In general, the real roots \(\be_k^\ii\) need not be positive or
pairwise distinct. We also write
\[
\hbe_{k,r} \seteq \be_{\kappa(k,r)}^\hii  \qquad (1\leq r\leq e_{i_k}).
\]
The folding map on root lattices satisfies
\[
\pi(\hbe_{k,r})=\be_k^\ii
\]
(\cite[Lemma~2.4]{JO26}).
In particular, the block
\[
\{\hbe_{k,r}\mid1\leq r\leq e_{i_k}\}
\]
is the fiber of the root \(\be_k^\ii\).  These blockwise fibers will
induce the fibers of \(\bfi\)-boxes and of the associated cluster
variables in the following subsections.

\subsection{\texorpdfstring{\(\bfi\)}{i}-box seeds and folding of
compatible pairs}
\label{subsec: folding of box seeds}

Keep the notation of
Subsection~\ref{subsec: admissible foldings and lifted sequences}.
We refer to \cite{KKOP23F,KK24} for further details on
\(\bfi\)-boxes and their associated seeds.

For \(1\leq p\leq l\), set
\[
p^+ \seteq \min\bigl(\{q>p\mid i_q=i_p\}\cup\{+\infty\}\bigr)
\qtq
p^- \seteq \max\bigl(\{q<p\mid i_q=i_p\}\cup\{-\infty\}\bigr).
\]
More generally, for \(j\in I\), set
\[
p(j)^+ \seteq \min\bigl(\{q\geq p\mid i_q=j\}\cup\{+\infty\}\bigr)
\qtq
p(j)^- \seteq \max\bigl(\{q\leq p\mid i_q=j\}\cup\{-\infty\}\bigr).
\]

A \emph{\(\bfi\)-box associated with \(\ii\)} is an interval
\[
\frakc=[x,y]\subset[1,l] \qquad\text{such that}\qquad i_x=i_y.
\]
We write
\[
i_\frakc\seteq i_x=i_y
\]
for its \emph{color}, and abbreviate a singleton \(\bfi\)-box
\([x,x]\) to \([x]\).

Two \(\bfi\)-boxes $\frakc=[x,y]$ and $\frakc'=[x',y']$
associated with \(\ii\) are said to \emph{commute} if
\begin{align} \label{eq: commuting boxes}
x^-<x'\leq y'<y^+ \qtoq (x')^-<x\leq y<(y')^+.
\end{align}
A collection of pairwise commuting \(\bfi\)-boxes associated with
\(\ii\) is called a \emph{commuting family}.  Such a family in
\([1,l]\) is called \emph{maximal} if it is maximal with respect to
inclusion.   Throughout this section, nonempty admissible chains have final
envelope \([1,l]\), called their \emph{extent}.
For the empty sequence, the chain and its underlying family are empty,
with extent \(\emptyset\). 

For an interval \([x,y]\), set
\[
[x,y\}\seteq[x,y(i_x)^-] \qtq \{x,y]\seteq[x(i_y)^+,y]
\]
whenever these intervals are nonempty.  An
\emph{admissible chain of \(\bfi\)-boxes associated with \(\ii\)}
of \emph{extent} \([1,l]\) is a sequence
\[
\frakC=(\frakc_k=[x_k,y_k])_{1\leq k\leq l}
\]
such that
\bnum
\item
$\tfrakc_k \seteq \bigcup_{1\leq s\leq k}\frakc_s = [\tx_k,\ty_k]$
is an interval of cardinality \(k\), called the \emph{\(k\)-th envelope}, and
\item
$\frakc_k = [\tx_k,\ty_k\}$ or $\frakc_k = \{\tx_k,\ty_k]$.
\ee
Its underlying family is
\[
\calF(\frakC)\seteq\{\frakc_k\mid1\leq k\leq l\}.
\]
It is a maximal commuting family of \(\bfi\)-boxes associated with
\(\ii\), and every maximal commuting family of \(\bfi\)-boxes arises in this way.  Note that a maximal commuting family of \(\bfi\)-boxes need not
determine a unique admissible chain. 

For \(l\ge1\), we describe box moves on admissible chains as follows. Since the first
envelope $\tfrakc_1$ is a singleton, we have
\[
\tfrakc_1=\frakc_1=[x].
\]
For \(1\leq k<l\), let \(T_k\in\{\calL,\calR\}\) record whether
\(\tfrakc_{k+1}\) is obtained from
\(\tfrakc_k\) by adjoining its left or right neighboring
point, respectively.  Thus the chain \(\frakC\) has the associated
pair
\[
(x,\frakT),
\qquad
\frakT=(T_1,\ldots,T_{l-1})\in\{\calL,\calR\}^{l-1}.
\]
Conversely, its extent \([1,l]\) together with \(\frakT\) determines
\(x\), all the envelopes, and hence \(\frakC\).

For \(1\leq k<l\), the member \(\frakc_k\) is called
\emph{movable} if
\[
k=1
\quad\text{or}\quad
T_{k-1}\ne T_k.
\]
For a movable \(\frakc_k\), the \emph{\(k\)-th box move}
\(\frakB_k(\frakC)\) is the admissible chain with associated pair
\((x',\frakT')\) defined as follows.  If \(k=1\), then
\[
(x',T'_1)
=
\begin{cases}
(x+1,\calL),&T_1=\calR,\\
(x-1,\calR),&T_1=\calL,
\end{cases}
\qquad
T'_s=T_s\quad(s\ne1).
\]
If \(k>1\), then \(x'=x\),
\[
T'_{k-1}=T_k,
\qquad
T'_k=T_{k-1},
\qquad
T'_s=T_s\quad(s\notin\{k-1,k\}).
\]
Equivalently, the move shifts the envelope \(\tfrakc_k\)
by one position inside \(\tfrakc_{k+1}\). 

\begin{lemma}[{\cite{KKOP24A,KK24}}]
\label{lem: two types of box moves}
Let \(\frakc_k\) be movable in \(\frakC\), and set
\[
\frakC' \seteq \frakB_k(\frakC) = (\frakc'_s)_{1\leq s\leq l}.
\]
Then exactly one of the following holds.
\bnum
\item Suppose that \(\tfrakc_{k+1}\) is not a
\(\bfi\)-box.  Then the box move interchanges the two consecutive
members \(\frakc_k\) and \(\frakc_{k+1}\); namely,
\[
\frakc'_k=\frakc_{k+1}, \qquad \frakc'_{k+1}=\frakc_k, \qquad \frakc'_s=\frakc_s \quad(s\notin\{k,k+1\}).
\]
Consequently,
\[
\calF(\frakC')=\calF(\frakC).
\]
This is called a \emph{permutation-type box move}.

\item Suppose that
\[
\tfrakc_{k+1}=[u,v]
\]
is an \(\bfi\)-box.  Then the box move replaces exactly one member of
the underlying family.  More precisely, there exist unique boxes
\[
\frakc\in\calF(\frakC) \qtq \frakc^\star\in\calF(\frakC')
\]
such that
\[
\{\frakc,\frakc^\star\} = \{[u^+,v],[u,v^-]\}
\qtq
\calF(\frakC') = \left( \calF(\frakC)\setminus\{\frakc\} \right) \sqcup\{\frakc^\star\}.
\]
This is called a \emph{mutation-type box move}.
\ee
\end{lemma}

Any two admissible chains of extent \([1,l]\) are connected
by a finite sequence of box moves (\cite[Lemma~5.10]{KKOP24A} and
\cite[Proposition~2.6]{KK24}).

We next recall the compatible pair attached to \(\calF(\frakC)\).
Let
\[
\frakE^\ii \seteq \soplus_{p=1}^l\Z\ve_p,
\]
and define a skew-symmetric pairing on \(\frakE^\ii\) by
\[
\lambda^\ii(\ve_p,\ve_q) \seteq 
\begin{cases}
-(\be_p^\ii,\be_q^\ii)_\g,&p<q,\\
0,&p=q,\\
(\be_p^\ii,\be_q^\ii)_\g,&p>q,
\end{cases}
\]
This is the negative of the pairing in~\cite[(5.4)]{KKOP23F}.
When \(\ii\) is reduced, the matrix \(L^\calF\) defined below is the
negative of the \(\Lambda\)-invariant matrix of the determinantial
modules indexed by \(\calF\); see Proposition~\ref{prop: Laurent maximum realization of box seed} below.

For an \(\bfi\)-box \(\frakc=[x,y]\) associated with \(\ii\), set
\[
\calE_\frakc
\seteq
\sum_{\substack{x\leq p\leq y\\i_p=i_\frakc}}\ve_p.
\]
We call $\calE_\frakc$ the \emph{box vector associated with $\frakc$}.  If \(\calF\) is a maximal commuting family, define
\[
L^\calF \seteq \bigl(l_{\frakc,\frakc'}^\calF\bigr)_{\frakc,\frakc'\in\calF},
\qquad
l_{\frakc,\frakc'}^\calF \seteq \lambda^\ii(\calE_\frakc,\calE_{\frakc'}).
\]

An \(\bfi\)-box \(\frakc=[x,y]\in\calF\) is called
\emph{frozen} if
\[
x^-=-\infty \qtq y^+=+\infty.
\]
Let \(\calF_{\ex}\) be the set of the remaining \(\bfi\)-boxes.
Following \cite[Section~3.3]{JO26}, we use the exchange matrix
\[
B^\calF = \bigl(b_{\frakc,\frakc'}^\calF\bigr)_{\frakc\in\calF,\,\frakc'\in\calF_{\ex}}
\seteq -B_{\mathrm{KK}}(\calF),
\]
where \(B_{\mathrm{KK}}(\calF)\) is the restriction to the
exchangeable columns of the matrix in \cite[(3.1)]{KK24}. By \cite[Theorem~3.5]{JO26}, 
\[
\text{the pair \((L^\calF,B^\calF)\) is compatible}
\]
and satisfies
\begin{equation}\label{eq: compatibility of box seed}
(B^\calF)^TL^\calF = (2D^\calF\mid0),  \qquad D^\calF \seteq \diag(d_{i_\frakc}\mid\frakc\in\calF_{\ex}).
\end{equation}
We denote the obtained \(\Uplambda\)-seed by
\[
\Sigma^\ii(\calF) \seteq (\bfx^\calF,L^\calF,B^\calF),
\qquad
\bfx^\calF = (x_\frakc\mid\frakc\in\calF).
\]

We now compare this seed with the corresponding lifted seed.  Choose an admissible chain \(\frakC\) realizing \(\calF\).
Let \(\hfrakC\) be its canonical lift
\cite[Section~3.4]{JO26}, and set
\[
\hcalF\seteq\calF(\hfrakC).
\] 
Thus \(\hcalF\) is a
maximal commuting family of \(\bfi\)-boxes associated with
\(\hii\).  By \cite[Section~2.2]{KK24} and \cite[Section~3.4]{JO26},
we have a surjection 
\[
\pi_{\bx}^{\frakC}: \hcalF\longrightarrow\calF. 
\]
For \(\frakc\in\calF\), define its \emph{box fiber} by
\[
\Fib_{\frakC}(\frakc) \seteq (\pi_{\bx}^{\frakC})^{-1}(\frakc). 
\]
Then \cite[Section~3.4]{JO26} implies
\begin{equation}\label{eq: cardinality of box fiber}
|\Fib_{\frakC}(\frakc)|=e_{i_\frakc}.
\end{equation}

Let $\hSigma^\hii(\hcalF) \seteq (\hbfx^{\hcalF},  \hL^{\hcalF},  \widehat{B}^{\hcalF})$ 
be the \(\Uplambda\)-seed associated with the lifted family.  The
exchange and compatibility matrices fold simultaneously as follows.

\begin{theorem}[{\cite[Theorem~3.15]{JO26}}]
\label{thm: folding of box compatible pairs}
Let \(\frakc,\frakc'\in\calF\), and choose
\(\hfrakc'\in\Fib_{\frakC}(\frakc')\).
Then
\[
b_{\frakc,\frakc'}^\calF = \sum_{\hfrake\in\Fib_{\frakC}(\frakc)} \hb_{\hfrake,\hfrakc'}^{\hcalF}
\quad(\frakc'\in\calF_{\ex}),
\qquad
l_{\frakc,\frakc'}^\calF = d_{i_\frakc} \sum_{\hfrake\in\Fib_{\frakC}(\frakc)} \hl_{\hfrake,\hfrakc'}^{\hcalF}.
\]
Both right-hand sides are independent of the choice of
\(\hfrakc'\).  In particular,
\[
(\widehat{B}^{\hcalF})^T \hL^{\hcalF} = (2\operatorname{Id}\mid0) 
\qtq
(B^\calF)^TL^\calF = (2D^\calF\mid0).
\]
\end{theorem}

The construction and
Theorem~\ref{thm: folding of box compatible pairs} require neither
reducedness of \(\ii\) nor finite type of the Cartan datum.  Moreover,
all maximal commuting families of \(\bfi\)-boxes associated with
\(\ii\) are connected by box moves.   A
permutation-type box move only relabels the corresponding seed,
whereas a mutation-type box move acts by simultaneous mutation of
\(B^\calF\) and \(L^\calF\).
Consequently, for any
two such families \(\calF\) and \(\calG\),
\[
\Sigma^\ii(\calF) \quad\text{and}\quad \Sigma^\ii(\calG)
\]
are mutation-equivalent up to seed relabeling.  In particular, 
the $\Uplambda$-cluster algebra
\[
\calA^\ii[1,l]
\seteq
\calA\bigl(\Sigma^\ii(\calF)\bigr),
\]
where \(\calF\) is any maximal commuting family of \(\bfi\)-boxes
associated with \(\ii\) and having extent \([1,l]\), is independent
of the choice of \(\calF\).

\subsection{Folding of cluster variables and \texorpdfstring{\(F\)}{F}-polynomials}

Keep the notation of
Subsection~\ref{subsec: folding of box seeds}. The box fibers induce a surjective homomorphism between the
corresponding seed Laurent rings:
\begin{equation}\label{eq: box Laurent folding map}
\begin{aligned}
\phi_{\frakC}: \calL\bigl(\hSigma^\hii(\hcalF)\bigr) &\longrightarrow \calL\bigl(\Sigma^\ii(\calF)\bigr),\\
\hx_{\hfrakc} &\longmapsto x_\frakc \qquad\qquad \bigl(\hfrakc\in\Fib_{\frakC}(\frakc)\bigr).
\end{aligned}
\end{equation}
They also define the folding map on exponent lattices
\[
\Fold_{\frakC}:
\Z^{\hcalF}\longrightarrow\Z^\calF,
\qquad
\bigl(\Fold_{\frakC}(\hbfg)\bigr)_\frakc
\seteq
\sum_{\hfrakc\in\Fib_{\frakC}(\frakc)}
\hg_{\hfrakc}.
\]
We use the same notation for its restriction to the exchangeable
coordinates.

For \(\frake\in\calF_{\ex}\) and
\(\hfrake\in\Fib_{\frakC}(\frake)\), define the corresponding
 Laurent monomials by 
\[
\dy_\frake
\seteq
\prod_{\frakc\in\calF}
x_\frakc^{b_{\frakc,\frake}^\calF}
\qtq
\hdy_{\hfrake}
\seteq
\prod_{\hfrakc\in\hcalF}
\hx_{\hfrakc}^
{\hb_{\hfrakc,\hfrake}^{\hcalF}}.
\]
Set  
\[
\dbfy\seteq(\dy_\frake)_{\frake\in\calF_{\ex}}
\qtq
\hdbfy\seteq(\hdy_{\hfrake})_{\hfrake\in\hcalF_{\ex}}.
\]
The folding of exchange matrices in
Theorem~\ref{thm: folding of box compatible pairs} gives
\begin{equation}\label{eq: folding of column monomials}
\phi_{\frakC}(\hdy_{\hfrake})
=
\dy_\frake.
\end{equation} 
Let
\[
\bfy=(y_\frake)_{\frake\in\calF_{\ex}}
\qtq
\widehat{\bfy}=(\hy_{\hfrake})_{\hfrake\in\hcalF_{\ex}}
\]
be families of indeterminates.
Define the algebra homomorphism
\[
\phi_{\frakC,y}:
\C[\hy_{\hfrake}\mid\hfrake\in\hcalF_{\ex}]
\longrightarrow
\C[y_\frake\mid\frake\in\calF_{\ex}]
\qquad 
\text{by}
\quad 
\phi_{\frakC,y}(\hy_{\hfrake})=y_\frake
\qquad
\bigl(\hfrake\in\Fib_{\frakC}(\frake)\bigr).
\]

The folding relation defined by \(\phi_{\frakC}\) is preserved at
the canonically lifted initial and final seeds for braid moves
and box moves.
Consequently, it folds all corresponding cluster variables in box seeds.

\begin{theorem}[{\cite[Theorem~6.2 and Corollary~6.3]{JO26}}]
\label{thm: folding of box cluster variables}
Let \(\ii,\jj\in\Seq(\ttb)\) and 
\bnum
\item let \(\frakD\) be an admissible chain of \(\bfi\)-boxes
associated with \(\jj\), with underlying maximal commuting family
\(\calG\seteq\calF(\frakD)\),
\item  let \(\hfrakD\) be the canonical lift of \(\frakD\) associated
with \(\hjj\), and set \(\hcalG=\calF(\hfrakD)\).  
\ee 
For \(\frake\in\calG\) and \(\hfrake\in\Fib_{\frakD}(\frake)\), we have
\[
\phi_{\frakC} \bigl(\hx_{\hfrake}^{\hcalG}\bigr) = x_\frake^\calG
\qtq
\deg_{\Sigma^\ii(\calF)} \bigl(x_\frake^\calG\bigr) = \Fold_{\frakC}\left( \deg_{\hSigma^\hii(\hcalF)} \bigl(\hx_{\hfrake}^{\hcalG}\bigr) \right).
\]
The same statements hold for the corresponding cluster monomials.
\end{theorem}

We now record the consequence for \(F\)-polynomials.  Let 
\[
\hz = \hbfx^{\hbfg}\, \hF_{\hz}(\hdbfy)
\]
be a cluster monomial associated with
\(\hSigma^\hjj(\hcalG)\), and set 
\[
z\seteq\phi_{\frakC}(\hz).
\]
By Theorem~\ref{thm: folding of box cluster variables}, we have
\[
\bfg_z=\Fold_{\frakC}(\hbfg),
\]
where both degrees are taken with respect to the fixed seeds
\(\Sigma^\ii(\calF)\) and \(\hSigma^\hii(\hcalF)\), respectively.
Applying \(\phi_{\frakC}\) to the pointed expression of
\(\hz\) and using
\eqref{eq: folding of column monomials}, we obtain
\[
z = \bfx^{\Fold_{\frakC}(\hbfg)} \left. \phi_{\frakC,y}\bigl( \hF_{\hz}(\widehat{\bfy}) \bigr)  \right|_{\bfy=\dbfy}.
\]
Comparing this with the downstairs pointed expression
\[
z=\bfx^{\bfg_z}F_z(\dbfy),
\]
the uniqueness of the pointed expression gives the following formal
polynomial identity.

\begin{proposition}\label{prop: folding of pointed F-parts}
With the notation above, we have 
\[
F_z(\bfy) = \phi_{\frakC,y}\bigl( \hF_{\hz}(\widehat{\bfy}) \bigr). 
\]
\end{proposition}

\subsection{Global folding and tropical invariants for positive braids}
\label{subsec: global folding for positive braids}

Let \(\ttb\in\ttB_\g^+\), and choose
\[
\ii=(i_1,\ldots,i_l)\in\Seq(\ttb).
\]

Let \(\calF\) be a maximal commuting family of \(\bfi\)-boxes
associated with \(\ii\), and choose an admissible chain \(\frakC\)
realizing \(\calF\).  We define
\[
\calA(\ttb) \seteq \calA\bigl(\Sigma^\ii(\calF)\bigr).
\]
Its mutation class depends only on \(\ttb\), and not on the choices of
\(\ii\), \(\frakC\), or \(\calF\).
Similarly, the blockwise lift \(\hii\) in~\eqref{eq: seq lift} represents \(\httb\in\ttB_\hg^+\), and we set
\[
\calA(\httb) \seteq \calA\bigl(\hSigma^\hii(\hcalF)\bigr).
\]

Keep the Laurent homomorphism \(\phi_{\frakC}\) defined in
\eqref{eq: box Laurent folding map}.  Although this map is defined
using the initial box fibers, Theorem~\ref{thm: folding of box cluster variables}
shows that it simultaneously folds all corresponding box seeds in
the mutation classes.  More precisely, let
\(\jj=(j_1,\ldots,j_l)\in\Seq(\ttb)\), let \(\calG\) be a maximal
commuting family of \(\bfi\)-boxes associated with \(\jj\), and
choose an admissible chain \(\frakD\) realizing \(\calG\).    Let \(\hjj\) be the blockwise lift of \(\jj\), let \(\hfrakD\)
be the canonical lift of \(\frakD\), and set
\(\hcalG=\calF(\hfrakD)\).  Then, for every
\(\frake\in\calG\) and
\(\hfrake\in\Fib_{\frakD}(\frake)\),
\begin{equation}\label{eq: global variable folding}
\phi_{\frakC}
\bigl(
\hx_{\hfrake}^{\hcalG}
\bigr)
=
x_{\frake}^{\calG}.
\end{equation}

Combining Theorem~\ref{thm: folding of box cluster variables} with
Theorem~\ref{thm: folding of box compatible pairs}, we obtain the
simultaneous folding of the entire lifted box seed:
\begin{equation}\label{eq: global box seed folding}
\hSigma^\hjj(\hcalG)
\longmapsto
\Sigma^\jj(\calG).
\end{equation}
 Here \eqref{eq: global box seed folding} combines the variable folding
\eqref{eq: global variable folding} with the folding of compatible pairs.

In particular, for
$\hcalM = \prod_{\hfrake\in\hcalG} \bigl(\hx_{\hfrake}^{\hcalG}\bigr)^{a_{\hfrake}}$,
we have
\[
\phi_{\frakC}(\hcalM) = \prod_{\frake\in\calG} \bigl(x_\frake^\calG\bigr)^{ \sum_{\hfrake\in\Fib_{\frakD}(\frake)} a_{\hfrake}}.
\]
Thus every cluster monomial associated with a lifted box seed is
sent to the corresponding downstairs cluster monomial.

The fixed Laurent homomorphism restricts to the cluster algebras as
follows.

\begin{theorem}[{\cite[Theorems~6.2 and~6.4]{JO26}}]
\label{thm: global folding for positive braids}
For every \(\ttb\in\ttB_\g^+\), the homomorphism
\(\phi_{\frakC}\) restricts to a surjective \(\C\)-algebra
homomorphism
\[
\phi_{\frakC}^{\ttb}:
\calA(\httb)\twoheadrightarrow\calA(\ttb).
\]
For every \(\jj\in\Seq(\ttb)\) and every maximal commuting family
\(\calG\) of \(\bfi\)-boxes associated with \(\jj\), it maps the
cluster variables and cluster monomials of
\(\hSigma^\hjj(\hcalG)\) to those of
\(\Sigma^\jj(\calG)\) according to their box fibers.
\end{theorem}

Theorem~\ref{thm: global folding for positive braids} requires
neither reducedness of the sequences nor finite type of the Cartan
data.  Every symmetrizable Cartan matrix admits an admissible
simply-laced realization \cite[Remark~1.3]{JO26}. Combining
Theorem~\ref{thm: global folding for positive braids}
with Proposition~\ref{prop: folding of pointed F-parts}, we compute the
\(F\)-polynomial and the partial \(F\)-invariant of every cluster variable
arising from a box seed by the corresponding fiber
evaluation upstairs. 

\begin{theorem}\label{thm: folding of F-invariants between box variables}
Let $\ii,\jj\in\Seq(\ttb)$,
and let \(\frakC\) and \(\frakD\) be admissible chains associated
with \(\ii\) and \(\jj\), respectively.  For $\frakc\in\calF(\frakC),\frakc'\in\calF(\frakD)$, and every
$\hfrakc'\in\Fib_{\frakD}(\frakc')$,
set
\[
\hx_{\Fib_{\frakC}(\frakc)}
\seteq
\prod_{\hfrakc\in\Fib_{\frakC}(\frakc)}\hx_{\hfrakc}.
\]
Then the two ordered tropical pairings fold separately:
\begin{align}
\label{eq: folding of first ordered tropical pairing}
\langle x_{\frakc},x_{\frakc'}\rangle_{\trop}
&=d_{i_{\frakc}}
\left\langle
\hx_{\Fib_{\frakC}(\frakc)},\hx_{\hfrakc'}
\right\rangle_{\trop},\\
\label{eq: folding of second ordered tropical pairing}
\langle x_{\frakc'},x_{\frakc}\rangle_{\trop}
&=d_{i_{\frakc}}
\left\langle
\hx_{\hfrakc'},\hx_{\Fib_{\frakC}(\frakc)}
\right\rangle_{\trop}.
\end{align}
Both right-hand sides are independent of the choice of
\(\hfrakc'\in\Fib_{\frakD}(\frakc')\).  Adding the two identities
gives
\begin{equation}\label{eq: folding of F-invariants between box variables}
\bigl(x_\frakc \Vert x_{\frakc'} \bigr)_F = d_{i_\frakc} \bigl( \hx_{\Fib_{\frakC}(\frakc)} \Vert \hx_{\hfrakc'} \bigr)_F,
\end{equation}
\end{theorem}
\begin{proof}
We compute the tropical pairings in the box seeds associated with
\(\frakC\) and its lift.  By
Theorem~\ref{thm: global folding for positive braids},
\[
\phi_{\frakC}^{\ttb}(\hx_{\hfrakc'})=x_{\frakc'}
\qquad
\bigl(\hfrakc'\in\Fib_{\frakD}(\frakc')\bigr).
\]
Write \(\bfg\) and \(\widehat{\bfg}\) for the degrees of
\(x_{\frakc'}\) and \(\hx_{\hfrakc'}\) in these seeds, and set
\[
\widehat{\ve}_{\frakc}
\seteq
\sum_{\hfrakc\in\Fib_{\frakC}(\frakc)}\ve_{\hfrakc}.
\]
The folding of cluster variables gives
\(\bfg=\Fold_{\frakC}(\widehat{\bfg})\).
Together with
Theorem~\ref{thm: folding of box compatible pairs}, this yields
\[
\ve_{\frakc}^{T}L\bfg
=d_{i_{\frakc}}
\widehat{\ve}_{\frakc}^{T}\widehat L\widehat{\bfg}.
\]

Since \(x_{\frakc}\) and
\(\hx_{\Fib_{\frakC}(\frakc)}\) are initial cluster monomials,
their degrees are \(\ve_{\frakc}\) and
\(\widehat{\ve}_{\frakc}\), respectively, and their
\(F\)-polynomials are \(1\).  Hence
\[
\begin{aligned}
\langle x_{\frakc},x_{\frakc'}\rangle_{\trop}
=\ve_{\frakc}^{T}L\bfg \qtq 
\left\langle
\hx_{\Fib_{\frakC}(\frakc)},\hx_{\hfrakc'}
\right\rangle_{\trop}
=\widehat{\ve}_{\frakc}^{T}\widehat L\widehat{\bfg}.
\end{aligned}
\]
This proves \eqref{eq: folding of first ordered tropical pairing}.

For the reverse order, let \(S\) and \(\widehat S\) denote the
maps that multiply the exchangeable coordinates by the
compatibility types \(2D^{\calF}\) and \(2\operatorname{Id}\),
respectively.  By definition and skew-symmetry,
\[
\begin{aligned}
\langle x_{\frakc'},x_{\frakc}\rangle_{\trop}
=-\ve_{\frakc}^{T}L\bfg
  +F_{x_{\frakc'}}[S\ve_{\frakc}] \qtq 
\left\langle
\hx_{\hfrakc'},\hx_{\Fib_{\frakC}(\frakc)}
\right\rangle_{\trop}
=-\widehat{\ve}_{\frakc}^{T}\widehat L\widehat{\bfg}
  +\widehat F_{\hx_{\hfrakc'}}
     [\widehat S\widehat{\ve}_{\frakc}].
\end{aligned}
\]
It remains to compare the partial \(F\)-invariants. 

If \(\frakc\) is exchangeable, then
\[
S\ve_{\frakc}=2d_{i_{\frakc}}\ve_{\frakc},
\qquad
\widehat S\widehat{\ve}_{\frakc}
=2\widehat{\ve}_{\frakc}.
\]
The coefficients of \(\widehat F_{\hx_{\hfrakc'}}\) are nonnegative by
universal positivity, so identifying variables along box fibers
cannot cancel monomials.  By Proposition~\ref{prop: folding of pointed F-parts}, we can take the
positive factor \(d_{i_{\frakc}}\) outside the maximum to obtain 
\[
F_{x_{\frakc'}}[2d_{i_{\frakc}}\ve_{\frakc}]
=d_{i_{\frakc}}\widehat F_{\hx_{\hfrakc'}}
  [2\widehat{\ve}_{\frakc}].
\]
If \(\frakc\) is frozen, every box in
\(\Fib_{\frakC}(\frakc)\) is frozen by
Theorem~\ref{thm: folding of box compatible pairs}. Thus
\[
S\ve_{\frakc}=0,
\qquad
\widehat S\widehat{\ve}_{\frakc}=0,
\] 
so both partial \(F\)-invariants vanish. This proves
\eqref{eq: folding of second ordered tropical pairing}.

The left-hand sides of the two identities are independent of
the choice of \(\hfrakc'\), so their right-hand sides are as well.
Adding the identities gives
\eqref{eq: folding of F-invariants between box variables}.
\end{proof}

\begin{example}\label{ex: A5 C3 folding of F-invariants}
Consider the folding \(A_5\to C_3\) with
\[
I_1=\{1,5\},\qquad I_2=\{2,4\},\qquad I_3=\{3\},
\qquad D^{C_3}=\diag(1,1,2).
\]
For the height function \(\xi=(1,2,3)\), with the convention of
Subsection~\ref{subsec: height functions and Q-adapted sequences},
the reduced \(Q\)-adapted sequence and its blockwise lift are 
\[
\ii=(1,2,1,3,2,1,3,2,3)\in\Seq(\ttb),
\qquad
\ttb=\sfr_\ii,
\qquad
\hii=(1,5,2,4,1,5,3,2,4,1,5,3,2,4,3).
\]
Let \(\frakC\) be the all-\(\calR\) chain and let \(\frakD\) be the
chain with associated sequence
\(\calR\calR\calR\calR\calL\calL\calL\calL\).  Their underlying
families are
\[
\begin{aligned}
\calF(\frakC)
&=
([1],[2],[1,3],[4],[2,5],[1,6],[4,7],[2,8],[4,9]),\\
\calF(\frakD)
&=
([5],[6],[7],[5,8],[7,9],[4,9],[3,6],[2,8],[1,6]).
\end{aligned}
\]
Set $\frakc=[2,5]$ and $\frakc'=[7,9]$.
Then \(i_\frakc=2\) and \(i_{\frakc'}=3\), and both
\(\frakc\) and \(\frakc'\) are non-singleton boxes.  The canonical
lifted chains give
\[
\Fib_{\frakC}([2,5]) = \{[3,8],[4,9]\},
\qquad
\Fib_{\frakD}([7,9]) = \{[12,15]\}.
\]
In particular,
\[
\phi_{\frakC}^{\ttb}(\hx_{[12,15]}) = x_{[7,9]},
\qquad
\hx_{\Fib_{\frakC}([2,5])} = \hx_{[3,8]}\hx_{[4,9]}.
\]

We record a direct cluster calculation.  With respect to the
all-\(\calR\) seed, whose exchangeable boxes are ordered as
\[
([1],[2],[1,3],[4],[2,5],[4,7]),
\]
the box-move mutation from \(\frakC\) to \(\frakD\) gives
\[
\bff_{x_{[7,9]};\frakC} = (2,2,2,1,2,1) 
\qtqh
f_{[2,5];x_{[7,9]}}=2.
\]
The lifted mutation calculation gives
\[
\hF_{\hx_{[12,15]}}
\bigl[\mathbf e_{[3,8]}+\mathbf e_{[4,9]}\bigr]=2.
\]
The compatibility types have diagonal entries \(2d_{i_\frake}\) at exchangeable boxes \(\frake\) downstairs and \(2\) upstairs.  Since \(d_2=1\), we obtain
\[
\begin{aligned}
\bigl(x_{[2,5]}\Vert x_{[7,9]}\bigr)_F
=2d_2 f_{[2,5];x_{[7,9]}}
=4,\  \
\bigl(\hx_{[3,8]}\hx_{[4,9]}
\Vert\hx_{[12,15]}\bigr)_F
=\hF_{\hx_{[12,15]}}
  \bigl[2(\mathbf e_{[3,8]}+\mathbf e_{[4,9]})\bigr]
=4.
\end{aligned}
\]
Thus
\[
\bigl(x_{[2,5]}\Vert x_{[7,9]}\bigr)_F
=d_2\bigl(\hx_{[3,8]}\hx_{[4,9]}
\Vert\hx_{[12,15]}\bigr)_F.
\]
\end{example}

\section{Categorical folding of
\texorpdfstring{\(\Lambda\)}{Lambda}- and
\texorpdfstring{\(\de\)}{d}-invariants}
\label{sec: categorical folding}

For a reduced sequence, the determinantial modules attached to \(\bfi\)-boxes realize the
box seeds of Section~\ref{sec: folding cluster algebras}.  We use
these realizations to convert the two ordered tropical folding
identities into formulas for the \(\Lambda\)-invariants involving sums over box fibers.
Symmetrizing the two formulas gives the folding formula for the
\(\de\)-invariant.

\subsection{Quiver Hecke categories and affreal modules}
\label{subsec: quiver Hecke box modules}
We refer to \cite{KKKO18} for the basic definitions and standard
results on quiver Hecke algebras and their finite-dimensional graded
module categories. 
Keep the Cartan and root data of
Section~\ref{sec: folding cluster algebras}.
Let \(\sfP\) be a weight lattice containing the fundamental weights
\((\varpi_i)_{i\in I}\), and set
\[
\sfP_+\seteq\sum_{i\in I}\Z_{\geq0}\varpi_i.
\]

Fix a quiver Hecke algebra \cite{KL1,KL2,R08}
\[
R=\soplus_{\beta\in\sfQ_+}R(\beta)
\]
over a field $\Bbbk$, associated with the Cartan datum of \(\g\).  We denote by
\(R\text{-gmod}\) the category of finite-dimensional graded
\(R\)-modules and by \(\conv\) its convolution product.

We generally suppress grading shifts from the notation and understand
isomorphisms of graded modules up to a grading shift unless stated
otherwise.  At the level of Grothendieck rings, we use the ungraded
specialization obtained by setting the grading parameter equal to
\(1\).  The grading is retained in the degrees of renormalized
\(R\)-matrices and hence in the invariants \(\Lambda(-,-)\) and
\(\de(-,-)\).

For simple \(R\)-modules \(M\) and \(N\), we say that they
\emph{commute} if
\[
M\conv N\simeq N\conv M
\]
up to a grading shift, and that they \emph{strongly commute} if, in
addition, \(M\conv N\) is simple.  A simple \(R\)-module \(M\) is
called \emph{real} if \(M\conv M\) is simple.

A simple \(R\)-module \(M\) is called \emph{affreal} if it is real
and admits an \emph{affinization} (\cite{KP18,KKOP21}).

For simple \(R\)-modules \(M\) and \(N\), suppose that
\[
\HOM_R(M\conv N,N\conv M)=\Bbbk\,\rmat{M,N}
\]
for a nonzero homogeneous homomorphism \(\rmat{M,N}\).  The map \(\rmat{M,N}\) is unique up to a nonzero scalar and is
called the \emph{renormalized \(R\)-matrix}.
We call the ordered pair \((M,N)\) \emph{\(\Lambda\)-definable}
and set
\[
\Lambda(M,N)\seteq\deg(\rmat{M,N}).
\]
If both ordered pairs $(M,N)$ and $(N,M)$ are
\(\Lambda\)-definable, we set
\[
\de(M,N)
\seteq
\frac12\bigl(\Lambda(M,N)+\Lambda(N,M)\bigr).
\]
In particular, if either \(M\) or \(N\) is affreal, then both
\((M,N)\) and \((N,M)\) are \(\Lambda\)-definable, and hence
\(\de(M,N)\) is defined.

A simple module \(S\) is called a \emph{convolution factor} of a
simple module \(M\) if there exists a simple module \(X\) such that
\[
M\simeq S\conv X
\]
up to a grading shift (\cite[Section~2.3]{KKOP26}).
 
We use the standard \(R\)-matrix calculus for affreal simple
modules, including composition and additivity under convolution
(\cite[Lemma~3.16]{KKOP21} and
\cite[Proposition~2.6]{KKOP26}).

We next recall normal sequences in the almost affreal form of
\cite[Section~2.3]{KKOP26}. 

\begin{definition}\label{def: normal sequence} \cite[Definition 2.5]{KK19}
Let \(M_1,\ldots,M_r\) be simple \(R\)-modules, all but at most one of
which are affreal. Let
\[
\begin{split}
\rmat{M_1,\ldots,M_r}
\seteq{}&
(\rmat{M_{r-1},M_r})\circ\cdots\circ
(\rmat{M_2,M_r}\circ\cdots\circ\rmat{M_2,M_3})\\
&\circ
(\rmat{M_1,M_r}\circ\cdots\circ\rmat{M_1,M_2}):
M_1\conv\cdots\conv M_r
\longrightarrow
M_r\conv\cdots\conv M_1.
\end{split}
\]
The ordered sequence
\( (M_1,\ldots,M_r) \) is called \emph{normal} if
\[
\rmat{M_1,\ldots,M_r}\ne0.
\]
\end{definition}

For the almost affreal form of the following lemma, see
\cite[Lemma~2.8 and Section~2.2]{KKOP23L}. 

\begin{lemma}[{\cite[Lemmas~2.6--2.8]{KK19}}]
\label{lem: basic properties of normal sequences}
Let \( (M_1,\ldots,M_r) \) be a normal sequence as in
Definition~\ref{def: normal sequence}.  Then
\[
{\rm Im}\bigl(\rmat{M_1,\ldots,M_r}\bigr) \simeq \head(M_1\conv\cdots\conv M_r) \simeq \soc(M_r\conv\cdots\conv M_1),
\]
and these modules are simple.  Moreover, every consecutive ordered 
subsequence of \( (M_1,\ldots,M_r) \) is normal, and
\begin{align*}
\Lambda\bigl( M_1,\head(M_2\conv\cdots\conv M_r) \bigr) &= \sum_{a=2}^r\Lambda(M_1,M_a),  \\
\Lambda\bigl(\head(M_1\conv\cdots\conv M_{r-1}),M_r \bigr) &= \sum_{a=1}^{r-1}\Lambda(M_a,M_r).
\end{align*}
\end{lemma}

For \(M\in R(\beta)\text{-gmod}\), its weight is
\(\wt(M)=-\beta\).  Set
\[
\mathrm W(M)
\seteq
\bigl\{
\gamma\in\sfQ_+\cap(\beta-\sfQ_+)
\mid
e(\gamma,\beta-\gamma)M\ne0
\bigr\}.
\]
Here \(e(\gamma,\beta-\gamma)\in R(\beta)\) is the standard quiver
Hecke idempotent associated with the decomposition
\(\beta=\gamma+(\beta-\gamma)\) (see \cite[Section~2.3]{KKOP18} for instances).
For \(w\in W_\g\), let \(\scrC_w\) be the full subcategory of
\(R\text{-gmod}\) consisting of the modules \(M\) such that
\[
\mathrm W(M)\subset\sfQ_+\cap w\sfQ_-.
\]
The category \(\scrC_w\) is stable under convolution products,
extensions, subquotients, and grading shifts.  The relation between
the Grothendieck ring of \(\scrC_w\) and the cluster algebra associated
with \(w\) will be recalled in the next subsection.
We set
\[
K(\scrC_w)\seteq\C\otimes_{\Z}K_0(\scrC_w).
\]

\subsection{Determinantial modules and categorical realizations}
\label{subsec: categorical realizations of box seeds}

We recall the determinantial modules attached to \(\bfi\)-boxes
and their realization of box seeds in \(\scrC_w\).
 
For \(\uplambda\in\sfP_+\) and \(\mu,\zeta\in W_\g\uplambda\) with
\(\mu\preceq\zeta\), let
\[
M(\mu,\zeta)
\]
denote the determinantial \(R(\zeta-\mu)\)-module associated with
\((\mu,\zeta)\).  It is a self-dual real simple module by
\cite[Propositions~4.1 and~4.2]{KKOP18} and admits an affinization by
\cite[Theorem~3.26]{KKOP21}.  Hence \(M(\mu,\zeta)\) is affreal.

Let \(w\in W_\g\), and set
\[
\Red(w) \seteq \left\{ (i_1,\ldots,i_l)\in I^l \mathrel{}\middle|\mathrel{} w=s_{i_1}\cdots s_{i_l},\quad l=\ell(w) \right\}.
\]
Fix $\ii=(i_1,\ldots,i_l)\in\Red(w)$ and define
\[
w_{\leq k}^\ii \seteq s_{i_1}\cdots s_{i_k}, \qquad w_{<k}^\ii \seteq s_{i_1}\cdots s_{i_{k-1}} \qquad(1\leq k\leq l).
\]

For \(1\leq k\leq l\), define the \emph{cuspidal module} associated with
the position \(k\) of \(\ii\) by
\begin{equation}\label{eq: cuspidal determinantial module}
S_k^\ii
\seteq
M\bigl(
w_{\leq k}^\ii\varpi_{i_k},
w_{<k}^\ii\varpi_{i_k}
\bigr).
\end{equation}

Let $\frakc=[x,y]$ be an \(\bfi\)-box associated with \(\ii\).  Thus
\[
1\leq x\leq y\leq l \qtq i_x=i_y=i_\frakc.
\]
We define the determinantial module associated with \(\frakc\) by
\begin{equation}\label{eq: determinantial module of a box}
M^\ii[\frakc] \seteq M^\ii[x,y] \seteq M\bigl( w_{\leq y}^\ii\varpi_{i_y}, w_{<x}^\ii\varpi_{i_x} \bigr).
\end{equation}
In particular, $M^\ii[k]=S_k^\ii.$ 

Let \(\frakc=[x,y]\) be an \(\bfi\)-box associated with \(\ii\), and
write
\[
x=x_0<x_1<\cdots<x_{r-1}=y
\]
for the positions of color \(i_\frakc\) in \([x,y]\).
 Then \cite[Proposition~4.3]{KKOP23F} and
\cite[Section~5.1.1]{KK24} give
\begin{equation}\label{eq: normal realization of box module}
\bigl(S_{x_{r-1}}^\ii,\ldots,S_{x_1}^\ii,S_{x_0}^\ii\bigr) \text{ is a normal sequence}
\qtq 
M^\ii[\frakc]
\simeq
\head\bigl( S_{x_{r-1}}^\ii\conv\cdots\conv S_{x_1}^\ii\conv S_{x_0}^\ii \bigr).
\end{equation}

The following formula for one order of the \(\Lambda\)-invariant is the form of
\cite[Proposition~5.7, (5.5)--(5.6)]{KKOP23F} in our sign
convention.

\begin{proposition}\label{prop: KKOP directional box Lambda formula}
Let \([x,y]\) and \([x',y']\) be two \(\ii\)-boxes for the same
reduced sequence \(\ii\), and put
\[
 M=M^\ii[x,y],\qquad N=M^\ii[x',y'].
\]
If $x>(x')^-$ or $y^+>y'$, then
\begin{equation}\label{eq: KKOP directional box Lambda formula}
 \Lambda(M,N)=-\lambda^\ii(\calE_{[x,y]},\calE_{[x',y']}).
\end{equation}
Consequently,
\[
 \Lambda(N,M)=2\de(M,N)
 +\lambda^\ii(\calE_{[x,y]},\calE_{[x',y']}).
\]
After interchanging the boxes, if
$ x'>x^-$ or $(y')^+>y$,
the same proposition computes \(\Lambda(N,M)\) directly.
\end{proposition}

\begin{remark}\label{rmk: one Lambda is computable}
Let \(M=M^\ii[x,y]\) and \(N=M^\ii[x',y']\) be two box modules
associated with the same reduced sequence \(\ii\).
At least one of the two orders satisfies the hypotheses of
Proposition~\ref{prop: KKOP directional box Lambda formula}.
Indeed, if \(x\ge x'\), then \(x>(x')^-\), while if \(x'>x\),
then \(x'>x^-\).
Thus the pairing \(\lambda^\ii\) determines one of
\(\Lambda(M,N)\) and \(\Lambda(N,M)\), and the identity
\[
\Lambda(M,N)+\Lambda(N,M)=2\de(M,N)
\]
determines the other once \(\de(M,N)\) is known.

If the hypotheses hold in both orders, skew-symmetry of
\(\lambda^\ii\) gives \(\de(M,N)=0\).
\end{remark}

Let \(K\) be a finite index set and let
\[
\sfM=\{M_k\mid k\in K\}
\]
be a family of affreal simple modules.  We call \(\sfM\) a
\emph{commuting family} if
\[
M_k\conv M_l\simeq M_l\conv M_k
\qquad(k,l\in K)
\]
up to a grading shift.  For
\(\bfa=(a_k)_{k\in K}\in\Z_{\geq0}^{\oplus K}\), the corresponding
convolution product is simple and independent of the ordering of its
factors, up to a grading shift.  We denote it by
\[
\sfM(\bfa)
\seteq
\mathop{\conv}_{k\in K}M_k^{\conv a_k}.
\]
Thus
\[
\sfM(\bfa)\conv\sfM(\bfb)
\simeq
\sfM(\bfa+\bfb)
\qquad
\bigl(\bfa,\bfb\in\Z_{\geq0}^{\oplus K}\bigr)
\]
up to a grading shift.

A simple module \(X\) is called an \emph{\(\sfM\)-monomial} if
\(X\simeq\sfM(\bfa)\) for some
\(\bfa\in\Z_{\geq0}^{\oplus K}\).  The commuting family \(\sfM\)
is called \emph{independent} if
\(\sfM(\bfa)\simeq\sfM(\bfb)\) implies \(\bfa=\bfb\). 

\smallskip

Let \(\frakC\) be an admissible chain of \(\bfi\)-boxes associated
with \(\ii\), and set
\[
\calF\seteq\calF(\frakC).
\]
For \(\frakc\in\calF\), write
\[
M_\frakc^\ii\seteq M^\ii[\frakc],
\qquad
\sfM^\ii_\calF
\seteq
\{M_\frakc^\ii\mid\frakc\in\calF\}.
\]

\begin{lemma}[{\cite[Theorem~4.21(b)]{KKOP24A};
\cite[Proposition~5.1]{KK24}}]
\label{lem: commuting box modules}
Let \(\frakc\) and \(\frakc'\) be commuting \(\bfi\)-boxes associated
with \(\ii\).  Then
\[
M^\ii[\frakc] \text{ and } M^\ii[\frakc'] \text{ strongly commute.}
\]
\end{lemma}

\begin{proposition}\label{prop: categorical box family}
The family \(\sfM^\ii_\calF\) is an independent commuting family of
affreal simple objects in \(\scrC_w\).  More precisely, for any
\(\frakc,\frakc'\in\calF\), the modules $M_\frakc^\ii$ and $M_{\frakc'}^\ii$ 
strongly commute.  Consequently, for every
$\bfa=(a_\frakc)_{\frakc\in\calF}
\in\Z_{\geq0}^{\oplus\calF}$,
the convolution product
\[
\sfM^\ii_\calF(\bfa)
\seteq
\mathop{\conv}_{\frakc\in\calF}
\bigl(M_\frakc^\ii\bigr)^{\conv a_\frakc} \text{ is simple and independent of the ordering of its factors,}
\]
up to a grading shift.
\end{proposition}

\begin{proof}
The affreal property and
Lemma~\ref{lem: commuting box modules}
give the assertions on convolution products.
The independence follows from the algebraic independence of
the corresponding cluster variables under the
Grothendieck-ring realization in
Theorem~\ref{thm: categorical realization of box seeds}.
\end{proof}

We next recall the cluster algebra structure on the Grothendieck
ring of \(\scrC_w\).
Let \(A_q(\mathfrak n(w))\) denote the quantum unipotent coordinate
algebra associated with \(w\), and write
\(A_q(\mathfrak n(w))_{\mathbb Z[q^{\pm1}]}\) for its
\(\mathbb Z[q^{\pm1}]\)-form.
Let \(K_q(\scrC_w)\) denote the \emph{graded} Grothendieck ring of
\(\scrC_w\), where a grading shift acts by a power of \(q\).
The quiver Hecke categorification induces an isomorphism of
\(\mathbb Z[q^{\pm1}]\)-algebras
\begin{equation}\label{eq: graded Grothendieck ring of Cw}
\Omega_w:
K_q(\scrC_w)
\xrightarrow{\sim}
A_q(\mathfrak n(w))_{\mathbb Z[q^{\pm1}]}
\end{equation}
(see \cite[Theorem~2.20(ii-a)]{KKOP18}).
Under this isomorphism, the class of a
determinantial module \(M(\mu,\zeta)\) is sent to the corresponding unipotent quantum minor
\(D(\mu,\zeta)\) (see \cite[Proposition~4.1]{KKOP18} for more detail):
\[
\Omega_w([M(\mu,\zeta)])=D(\mu,\zeta).
\]

The algebra \(A_q(\mathfrak n(w))\) admits a quantum cluster algebra
structure associated with every reduced sequence
\(\ii\in\Red(w)\).  In the symmetric setting, this structure
originates from the work of Geiß--Leclerc--Schröer, while its
symmetrizable counterpart is due to Goodearl--Yakimov \cite{GLS11,GLS13,GY17,GY21}.

Let
\[
\frakC_R^\ii \seteq \bigl(\{1,k]\bigr)_{1\leq k\leq l} 
\qtq
\calF_R^\ii \seteq \bigl\{\{1,k]\mid1\leq k\leq l\bigr\}
\]
be the \emph{constant \(\calR\)-chain} and its underlying maximal commuting
family.  For \(1\leq k\leq l\), one has
\[
M^\ii[\{1,k]] \simeq M\bigl( w_{\leq k}^\ii\varpi_{i_k}, \varpi_{i_k} \bigr).
\]
Consequently,
\[
\Omega_w\bigl([M^\ii[\{1,k]]\bigr)
= 
D\bigl( w_{\leq k}^\ii\varpi_{i_k}, \varpi_{i_k} \bigr).
\]
These unipotent quantum minors form the initial cluster associated with \(\ii\).

We transport this quantum cluster algebra structure along
\(\Omega_w\).  Specializing the grading parameter to \(1\), we
obtain a \(\C\)-algebra isomorphism
\[
\Phi_w:
K(\scrC_w)
\xrightarrow{\sim}
\calA(w),
\]
where \(\calA(w)\) is the corresponding classical cluster algebra.
We use the opposite exchange matrices in the specialized
classical seed pattern, in accordance with our convention for
\(B^\calF\). 
Under this isomorphism,
\begin{equation}\label{eq: constant R categorical seed}
\Phi_w\bigl([M_\frakc^\ii]\bigr)
=
x_\frakc
\qquad
(\frakc\in\calF_R^\ii).
\end{equation}
Thus the family
\[
\sfM_R^\ii
\seteq
\{M_\frakc^\ii\mid\frakc\in\calF_R^\ii\} \text{ realizes the initial box seed $\Sigma^\ii(\calF_R^\ii)$.}
\]

 We next show that mutation-type box moves are realized
by categorical mutations.  Let
$\frakC=(\frakc_k)_{1\leq k\leq l}$
be an admissible chain of \(\bfi\)-boxes associated with \(\ii\),
and set
$\calF\seteq\calF(\frakC)$
and
$\sfM^\ii_\calF = \{M_\frakc^\ii\mid\frakc\in\calF\}$.

The box-move terminology and the exchanged pair
\(\{[u^+,v],[u,v^-]\}\) are as in
Lemma~\ref{lem: two types of box moves}.

\begin{lemma}[{\cite[Proposition~5.2]{KK24}}]
\label{lem: box T-system}
 Let \([u,v]\) be an \(\bfi\)-box of color \(i\) with \(u<v\).  Up to grading
shifts, there is a short exact sequence 
\begin{equation} \label{eq: T-system}
\begin{aligned}
0\longrightarrow
\mathop{\conv}_{\substack{j\in I\setminus\{i\}\\c_{j,i}<0}}
M^\ii[u(j)^+,v(j)^-]^{\conv(-c_{j,i})}
& \longrightarrow 
M^\ii[u^+,v]\conv M^\ii[u,v^-] \\
& \longrightarrow
M^\ii[u^+,v^-]\conv M^\ii[u,v]
\longrightarrow0.    
\end{aligned}
\end{equation}
An undefined or empty interval is omitted, and an empty convolution
product is understood to be the monoidal unit. We call~\eqref{eq: T-system} the \emph{\(T\)-system}. 
\end{lemma}

\begin{proposition}\label{prop: box moves are categorical mutations}
Let $\frakC'\seteq\frakB_k(\frakC)$ be a mutation-type box move replacing
\[
\frakc\in\calF(\frakC) \quad\text{with}\quad \frakc^\star\in\calF(\frakC').
\]
Assume that
\[
\left(\bigl([M_\frake^\ii]\bigr)_{\frake\in\calF(\frakC)},B^{\calF(\frakC)}\right) \text{ is a seed of \(K(\scrC_w)\). }
\]
Then, up to the natural relabeling,
\[
\mu_\frakc \left( \bigl([M_\frake^\ii]\bigr)_{\frake\in\calF(\frakC)}, B^{\calF(\frakC)} \right)
= \left( \bigl([M_\frake^\ii]\bigr)_{\frake\in\calF(\frakC')}, B^{\calF(\frakC')} \right).
\]
\end{proposition}
\begin{proof}
Put \(\calF=\calF(\frakC)\) and \(\calF'=\calF(\frakC')\).
By Lemma~\ref{lem: two types of box moves}, there is an
\(\bfi\)-box \([u,v]\) of color \(i\) such that
\[
\{\frakc,\frakc^\star\}=\{[u^+,v],[u,v^-]\},
\qquad
\frakc_{k+1}=[u,v].
\]
Both box modules are affreal, and each strongly commutes with
the modules indexed by \(\calF\cap\calF'\).

The first \(k-1\) boxes are unchanged by the move and have
envelope \([u+1,v-1]\).
For each color \(j\ne i\) occurring in this interval, set  
\[
r_j=\max\{r<k\mid i_{\frakc_r}=j\}.
\]
By the definition of an admissible chain,
\[
\frakc_{r_j}=[u(j)^+,v(j)^-].
\]
If \([u^+,v^-]\) is nonempty, it is the last box of color \(i\)
among the first \(k-1\) boxes.
Thus all nonempty boxes appearing in the two end terms of
\eqref{eq: T-system} belong to \(\calF\cap\calF'\).

By \cite[(3.1)]{KK24}, these end terms correspond
to the two exchange monomials.  After specializing the grading
parameter to \(1\), we obtain
\[
[M_\frakc^\ii][M_{\frakc^\star}^\ii]
=
\prod_{\frake\in\calF}
[M_\frake^\ii]^{[b_{\frake,\frakc}^\calF]_+}
+
\prod_{\frake\in\calF}
[M_\frake^\ii]^{[-b_{\frake,\frakc}^\calF]_+}.
\]
Hence the new cluster variable is \([M_{\frakc^\star}^\ii]\).
The box move also gives
\(B^{\calF'}=\mu_\frakc(B^\calF)\), proving the asserted
seed identity.
\end{proof}

\begin{theorem}[{\cite[Theorem~5.20]{KK24}}]
\label{thm: categorical realization of box seeds}
Let \(\frakC\) be an admissible chain of \(\bfi\)-boxes associated
with \(\ii\), and set \(\calF\seteq\calF(\frakC)\).
Then
\[
\left(
\bigl([M_\frakc^\ii]\bigr)_{\frakc\in\calF},
B^\calF
\right)
\]
is a seed of \(K(\scrC_w)\), identified under \(\Phi_w\) with
the box seed \(\Sigma^\ii(\calF)\).  Moreover,
\[
\Phi_w([M_\frakc^\ii])=x_\frakc
\quad(\frakc\in\calF),
\qquad
L^\calF
=-\bigl(
\Lambda(M_\frakc^\ii,M_{\frakc'}^\ii)
\bigr)_{\frakc,\frakc'\in\calF}.
\]
\end{theorem}

\begin{proof}
The seed identification holds for the constant \(\calR\)-chain by
\eqref{eq: constant R categorical seed}.
Every admissible chain associated with \(\ii\) is connected to
\(\frakC_R^\ii\) by box moves.
Permutation-type moves only relabel the seed, and mutation-type
moves are treated by
Proposition~\ref{prop: box moves are categorical mutations}.
Induction along these moves proves the seed identification.
The identity for \(L^\calF\) follows from
\cite[Propositions~5.7 and~5.9]{KKOP23F}
and our sign convention for \(\lambda^\ii\).
\end{proof}

\begin{corollary}\label{cor: every box module is a cluster variable}
For every \(\bfi\)-box \(\frakc\) associated with \(\ii\), 
\[
\text{the class $[M^\ii[\frakc]]$ is a cluster variable of \(K(\scrC_w)\).}
\]
More precisely, there
exists an admissible chain \(\frakC\) such that
\[
\frakc\in\calF(\frakC) \qtq  \Phi_w([M^\ii[\frakc]])=x_\frakc
\]
in the corresponding box seed.
\end{corollary}

\begin{proof}
Every \(\bfi\)-box associated with \(\ii\) is contained in the
underlying maximal commuting family of an admissible chain.  The
assertion follows from
Theorem~\ref{thm: categorical realization of box seeds}.
\end{proof}

\subsection{Quasi-Laurent families and the categorical
\texorpdfstring{\(F\)}{F}-invariant}
\label{subsec: quasi-Laurent families}

Following \cite[Definition~3.2]{KKOP23F},   an independent commuting
family \(\sfM\) of affreal simple objects of a full monoidal
subcategory \(\scrC\subset R\text{-gmod}\) is called a
\emph{quasi-Laurent family in \(\scrC\)}
if, for every simple object \(X\in\scrC\) commuting with all
\(M_k\), there exist 
\(\bfa,\bfb\in\Z_{\geq0}^{\oplus K}\) such that
\[
X\conv\sfM(\bfa)\simeq\sfM(\bfb)
\]
up to a grading shift.  The family \(\sfM\) is called a \emph{Laurent family in \(\scrC\)}
if every simple object of \(\scrC\) commuting with all \(M_k\)
is an \(\sfM\)-monomial.  Every Laurent family is, in particular, a
quasi-Laurent family.

\begin{proposition}[{\cite[Proposition~3.6 and Corollary~3.4]{KKOP23F}}]
\label{prop: Laurent maximum realization of box seed}
 Let \(\frakC\) be an admissible chain of \(\bfi\)-boxes associated
with \(\ii\), and set
$\calF\seteq\calF(\frakC)$
and
$\sfM^\ii_\calF = \{M_\frakc^\ii\mid\frakc\in\calF\}$.
Then \(\sfM^\ii_\calF\) is a quasi-Laurent family in
\(\scrC_w\), and 
\[
\Phi_w([M_\frakc^\ii])=x_\frakc
\qquad(\frakc\in\calF),
\qtq
L^\calF=-\Lambda_{\sfM^\ii_\calF}.
\]
In particular, every simple object \(X\in\scrC_w\) has a positive
Laurent expansion
\begin{equation}\label{eq: categorical Laurent expansion for box seed}
\Phi_w([X])
=
\sum_{\nu=1}^r c_\nu(\bfx^\calF)^{\bfv_\nu},
\qquad
c_\nu\in\Z_{>0}.
\end{equation}
\end{proposition}
\begin{proof}
The identification \(\Phi_w([M_\frakc^\ii])=x_\frakc\) follows
from Theorem~\ref{thm: categorical realization of box seeds},
and \(L^\calF=-\Lambda_{\sfM^\ii_\calF}\) follows from
\cite[Propositions~5.7 and~5.9]{KKOP23F}.
The quasi-Laurent property and the positive Laurent expansion
then follow from \cite[Proposition~3.6 and Corollary~3.4]{KKOP23F},
using the quantum seed associated with the box family in
\cite[Section~5.1.4 and Theorem~5.20]{KK24}.
\end{proof}

Let \(\scrC\) be a full monoidal subcategory of \(R\text{-gmod}\)
stable under subquotients, extensions, and grading shifts.
Let \(K_q(\scrC)\) denote its graded Grothendieck ring.
For a quasi-Laurent family
\(\sfM=\{M_k\mid k\in K\}\) in \(\scrC\), let
\(\calT_{\sfM}\) be the quantum torus over \(\Z[q^{\pm1/2}]\)
generated by the classes \([M_k]\) and their inverses.
We denote by
\[
\phi_{\sfM}:
\Z[q^{\pm1/2}]\tens_{\Z[q^{\pm1}]}K_q(\scrC)
\hookrightarrow \calT_{\sfM}
\]
the injective algebra homomorphism given by
\cite[Lemma~2.12]{KKOP26}. 
For \(\bfa\in\Z^{\oplus K}\), the notation
\([\sfM(\bfa)]\) denotes the corresponding Laurent monomial in
\(\calT_{\sfM}\); when
\(\bfa\in\Z_{\geq0}^{\oplus K}\), it is represented by the class
of the corresponding convolution module.

\begin{theorem}[{\cite[Theorem~2.13]{KKOP26}}]
\label{thm: KKOP Laurent maximum}
With \(\scrC\) and \(\sfM\) as above, let \(X,L\in\scrC\) be simple
modules such that \(L\) is affreal, and suppose that the Laurent
expansion of \([X]\) with respect to \(\sfM\) is
\[
\phi_{\sfM}([X])
=
\sum_\nu q^{c_\nu}[\sfM(\bfa_\nu)],
\qquad
c_\nu\in\Z, \quad
\bfa_\nu\in\Z^{\oplus K}.
\]
Then we have 
\[
\Lambda(L,X) = \max_\nu\Lambda\bigl(L,\sfM(\bfa_\nu)\bigr)
\qtq
\Lambda(X,L) = \max_\nu\Lambda\bigl(\sfM(\bfa_\nu),L\bigr).
\]
\end{theorem}

\begin{theorem}\label{thm: categorical interpretation of F-invariant}
 Let \(X\) be an affreal simple object of \(\scrC_w\) such that
\(z_X\seteq\Phi_w([X])\) is good and bipointed.  For
\(\bfa\in\Z_{\geq0}^{\oplus\calF}\), set
\[
N\seteq\sfM^\ii_\calF(\bfa).
\]
Then
\[
\Lambda(N,X) = \langle z_X,\Phi_w([N])\rangle_{\trop}
\qtq
\Lambda(X,N) = \langle\Phi_w([N]),z_X\rangle_{\trop},
\]
and consequently
\[
2\de(X,N) = \bigl(z_X\Vert\Phi_w([N])\bigr)_F
\in2\Z_{\geq0}.
\]
\end{theorem}

\begin{proof}
Choose a graded positive Laurent expansion supplied by
\cite[Corollary~3.4]{KKOP23F}.   After specializing \(q=1\) and
collecting equal Laurent monomials, it gives
\eqref{eq: categorical Laurent expansion for box seed}.    
Since \(z_X\) is good and bipointed, its exponents have the form
\[
\bfv_\nu
=
\bfg_{X;\calF}+B^\calF\bfw_\nu,
\qquad
\bfw_\nu\in\Z_{\geq0}^{\oplus\calF_{\ex}},
\]
where both \(\bfw_\nu=0\) and
\(\bfw_\nu=\bff_{X;\calF}\) occur, and
\[
0\leq(\bfw_\nu)_\frakc
\leq f_{\frakc;X;\calF}
\qquad(\frakc\in\calF_{\ex}).
\]

For each \(\frakc\in\calF\) with \(a_\frakc>0\), the affreal
simple module \(M_\frakc^\ii\) is a convolution factor of \(N\).
Apply Theorem~\ref{thm: KKOP Laurent maximum} with
\(L=M_\frakc^\ii\), and then sum with multiplicity
\(a_\frakc\).  Since \(N=\sfM^\ii_\calF(\bfa)\) is simple,
\cite[Proposition~2.6]{KKOP26} gives
\begin{equation}\label{eq: Lambda additivity for seed monomials}
\Lambda(N,X)
=
\sum_{\frakc\in\calF}
a_\frakc\Lambda(M_\frakc^\ii,X) \qtq 
\Lambda(X,N)
=
\sum_{\frakc\in\calF}
a_\frakc\Lambda(X,M_\frakc^\ii).
\end{equation} 

Set
\(\bfa_{\ex}\seteq (a_\frakc)_{\frakc\in\calF_{\ex}}\). 
The compatibility relation
\[
(B^\calF)^TL^\calF=(2D^\calF\mid0)
\quad
\text{gives}
\quad 
\bfv_\nu^TL^\calF\bfa
=
\bfg_{X;\calF}^TL^\calF\bfa
+\bfw_\nu^T(2D^\calF)\bfa_{\ex}.
\]
Thus the term \(\bfw_\nu=\bff_{X;\calF}\) is a common
maximizer for all these convolution factors and their
multiplicities.  Theorem~\ref{thm: KKOP Laurent maximum}
and~\eqref{eq: Lambda additivity for seed monomials} give
\begin{align}\label{eq: Lambda N X as tropical invariant}
\Lambda(N,X)
&=
\bfg_{X;\calF}^TL^\calF\bfa
+F_{X;\calF}[2D^\calF\bfa_{\ex}]
=
\langle z_X,\Phi_w([N])\rangle_{\trop}.
\end{align}

On the other hand, skew-symmetry and compatibility give
\[
L^\calF B^\calF
=
-\begin{pmatrix}2D^\calF\\0\end{pmatrix},
\qtqh
\bfa^TL^\calF\bfv_\nu
=
\bfa^TL^\calF\bfg_{X;\calF}
-\bfa_{\ex}^T(2D^\calF)\bfw_\nu.
\]
Here \(\bfw_\nu=0\) is a common maximizer, and therefore
\begin{align}\label{eq: Lambda X N as tropical invariant}
\Lambda(X,N)
&=
\bfa^TL^\calF\bfg_{X;\calF}
=
\langle\Phi_w([N]),z_X\rangle_{\trop}.
\end{align}
The last assertion follows by adding the two identities and applying Theorem~\ref{thm: Cao F-invariant}.
\end{proof}

By Proposition~\ref{prop: cluster monomials are bipointed} and
universal positivity, the good and bipointed hypothesis is
automatic when \(z_X\) is a cluster variable or, more generally,
a cluster monomial.

\begin{corollary}[{cf.~\cite[Proposition~4.13]{Cao23}}]
\label{cor: d-invariant from f-vector}
Let \(\frakc\in\calF_{\ex}\), and let \(X\) be an affreal
simple object of \(\scrC_w\) such that \(\Phi_w([X])\) is good
and bipointed.  Then
\[
\de(M_\frakc^\ii,X)
=
d_{i_\frakc}f_{\frakc;\Phi_w([X]);\calF}.
\]
\end{corollary}

\subsection{Folding formulas for
\texorpdfstring{\(\Lambda\)}{Lambda}- and
\texorpdfstring{\(\de\)}{d}-invariants}
\label{subsec: categorical folding formulas}

We now combine the categorical interpretation of the
\(F\)-invariant with the folding of \(F\)-polynomials.

Let \(\frakC\) be an admissible chain of \(\bfi\)-boxes associated
with \(\ii\), and let $\hfrakC$
be its canonical lift.  Set
\[
\calF\seteq\calF(\frakC) \qtq \hcalF\seteq\calF(\hfrakC).
\]
For \(\frakc\in\calF\), recall that $\Fib_{\frakC}(\frakc) \subset\hcalF$ 
denotes its box fiber.

Set
$\hw \seteq \hs_{i_1}\cdots\hs_{i_l}\in W_\hg$.
Let $\widehat{\scrC}_{\hw}$ be the quiver Hecke category associated with the lifted reduced
sequence \(\hii\).  For
\(\hfrakc\in\hcalF\), write
\[
\hM_{\hfrakc} \seteq M^\hii[\hfrakc].
\]
Let
$\widehat\Phi_{\hw}: K(\widehat{\scrC}_{\hw}) \xrightarrow{\sim} \calA(\hw)$
denote the upstairs counterpart of \(\Phi_w\).  The lifted box family categorically realizes
\(\hSigma^\hii(\hcalF)\), and hence
\[
\widehat\Phi_{\hw} \bigl([\hM_{\hfrakc}]\bigr) = \hx_{\hfrakc}.
\]

For \(\frakc\in\calF\),   define the convolution product 
\[
\hM_{\Fib_{\frakC}(\frakc)} \seteq \mathop{\conv}_{\hfrakc\in\Fib_{\frakC}(\frakc)} \hM_{\hfrakc}.
\]
The modules occurring in this convolution belong to a common lifted
box seed and hence pairwise strongly commute. 
Thus 
\[
\widehat\Phi_{\hw} \bigl([\hM_{\Fib_{\frakC}(\frakc)}]\bigr) = \hx_{\Fib_{\frakC}(\frakc)} \seteq
\prod_{\hfrakc\in\Fib_{\frakC}(\frakc)} \hx_{\hfrakc}.
\]

\begin{theorem}\label{thm: categorical folding of directional Lambda invariants}
Keep the preceding notation.  Let \(\frakc\in\calF\), and let
\(N\in\scrC_w\) be simple such that \(z_N=\Phi_w([N])\) is a
cluster variable.   Suppose that \(z_N\) and \(\hz_N\) are corresponding cluster variables
in the box seeds of~\eqref{eq: global box seed folding}, with
\(\phi_{\frakC}(\hz_N)=z_N\).  Let \(\hN\) be the simple object of
\(\widehat{\scrC}_{\hw}\) with
\(\widehat\Phi_{\hw}([\hN])=\hz_N\). Then
\begin{subequations} \label{eq: categorical Lambda fiber sum}
\begin{align}
\label{eq: categorical first Lambda fiber sum}
\Lambda(M_\frakc^\ii,N)
&=d_{i_\frakc}\sum_{\hfrakc\in\Fib_{\frakC}(\frakc)}
\Lambda(\hM_{\hfrakc},\hN),\\
\label{eq: categorical second Lambda fiber sum}
\Lambda(N,M_\frakc^\ii)
&=d_{i_\frakc}\sum_{\hfrakc\in\Fib_{\frakC}(\frakc)}
\Lambda(\hN,\hM_{\hfrakc}).
\end{align}
\end{subequations}
\end{theorem}

\begin{proof}
The cluster variables in the statement are good and bipointed, and
their categorical realizations are affreal.  We may therefore use
Theorem~\ref{thm: categorical interpretation of F-invariant} in
both box seeds.  With the order convention in that theorem, we have
\[
 \Lambda(M_\frakc^\ii,N)=\langle z_N,x_\frakc\rangle_{\trop}
 \qtq 
 \Lambda(N,M_\frakc^\ii)=\langle x_\frakc,z_N\rangle_{\trop},
\]
and the analogous identities upstairs.  Applying respectively
\eqref{eq: folding of second ordered tropical pairing} and
\eqref{eq: folding of first ordered tropical pairing} from
Theorem~\ref{thm: folding of F-invariants between box variables}
therefore gives
\[
\begin{aligned}
 \Lambda(M_\frakc^\ii,N)
 =d_{i_\frakc}\Lambda(\hM_{\Fib_{\frakC}(\frakc)},\hN) \qtq
 \Lambda(N,M_\frakc^\ii)
 =d_{i_\frakc}\Lambda(\hN,\hM_{\Fib_{\frakC}(\frakc)}).
\end{aligned}
\]
Finally, applying \cite[Proposition~2.6]{KKOP26}
to \(\hM_{\Fib_{\frakC}(\frakc)}
=\displaystyle\mathop{\conv}_{\hfrakc\in\Fib_{\frakC}(\frakc)}
\hM_{\hfrakc}\) in both orders gives
\eqref{eq: categorical Lambda fiber sum}.  
\end{proof}

\begin{theorem}\label{thm: categorical folding of d-invariants}
Let \(\frakc\in\calF\), and let \(N\) be a simple object of
\(\scrC_w\) such that
\[
z_N\seteq\Phi_w([N])
\]
is a cluster variable. Suppose that \(z_N\) and \(\hz_N\) are corresponding cluster variables
in the box seeds of~\eqref{eq: global box seed folding} and satisfy 
\[
\phi_{\frakC}(\hz_N)=z_N,
\]
and let \(\hN\) be the corresponding simple object of
\(\widehat{\scrC}_{\hw}\), determined by
\[
\widehat\Phi_{\hw}([\hN]) = \hz_N.
\]
Then we have 
\begin{equation}\label{eq: categorical fiber sum formula}
\de\bigl(M_\frakc^\ii,N\bigr) = d_{i_\frakc} \sum_{\hfrakc\in\Fib_{\frakC}(\frakc)} \de\bigl(\hM_{\hfrakc},\hN\bigr) = d_{i_\frakc}\, \de\bigl(\hM_{\Fib_{\frakC}(\frakc)},\hN\bigr).
\end{equation}
\end{theorem}

\begin{proof}
Add equations~\eqref{eq: categorical first Lambda fiber sum} and
\eqref{eq: categorical second Lambda fiber sum}, and divide by $2$.
The additivity already used in their proof identifies the
symmetrized upstairs sum with
\(\de(\hM_{\Fib_{\frakC}(\frakc)},\hN)\). 
\end{proof}

\begin{corollary}\label{cor: folding of d-values between box modules}
Let \(\ii,\jj\in\Red(w)\), and let \(\frakc\) and \(\frakc'\) be
\(\bfi\)-boxes associated with \(\ii\) and \(\jj\), respectively.
Choose admissible chains \(\frakC\) and \(\frakD\) whose underlying
families contain \(\frakc\) and \(\frakc'\), respectively,
and choose
$\hfrakc'\in\Fib_{\frakD}(\frakc')$.
Then we have formulas for both orders of the \(\Lambda\)-invariant:
\begin{align}
\label{eq: folded first Lambda between box modules}
\Lambda\bigl(M^\ii[\frakc],M^\jj[\frakc']\bigr)
&=d_{i_\frakc}\sum_{\hfrake\in\Fib_{\frakC}(\frakc)}
\Lambda\bigl(M^\hii[\hfrake],M^\hjj[\hfrakc']\bigr),\\
\label{eq: folded second Lambda between box modules}
\Lambda\bigl(M^\jj[\frakc'],M^\ii[\frakc]\bigr)
&=d_{i_\frakc}\sum_{\hfrake\in\Fib_{\frakC}(\frakc)}
\Lambda\bigl(M^\hjj[\hfrakc'],M^\hii[\hfrake]\bigr).
\end{align}
Their symmetrization gives
\begin{equation}\label{eq: folded d-value between box modules}
\de\bigl( M^\ii[\frakc], M^\jj[\frakc'] \bigr)
= d_{i_\frakc} \sum_{\hfrake\in\Fib_{\frakC}(\frakc)} \de\bigl( M^\hii[\hfrake], M^\hjj[\hfrakc'] \bigr).
\end{equation}
The right-hand side is independent of the choice of the lifted box \(\hfrakc'\).
\end{corollary}

\begin{proof}
By Corollary~\ref{cor: every box module is a cluster variable},
$z_{\frakc'} \seteq \Phi_w\bigl([M^\jj[\frakc']]\bigr)$
is a cluster variable.  For every
\(\hfrakc'\in\Fib_{\frakD}(\frakc')\), set
\[
\hz_{\hfrakc'} \seteq \widehat\Phi_{\hw} \bigl([M^\hjj[\hfrakc']]\bigr).
\]
Theorem~\ref{thm: global folding for positive braids} gives
\[
\phi_{\frakC}^{\ttb} \bigl(\hz_{\hfrakc'}\bigr)
= z_{\frakc'}
\qquad
\bigl(\hfrakc'\in\Fib_{\frakD}(\frakc')\bigr).
\]
Hence Theorem~\ref{thm: categorical folding of directional Lambda invariants},
applied with \(N=M^\jj[\frakc']\) and
\(\hN=M^\hjj[\hfrakc']\), gives
\eqref{eq: folded first Lambda between box modules}
and~\eqref{eq: folded second Lambda between box modules}.
Adding these identities and dividing by $2$ gives
\eqref{eq: folded d-value between box modules}.
Since the left-hand side is independent of
\(\hfrakc'\), the right-hand side has the same value for every
\(\hfrakc'\in\Fib_{\frakD}(\frakc')\).
\end{proof}

In the notation of
Example~\ref{ex: A5 C3 folding of F-invariants}, the categorical
counterpart of the two box variables is
\[
x_{[2,5]}\longleftrightarrow M^\ii[2,5],
\qquad
x_{[7,9]}\longleftrightarrow M^\ii[7,9],
\]
and the relevant lifted variables correspond to
\(M^\hii[3,8]\), \(M^\hii[4,9]\), and \(M^\hii[12,15]\).
The \(A_5\) involution exchanges \(\hx_{[3,8]}\) and
\(\hx_{[4,9]}\), fixes \(\hx_{[12,15]}\), and preserves the lifted
seed data by \cite[Lemma~3.14]{JO26}.
The categorical interpretation of the \(F\)-invariant therefore
gives equal \(\de\)-values for the two upstairs pairs.
Since the \(F\)-invariant for \(x_{[2,5]}\) and \(x_{[7,9]}\)
is \(4\), Theorem~\ref{thm: categorical interpretation of F-invariant}
gives 
\[
\de\bigl(M^\ii[2,5],M^\ii[7,9]\bigr)
=d_2\!\sum_{\hfrake\in\{[3,8],[4,9]\}}
\de\bigl(M^\hii[\hfrake],M^\hii[12,15]\bigr)\\
=1\cdot(1+1)=2.
\]
Unlike the cuspidal modules in
\eqref{eq: cuspidal determinantial module}, both determinantial
modules in this example arise from boxes containing two occurrences
of their respective colors.   

\section{KR-modules over quiver Hecke algebras and symmetrized inverse
\texorpdfstring{\(t\)}{t}-Cartan matrices}
\label{sec: KR modules and inverse}

In this final section, we apply the folding formulas established
above to \emph{Kirillov--Reshetikhin modules over quiver Hecke algebras}.
We first introduce these modules as determinantial modules attached
to \(\bfi\)-boxes arising from adapted reduced sequences.
We then relate their \(\de\)-invariants to the symmetrized inverse
\(t\)-Cartan matrix.
In finite type, this relation, together with folding, yields
explicit formulas in all finite types and closed formulas in the
classical finite types.
Together with the formula in
Proposition~\ref{prop: KKOP directional box Lambda formula}, these
\(\de\)-formulas also permit the computation of both orders of the
\(\Lambda\)-invariant for boxes in the same \(Q\)-adapted reduced
sequence.

We begin by describing adapted reduced sequences in terms of height
functions and repetition quivers.
Although the usual Auslander--Reiten interpretation is restricted
to finite type, the construction below makes sense for every
symmetrizable Cartan datum whose Dynkin diagram admits a height
function.

\subsection{Height functions and
\texorpdfstring{\(Q\)}{Q}-adapted sequences}
\label{subsec: height functions and Q-adapted sequences}

We use coordinates on repetition quivers to define the
KR-module parameters and show that these parameters are
preserved under folding.

Keep the symmetrizable Cartan datum
\[
C=(c_{i,j})_{i,j\in I}, \qquad D=\diag(d_i\mid i\in I),
\]
and let \(\Dynkin\) be the Dynkin diagram of $\g$ with vertex set \(I\). 
Let $\{\alpha_i\mid i\in I\}$
be the set of simple roots, normalized by
$(\alpha_i,\alpha_i)_\g=2d_i$.

For \(i,j\in I\), we write $i\sim j$
if $c_{i,j}c_{j,i}>0$. Thus \(i\sim j\) records adjacency in \(\Dynkin\).

A \emph{Dynkin quiver} is a pair $Q=(\Dynkin,\xi)$,
where a \emph{height function}
\[
\xi:I\longrightarrow\Z
\]
satisfies $|\xi_i-\xi_j|=1$ whenever $i\sim j$.

The height function determines an orientation of the underlying
graph of \(\Dynkin\) by
\[
i\longrightarrow j \quad\Longleftrightarrow\quad \xi_i=\xi_j+1.
\]
This quiver orientation should not be confused with the arrows in
the usual Dynkin diagram that indicate relative root lengths.

\begin{remark} \label{rmk: not every Dynkin quiver allowed}
Not every Dynkin diagram admits a height function.  Indeed, such a
function exists if and only if the underlying graph of \(\Dynkin\)
is bipartite.  Thus every finite-type Dynkin diagram admits a height
function, whereas the affine Dynkin diagram of type
\(A_{2n}^{(1)}\), whose underlying graph is an odd cycle, does not.
\end{remark}

\smallskip

\emph{Throughout this section, we assume that \(\Dynkin\) is connected
and admits a height function \(\xi\), and fix the corresponding
Dynkin quiver \(Q=(\Dynkin,\xi)\).}
\smallskip

A vertex \(i\in I\) is called a \emph{sink} of \(Q\) if
\[
\xi_i=\xi_j-1 \qquad\text{for every }j\sim i.
\]
If \(i\) is a sink, let \(s_iQ\) be the Dynkin quiver obtained by
reversing all arrows incident with \(i\).  Its height function is
given by
\[
(s_i\xi)_j \seteq \xi_j+2\delta_{i,j} \qquad(j\in I).
\]

A sequence $\ii=(i_1,\ldots,i_l)$ is called \emph{\(Q\)-adapted} if
\[
i_k \text{ is a sink of } s_{i_{k-1}}\cdots s_{i_1}Q \qquad(1\leq k\leq l).
\]
If moreover \(\ii\in\Red(w)\), we call it a
\emph{\(Q\)-adapted reduced sequence of \(w\)}.

Fix a parity function $\epsilon:I\longrightarrow\{0,1\}$
such that
\[
\epsilon_i\ne\epsilon_j \qquad\text{whenever}\qquad
i\sim j.
\]
After adding a constant to \(\xi\)  if
necessary, we assume that
$\xi_i\equiv\epsilon_i\pmod2$ $(i\in I)$.

The \emph{repetition quiver} associated with \(\Dynkin\) is the
colored quiver \(\tDynkin\) with vertex set
\[
\tDynkin_0 \seteq \bigl\{ (i,p)\in I\times\Z \mid p\equiv\epsilon_i\pmod2 \bigr\}
\]
and arrow set
\[
\tDynkin_1 \seteq \bigl\{ (i,p)\longrightarrow(j,p-1) \mid i\sim j \bigr\}.
\]
Here every arrow is drawn as a single arrow.  The lacing data of
\(\Dynkin\) are retained by the colors of the vertices (see Example~\ref{ex: one-sided repetition quivers} below).

The height function \(\xi\) determines the full subquiver
$\tDynkin_{\geq\xi}$
of \(\tDynkin\) with vertex set
\[
\bigl(\tDynkin_{\geq\xi}\bigr)_0 \seteq \bigl\{ (i,p)\in\tDynkin_0 \mid p\geq\xi_i \bigr\}.
\]
Thus the vertices of color \(i\) are
\[
(i,\xi_i),\  (i,\xi_i+2),\  (i,\xi_i+4),\ldots.
\]
The boundary vertices
$\{(i,\xi_i)\mid i\in I\}$
form a copy of the Dynkin quiver \(Q\).

\smallskip
The same one-sided indexing is used by Fujita and Murakami
in their definition of $\Phi_\xi$ for bipartite generalized Cartan
matrices~\cite[Definition~2.13]{FM23general}.

\begin{example}\label{ex: one-sided repetition quivers}
We give three examples of the one-sided repetition quiver \(\tDynkin_{\geq\xi}\).

\begin{enumerate}[(i)]
\item Let $Q=(A_3,\xi)$ with $(\xi_1,\xi_2,\xi_3)=(1,2,3)$.
Thus the boundary orientation is $3\to2\to1$.
The one-sided repetition quiver is
\[
\raisebox{7mm}{
\scalebox{0.78}{\xymatrix@!C=4mm@R=3mm{
(i\setminus p)&1&2&3&4&5&6&7&8&9\\
1&\bullet&&\bullet\ar[dl]&&\bullet\ar[dl]&&
  \bullet\ar[dl]&&\bullet\ar[dl]\\
2&&\bullet\ar[ul]&&\bullet\ar[ul]\ar[dl]&&
  \bullet\ar[ul]\ar[dl]&&\bullet\ar[ul]\ar[dl]&\\
3&&&\bullet\ar[ul]&&\bullet\ar[ul]&&
  \bullet\ar[ul]&&\bullet\ar[ul]
}}}.
\]

\item Let $Q=(C_3,\xi)$ with $(\xi_1,\xi_2,\xi_3)=(1,2,3)$.
Again, the boundary orientation is $3\to2\to1$.
The one-sided repetition quiver is
\[
\raisebox{7mm}{
\scalebox{0.78}{\xymatrix@!C=4mm@R=3mm{
(i\setminus p)&1&2&3&4&5&6&7&8&9\\
1&\bullet&&\bullet\ar[dl]&&\bullet\ar[dl]&&
  \bullet\ar[dl]&&\bullet\ar[dl]\\
2&&\bullet\ar[ul]&&\bullet\ar[ul]\ar[dl]&&
  \bullet\ar[ul]\ar[dl]&&\bullet\ar[ul]\ar[dl]&\\
3&&&\circ\ar[ul]&&\circ\ar[ul]&&
  \circ\ar[ul]&&\circ\ar[ul]
}}}.
\]
Here the two vertex styles distinguish the two root lengths.  As in
the simply-laced case, only one arrow is drawn between adjacent
colored vertices.
\item 
Let $Q=(A_3^{(1)},\xi)$ with $(\xi_0,\xi_1,\xi_2,\xi_3)=(1,2,1,2)$.
Thus the boundary orientation is
\[
1\longrightarrow0\longleftarrow3,
\qquad
1\longrightarrow2\longleftarrow3.
\]
The one-sided repetition quiver is represented by the following
fundamental strip:
\[
\begin{tikzpicture}[
  x=0.72cm,
  y=0.72cm,
  baseline=(current bounding box.center),
  every node/.style={font=\small},
  vertex/.style={circle,fill=black,inner sep=1.45pt},
  arrow/.style={->,>=stealth,thin}
]

\node at (0,0.65) {$(i\backslash p)$};
\foreach \p in {1,...,9}{
  \node at (\p,0.65) {\p};
}

\node at (0,0)   {$0$};
\node at (0,-1)  {$1$};
\node at (0,-2)  {$2$};
\node at (0,-3)  {$3$};
\node at (0,-4)  {$0$};

\foreach \p in {1,3,5,7,9}{
  \node[vertex] (zt\p) at (\p,0) {};
  \node[vertex] (tw\p) at (\p,-2) {};
  \node[vertex] (zb\p) at (\p,-4) {};
}

\foreach \p in {2,4,6,8}{
  \node[vertex] (on\p) at (\p,-1) {};
  \node[vertex] (th\p) at (\p,-3) {};
}

\foreach \p/\pm/\pp in {2/1/3,4/3/5,6/5/7,8/7/9}{
  \draw[arrow] (on\p) -- (zt\pm);
  \draw[arrow] (zt\pp) -- (on\p);

  \draw[arrow] (on\p) -- (tw\pm);
  \draw[arrow] (tw\pp) -- (on\p);

  \draw[arrow] (th\p) -- (tw\pm);
  \draw[arrow] (tw\pp) -- (th\p);

  \draw[arrow] (th\p) -- (zb\pm);
  \draw[arrow] (zb\pp) -- (th\p);
}

\draw[<->,dashed,>=stealth]
  (9.65,-0.05)
  --
  node[right,align=left,font=\scriptsize]
  {identified}
  (9.65,-3.95);

\node at (9.25,0)  {$\cdots$};
\node at (9.25,-1) {$\cdots$};
\node at (9.25,-2) {$\cdots$};
\node at (9.25,-3) {$\cdots$};
\node at (9.25,-4) {$\cdots$};

\end{tikzpicture}
\]
Here the first and last rows are two copies of the color-\(0\)
boundary and are identified.  Thus the displayed fundamental strip
is naturally regarded as a quiver drawn on a cylinder.
\end{enumerate}
\end{example}
Let \(\Gamma\) be a finite full subquiver of
\(\tDynkin_{\geq\xi}\).  A sequence
\[
\ttR=(v_1,\ldots,v_r)
\]
containing every vertex of \(\Gamma\) exactly once is called a
\emph{compatible reading} of \(\Gamma\) if
\[
v_a\longrightarrow v_b \quad\Longrightarrow\quad b<a.
\]
Thus the target of every arrow is read before its source.

Write \(v_a=(i_a,p_a)\) for \(1\leq a\leq r\). For \(1\leq k\leq r\), set  
\[
m_i(k) \seteq \bigl|\{1\leq a<k\mid i_a=i\}\bigr|,
\qquad \ii(\ttR)\seteq(i_1,\ldots,i_r).
\]

After \(v_1,\ldots,v_{k-1}\) have been read, set  
\[
\partial_k \seteq \bigl\{ (i,\xi_i+2m_i(k)) \mid i\in I \bigr\}
\]
and for a vertex  
\[
(i,\xi_i+2m_i(k))\in\partial_k
\]
consider the condition 
\[
\xi_i+2m_i(k) = \xi_j+2m_j(k)-1
\qquad
\text{for every }j\sim i.
\]

We consider compatible readings \(\ttR\) such that, for every \(1\leq k\leq r\), 
\[
(i_k,p_k) = \bigl(i_k,\xi_{i_k}+2m_{i_k}(k)\bigr) \in\partial_k
\]
and the equality for adjacent colors holds with \(i=i_k\). 

\begin{lemma}\label{lem: boundary readings are Q-adapted}
Let $\ttR=((i_1,p_1),\ldots,(i_r,p_r))$ be a compatible reading of a finite
full subquiver of \(\tDynkin_{\geq\xi}\), satisfying the two conditions above.
Then 
\[
\ii(\ttR)\seteq(i_1,\ldots,i_r) \text{ is \(Q\)-adapted.}
\]
\end{lemma}

\begin{proof}
Reflection at a sink increases its height by \(2\) and leaves the
other heights unchanged.  Thus, just before the \(k\)-th vertex is
read, the height function is
\[
\xi^{(k-1)}=s_{i_{k-1}}\cdots s_{i_1}\xi,
\qquad
\xi_i^{(k-1)}=\xi_i+2m_i(k)\quad(i\in I).
\]
Thus \(i_k\) is a sink of \(s_{i_{k-1}}\cdots s_{i_1}Q\), so
\(\ii(\ttR)\) is \(Q\)-adapted. 
\end{proof}

Let
\[
\ttR=((i_1,p_1),\ldots,(i_l,p_l))
\]
be a finite compatible reading as in
Lemma~\ref{lem: boundary readings are Q-adapted}, and set 
\[
\ii\seteq\ii(\ttR)=(i_1,\ldots,i_l).
\]
We denote by $\Upsilon_Q^\ii$
the full subquiver of \(\tDynkin_{\geq\xi}\) whose vertex set is
\[
(\Upsilon_Q^\ii)_0 \seteq \{(i_k,p_k)\mid1\leq k\leq l\}.
\]
We call \(\Upsilon_Q^\ii\) the
\emph{\(Q\)-adapted quiver associated with \(\ii\)}.

The position of the \(k\)-th vertex is determined by \(Q\) and
\(\ii\).  Since \(Q\) is fixed, we write
\[
p_k^\ii \seteq \xi_{i_k} + 2\bigl| \{1\leq a<k\mid i_a=i_k\} \bigr| \qquad(1\leq k\leq l). 
\]
Thus the compatible reading associated with \(\ii\) is
\[
\bigl( (i_1,p_1^\ii),\ldots,(i_l,p_l^\ii) \bigr).
\]

We shall also use Dynkin quivers compatible with the folding datum
of Section~\ref{sec: folding cluster algebras}.  Let
\[
\hQ=(\hDynkin,\widehat\xi)
\]
be a Dynkin quiver for the simply-laced Cartan datum of \(\hg\).
We call \(\hQ\) \emph{\(\Gamma\)-foldable} if 
\[
\widehat\xi_{g(\im)} = \widehat\xi_\im
\qquad (g\in \Gamma,\ \im\in\hI).
\]
In this case, \(\widehat\xi\) induces a height function on the
quotient Dynkin diagram by
\[
\xi_i \seteq \widehat\xi_\im \qquad (\im\in I_i).
\]
We denote the quotient Dynkin quiver by
\[
Q\seteq\hQ/\Gamma=(\Dynkin,\xi).
\]

\begin{lemma}\label{lem: folding of adapted repetition quivers}
Let \(\hQ\) be a \(\Gamma\)-foldable simply-laced Dynkin quiver and set
$Q\seteq\hQ/\Gamma.$
Then $(\im,p)\longmapsto(\pi(\im),p)$
defines the natural folding of repetition quivers
\[
\widetilde{\hDynkin}_{\geq\widehat\xi} \longrightarrow \tDynkin_{\geq\xi}.
\]
Moreover, if $\ii=(i_1,\ldots,i_l)$
is \(Q\)-adapted, then its blockwise lift
$\hii$ is \(\hQ\)-adapted
(\cite[Lemma~8.1]{JO26}, stated there in the opposite source-adapted convention).  If \(\ii\) and \(\hii\) are reduced, the
preceding folding restricts to
\[
\Upsilon_{\hQ}^{\hii} \longrightarrow \Upsilon_Q^\ii.
\]
Under the canonical lift of an admissible chain, the obtained
repetition-quiver vertex fibers induce the box fibers used in
Section~\ref{sec: folding cluster algebras}.
\end{lemma}

\begin{lemma}\label{lem: folding preserves KR box parameters}
Assume that \(\ii\) is a \(Q\)-adapted reduced sequence.  Let
\(\frakC\) be an admissible chain associated with \(\ii\), and let
\(\hfrakC\) be its canonical lift.
\bnum
\item For $\frakc=[x,y]\in\calF(\frakC)$, set 
$k_\frakc \seteq \bigl| \{a\in[x,y]\mid i_a=i_\frakc\} \bigr|$ and $p_\frakc \seteq \dfrac{p_x^\ii+p_y^\ii}{2}$.
\item For $\hfrakc=[\hx,\hy] \in\Fib_{\frakC}(\frakc)$,
let $\im_{\hfrakc}\in\hI$ be its color, so that
$\pi(\im_{\hfrakc})=i_\frakc$, and set
\[k_{\hfrakc} \seteq \bigl| \{ \ha\in[\hx,\hy] \mid \hi_{\ha}=\im_{\hfrakc} \} \bigr| \qtq
p_{\hfrakc} \seteq \frac{ p_{\hx}^{\hii} + p_{\hy}^{\hii} }{2}.
\] 
\ee 
Then we have 
\[
\pi(\im_{\hfrakc})=i_\frakc, \quad k_{\hfrakc}=k_\frakc, \quad p_{\hfrakc}=p_\frakc  \quad \text{for every $\hfrakc\in\Fib_{\frakC}(\frakc)$.}
\]
\end{lemma}

\begin{proof}
Let \(\rho:[1,N_l]\longrightarrow[1,l]\) be defined by
\(\rho(\kappa(a,r))=a\).
By the construction of the canonical lift
\cite[Lemma~3.13]{JO26}, \(\rho\) induces an
order-preserving bijection
\[
\{\ha\in[\hx,\hy]\mid\hi_{\ha}=\im_{\hfrakc}\}
\xrightarrow{\sim}
\{a\in[x,y]\mid i_a=i_\frakc\}.
\]
Hence \(k_{\hfrakc}=k_\frakc\),
\(\rho(\hx)=x\), and \(\rho(\hy)=y\).
Since the repetition-quiver folding preserves the second
coordinate,
\[
p_{\hfrakc}
=\frac{p_{\hx}^{\hii}+p_{\hy}^{\hii}}2
=\frac{p_x^\ii+p_y^\ii}2
=p_\frakc.
\]
The equality \(\pi(\im_{\hfrakc})=i_\frakc\) holds by construction.
\end{proof}

Suppose now that $\ii=(i_1,\ldots,i_l)\in\Red(w)$ is
\(Q\)-adapted, where \(w\in W_\g\). Recall from \eqref{eq: root sequence} that
\[ 
\be_k^\ii = s_{i_1}\cdots s_{i_{k-1}}(\alpha_{i_k}).
\]
Since \(\ii\) is reduced, one has
\[
\be_k^\ii\in\Delta_+ \qquad(1\leq k\leq l) 
\qtq 
\{\be_k^\ii\mid1\leq k\leq l\}
=
\Delta_+\cap w(\Delta_-).
\]
We label the vertices of \(\Upsilon_Q^\ii\) by
$(i_k,p_k)
\longleftrightarrow
\be_k^\ii$.

Different compatible readings of the same quiver
\(\Upsilon_Q^\ii\) are related by commutation moves, that is, by
interchanging consecutive vertices with distinct, nonadjacent colors in $\Dynkin$. Since the
corresponding simple reflections commute, the positive-root
label of each vertex is independent of the chosen compatible
reading.

When \(\g\) is of finite type and \(w=w_0\), we write
\[
\Upsilon_Q\seteq\Upsilon_Q^\ii.
\]
Then \(\Upsilon_Q\) is the usual Auslander--Reiten quiver associated
with \(Q\), its vertex set is labeled by \(\Delta_+\), and
\[
|(\Upsilon_Q)_0|=|\Delta_+|.
\]
For a general symmetrizable Cartan datum, we use the finite
\(Q\)-adapted quiver \(\Upsilon_Q^\ii\) associated with a chosen
\(Q\)-adapted sequence, imposing reducedness only when a
positive-root labeling is required.

\subsection{KR-modules over quiver Hecke algebras}
\label{subsec: KR modules}

Let $\ii=(i_1,\ldots,i_l)\in\Red(w)$
be a \(Q\)-adapted reduced sequence, and let
\[
(i_k,p_k)\longleftrightarrow\be_k^\ii
\qquad(1\leq k\leq l)
\]
be the corresponding labeling of \(\Upsilon_Q^\ii\).
Recall from \eqref{eq: cuspidal determinantial module} that the
cuspidal module associated with the \(k\)-th position is \(S_k^\ii\).
For a vertex \((i_k,p_k)\in\Upsilon_Q^\ii\), we write
\[
S_Q(i_k,p_k)\seteq S_k^\ii.
\]

Let \(r\in\Z_{\geq0}\), and suppose that the vertices
\[
(i,p-2r),\ (i,p-2r+2),\ldots,(i,p-2),\ (i,p)
\]
occur in \(\Upsilon_Q^\ii\). Then
\[
\bigl(
S_Q(i,p),S_Q(i,p-2),\ldots,S_Q(i,p-2r)
\bigr)
\]
is a normal sequence (\cite[Proposition 4.3]{KKOP23F}).  In particular, 
\[
\text{ the convolution product } \  S_Q(i,p)\conv S_Q(i,p-2)\conv\cdots\conv S_Q(i,p-2r) \  
\text{ has a simple head.}
\]

We define the corresponding \emph{KR-module over the quiver Hecke
algebra} by
\begin{align} \label{eq: KR-module over quiver Hecke}
\Vkm{i^{r+1}}_{p-r}
\seteq
\head\bigl(
S_Q(i,p)\conv S_Q(i,p-2)\conv\cdots
\conv S_Q(i,p-2r)
\bigr).
\end{align}
Here \(p-r\) is the center of the parameters
\[
p-2r,\ p-2r+2,\ldots,p-2,\ p.
\]
Equivalently, for \(k\geq1\), the centered notation
$\Vkm{i^{k}}_{a}$ is characterized by
\[
\Vkm{i^{k}}_{a} = \head\bigl( S_Q(i,a+k-1)\conv S_Q(i,a+k-3)\conv\cdots\conv S_Q(i,a-k+1) \bigr),
\]
whenever all the corresponding vertices occur in
\(\Upsilon_Q^\ii\).

More precisely, suppose that \((i,p-2r)\) and \((i,p)\) are the
\(x\)-th and \(y\)-th vertices, respectively, in the compatible
reading defining \(\ii\).  Then
\[
i_x=i_y=i
\]
and the positions of color \(i\) in \([x,y]\) correspond exactly to
\[
(i,p-2r),\ (i,p-2r+2),\ldots,(i,p).
\]
Consequently,
\[
\begin{aligned}
\Vkm{i^{r+1}}_{p-r}
\simeq M^\ii[x,y]=
M\bigl(w_{\leq y}^\ii\varpi_i,w_{<x}^\ii\varpi_i \bigr).
\end{aligned}
\]
Thus the KR-modules defined above are precisely the determinantial
modules attached to the \(\bfi\)-boxes associated with the
\(Q\)-adapted reduced sequence \(\ii\).

\subsection{Symmetrized inverse \texorpdfstring{\(t\)}{t}-Cartan matrices}
\label{subsec: inverse t-Cartan}

Let
\[
A=(a_{i,j})_{i,j\in I},
\qquad
a_{i,j}\seteq2\delta_{i,j}-c_{i,j}.
\]
Thus
\[
a_{i,i}=0 \qtq a_{i,j}=-c_{i,j}\in\Z_{\geq0}
\qquad(i\ne j).
\]
The \emph{\(t\)-quantized Cartan matrix} is
\[
C(t)
\seteq
(t+t^{-1})\operatorname{Id}-A
=
\bigl(c_{i,j}(t)\bigr)_{i,j\in I} 
\Longleftrightarrow
c_{i,j}(t)
=
\begin{cases}
t+t^{-1},&i=j,\\
c_{i,j},&i\ne j.
\end{cases}
\]

Since
\[
C(t) = t^{-1}
\bigl(\operatorname{Id}-tA+t^2\operatorname{Id}\bigr),
\]
the inverse \(t\)-quantized Cartan matrix has the \(t\)-adic expansion
\[
C(t)^{-1}=
t \bigl(\operatorname{Id}-tA+t^2\operatorname{Id}\bigr)^{-1}
\in t\operatorname{Mat}_I(\Z[[t]]).
\]
We refer to
\[
\tC(t) \seteq D C(t)^{-1} = \bigl(\tC_{i,j}(t)\bigr)_{i,j\in I}
\]
as the \emph{symmetrized inverse \(t\)-quantized Cartan matrix}, or
simply the \emph{symmetrized inverse \(t\)-Cartan matrix}.
Note that \(\tC(t)\) is symmetric.

For \(i,j\in I\), write
\begin{align} \label{eq: coeff}
\tC_{i,j}(t) \seteq \sum_{u\in\Z}\tc_{i,j}(u)t^u.
\end{align}
Hence
$\tc_{i,j}(u) = \tc_{j,i}(u)$ for $i,j\in I$ and $u\in\Z$.

For \(i,j\in I\), let \(d(i,j)\) denote the distance between \(i\)
and \(j\) in $\Dynkin$.

\begin{proposition}[{\cite{Fuj22a,FO21,KO23,FM21,FM23general}}]
\label{prop: coefficients of inverse t-Cartan}
For \(i,j\in I\) and \(u\in\Z\), the coefficients $\tc_{i,j}(u)$ in~\eqref{eq: coeff}
satisfy the following properties.
\bnum
\item \label{it: 0-range}
We have $\tc_{i,j}(u)=0$ for $u\leq d(i,j)$ and, in particular,
$\tc_{i,j}(1)=d_i\delta_{i,j}$.

\item \label{it: parity} If  $\tc_{i,j}(u)\ne0$
then $u+\epsilon_i+\epsilon_j\equiv1\pmod2$.

\item \label{it: positivity} If \(C\) is indecomposable and of infinite type, then
$\tc_{i,j}(u)\in\Z_{\geq0}$
for $u\in\Z$, 
and hence $\tC_{i,j}(t) \in t\Z_{\geq0}[[t]]$.
\ee
\end{proposition}

\begin{proof}
The finite-type case of~\eqref{it: 0-range}, as well as~\eqref{it: parity} and~\eqref{it: positivity}, is covered by
the cited references.  For completeness, we give a direct proof of~\eqref{it: 0-range}, valid in arbitrary type.  Write
\[
 (I-tA+t^2I)^{-1}=\sum_{r\geq0}Q_rt^r.
\]
Comparing coefficients gives
\[
 Q_0=I,\qquad Q_1=A,\qquad
 Q_r=AQ_{r-1}-Q_{r-2}\quad(r\geq2).
\]
It follows inductively that each \(Q_r\) is a polynomial in \(A\) of
degree at most \(r\).  Since
\[
 (A^k)_{i,j}=0\qquad\text{whenever }k<d(i,j),
\]
we have \((Q_r)_{i,j}=0\) for \(r<d(i,j)\).  Now
\[
 \tC(t)=Dt\sum_{r\geq0}Q_rt^r.
\]
Thus there are no terms of nonpositive degree, and for
\(1\leq u\leq d(i,j)\) the coefficient of \(t^u\) vanishes because
\(u-1<d(i,j)\).  Hence
\[
 \tc_{i,j}(u)=0\qquad\text{for }u\leq d(i,j).
\]
Moreover, \(\tC(t)=tD+O(t^2)\), and therefore
\[
 \tc_{i,j}(1)=d_i\delta_{i,j}. \qedhere
\]
\end{proof}

We next record the compatibility of
\(\tC(t)\) with folding.  Keep the folding datum of
Section~\ref{sec: folding cluster algebras}.  
Let
\[
\tC^\hg(t) = \bigl( \tC_{\im,\jm}^\hg(t) \bigr)_{\im,\jm\in\hI} 
\qtq 
\tC^\g(t) = \bigl( \tC_{i,j}^\g(t) \bigr)_{i,j\in I} 
\]
be the corresponding symmetrized inverse \(t\)-Cartan matrices.
Since the upstairs datum is simply-laced, we have 
\[
\tC^\hg(t) = C^\hg(t)^{-1}.
\]

Let $R:\Q(t)^I\longrightarrow\Q(t)^{\hI}$
be the lifting map onto the fiber-constant subspace, defined by
\[
(R\bfv)_\im \seteq v_{\pi(\im)}.
\]
With respect to the standard bases, its matrix is given by
\[
R_{\im,j} = \delta_{\pi(\im),j}  \qquad (\im\in\hI,\ j\in I).
\]
The quotient definition
\(c_{j,i}=\sum_{\jm\in I_j}\hc_{\jm,\im}\) for \(\im\in I_i\),
together with the symmetry of \(\hC\), gives 
\begin{equation}\label{eq: folding of t-Cartan matrices}
C^\hg(t)R = RC^\g(t)^T.
\end{equation}

\begin{example}\label{ex: A5 C3 t-Cartan folding}
Consider the standard folding
$A_5\longrightarrow C_3$
induced by the involution
$1\longleftrightarrow5$, $2\longleftrightarrow4$, $3\longleftrightarrow3$.
The vertex fibers and the downstairs symmetrizer are
\[
I_1=\{1,5\}, \qquad I_2=\{2,4\}, \qquad I_3=\{3\}, \qquad D^{C_3}=\diag(1,1,2).
\]
Set $a\seteq t+t^{-1}$. Then
$C^{A_5}(t)$, $C^{C_3}(t)$ and the lifting map $R$ are given as follows:
\[
C^{A_5}(t) =
{\scriptsize\begin{pmatrix}
a&-1&0&0&0\\
-1&a&-1&0&0\\
0&-1&a&-1&0\\
0&0&-1&a&-1\\
0&0&0&-1&a
\end{pmatrix}},
\quad
C^{C_3}(t)
=
{\scriptsize\begin{pmatrix}
a&-1&0\\
-1&a&-2\\
0&-1&a
\end{pmatrix}},
\quad 
R =
{\scriptsize\begin{pmatrix}
1&0&0\\
0&1&0\\
0&0&1\\
0&1&0\\
1&0&0
\end{pmatrix}},
\]
satisfying 
\[
C^{A_5}(t)R =
{\scriptsize\begin{pmatrix}
a&-1&0\\
-1&a&-1\\
0&-2&a\\
-1&a&-1\\
a&-1&0
\end{pmatrix}}
= R\,C^{C_3}(t)^T,
\]
which illustrates \eqref{eq: folding of t-Cartan matrices}.
\end{example}

\begin{theorem}\label{thm: folding of inverse t-Cartan matrices}
Let \(i,j\in I\), and choose $\im\in I_i$ and $\jm\in I_j$. 
Then
\begin{equation}\label{eq: folding inverse t-Cartan first fiber}
\tC_{i,j}^\g(t) = d_i \sum_{\im'\in I_i} \tC_{\im',\jm}^\hg(t)
= \frac{m}{|\Gamma|} \sum_{g\in \Gamma} \tC_{g(\im),\jm}^\hg(t).
\end{equation}
These expressions are independent of the chosen lifts.  By symmetry,
one also has
\begin{equation}\label{eq: folding inverse t-Cartan second fiber}
\tC_{i,j}^\g(t) = d_j \sum_{\jm'\in I_j} \tC_{\im,\jm'}^\hg(t)
= \frac{m}{|\Gamma|} \sum_{g\in \Gamma} \tC_{\im,g(\jm)}^\hg(t).
\end{equation}
\end{theorem}

\begin{proof}
Transposing \eqref{eq: folding of t-Cartan matrices} and using the
symmetry of \(C^\hg(t)\) gives
\(R^TC^\hg(t)=C^\g(t)R^T\).  Since both Cartan matrices are
invertible as \(t\)-adic formal matrices,
\[
R^TC^\hg(t)^{-1}=C^\g(t)^{-1}R^T.
\]
Using \(\tC^\hg(t)=C^\hg(t)^{-1}\) and comparing the
\((i,\jm)\)-entries for \(\jm\in I_j\), we obtain
\[
\sum_{\im'\in I_i}\tC_{\im',\jm}^\hg(t)
=\bigl(C^\g(t)^{-1}\bigr)_{i,j}.
\]
Multiplication by \(d_i\) gives the first fiber formula in
\eqref{eq: folding inverse t-Cartan first fiber}, independently of
the chosen lift \(\jm\).
For any \(\im\in I_i\), orbit--stabilizer gives
\[
\sum_{g\in\Gamma}\tC_{g(\im),\jm}^\hg(t)
=\frac{|\Gamma|}{e_i}\sum_{\im'\in I_i}\tC_{\im',\jm}^\hg(t).
\]
Together with \(d_ie_i=m\), this gives the group-sum factor
\(m/|\Gamma|\), independently of \(\im\).
The second fiber formula
\eqref{eq: folding inverse t-Cartan second fiber} follows by
symmetry of \(\tC^\g(t)\) and \(\tC^\hg(t)\).
\end{proof}

The folding identity also holds coefficientwise:
\[
\tc_{i,j}^\g(u)
= d_i \sum_{\im'\in I_i} \tc_{\im',\jm}^\hg(u)
= \frac{m}{|\Gamma|} \sum_{g\in \Gamma} \tc_{g(\im),\jm}^\hg(u)
\qquad(u\geq1).
\]

\begin{remark}
Theorem~\ref{thm: folding of inverse t-Cartan matrices} applies to
admissible foldings of arbitrary Cartan type, including indefinite type.
\end{remark}

\begin{example}\label{ex: affine D4 to B3 inverse t-Cartan folding}
Consider the standard folding
$D_4^{(1)}\longrightarrow B_3^{(1)}$.
We label the vertices of \(D_4^{(1)}\) so that \(2\) is the central
vertex and $0,\ 1,\ 3,\ 4$
are its four adjacent vertices.  Let \(\Gamma\) be generated by the
diagram automorphism interchanging \(3\) and \(4\).  The vertex
fibers are
$I_0=\{0\}, I_1=\{1\}, I_2=\{2\}, I_3=\{3,4\}$,
and the symmetrizer of \(B_3^{(1)}\) is $D=\diag(2,2,2,1).$

The first terms of the relevant entries upstairs are
\[
\begin{aligned}
\tC_{0,0}^{D_4^{(1)}}(t)
= t+2t^5+t^7+O(t^9),\quad
\tC_{0,1}^{D_4^{(1)}}(t)
= t^3+t^5+2t^7+O(t^9).
\end{aligned}
\]
The corresponding entries downstairs are
\[
\begin{aligned}
\tC_{0,0}^{B_3^{(1)}}(t)
=
2t+4t^5+2t^7+O(t^9),\quad
\tC_{0,1}^{B_3^{(1)}}(t)
=
2t^3+2t^5+4t^7+O(t^9).
\end{aligned}
\]
Thus Theorem~\ref{thm: folding of inverse t-Cartan matrices} gives
\[
\begin{aligned}
\tC_{0,0}^{B_3^{(1)}}(t)
=
d_0
\sum_{\im\in I_0}
\tC_{\im,0}^{D_4^{(1)}}(t)
=
2\tC_{0,0}^{D_4^{(1)}}(t),\quad
\tC_{0,1}^{B_3^{(1)}}(t)
=
d_0
\sum_{\im\in I_0}
\tC_{\im,1}^{D_4^{(1)}}(t)
=
2\tC_{0,1}^{D_4^{(1)}}(t).
\end{aligned}
\]

The nontrivial vertex fiber gives a sum of two upstairs entries.  Indeed, 
\[
\begin{aligned}
\tC_{3,0}^{D_4^{(1)}}(t)
=
t^3+t^5+2t^7+O(t^9),\quad
\tC_{4,0}^{D_4^{(1)}}(t)
=
t^3+t^5+2t^7+O(t^9),
\end{aligned}
\]
and hence
\[
\tC_{3,0}^{B_3^{(1)}}(t)
= \sum_{\im\in I_3} \tC_{\im,0}^{D_4^{(1)}}(t)= \tC_{3,0}^{D_4^{(1)}}(t) + \tC_{4,0}^{D_4^{(1)}}(t)
= 2t^3+2t^5+4t^7+O(t^9). 
\]
This illustrates both the factor \(d_0=2\) on a singleton fiber
and the sum over a nontrivial orbit. 
\end{example}

We next turn to the explicit formulas for the
\(\de\)-invariants between \(Q\)-adapted KR-modules of finite types and their compatibility with folding.

\subsection{Explicit finite-type formulas and folding}
\label{subsec: finite-type formulas and folding}

Assume in this subsection that \(C\) is of finite type.  Let
$Q=(\Dynkin,\xi)$ be a Dynkin quiver, and let
$\ii=(i_1,\ldots,i_N)\in\Red(w_0)$ $(N=|\Delta_+|)$
be a \(Q\)-adapted reduced sequence.  Then the associated quiver
\[
\Upsilon_Q=\Upsilon_Q^\ii
\]
is the usual Auslander--Reiten quiver of \(Q\).  Its vertices are
labeled by the positive roots through
\[
(i_k,p_k)\longleftrightarrow\be_k^\ii,
\]
where \(\be_k^\ii\) is given by \eqref{eq: root sequence},
and the corresponding cuspidal module is
\[
S_Q(i_k,p_k)=S_k^\ii.
\]

The \(\de\)-values between these cuspidal modules are completely
determined by the symmetrized inverse \(t\)-Cartan matrix.

\begin{theorem}[{\cite[Theorem~8.2]{KO24}}]
\label{thm: finite-type cuspidal d-values}
For any $(i,p),(j,q)\in(\Upsilon_Q)_0$,
one has
\[
\de\bigl(S_Q(i,p),S_Q(j,q)\bigr) = \tc_{i,j}(|p-q|-1).
\]
\end{theorem}

The formula depends only on the colors \(i,j\) and the relative
position \(|p-q|\), and hence applies directly to the sink-adapted
convention used here.

Although our main objects are modules over quiver Hecke algebras,
in simply-laced finite type $\g$ they are related to modules over quantum affine algebras $U_q'(\g^{(1)})$
through generalized quantum affine Schur--Weyl duality.
For a Dynkin quiver \(Q\), the construction in
\cite{KKK2,KKK18,KKKOIV} gives an exact monoidal functor
\[
\frakF_Q: R\text{-gmod} \longrightarrow \scrC_Q,
\]
where \(\scrC_Q\) is the heart subcategory associated with \(Q\)
inside the category of finite-dimensional integrable 
\(U_q'(\g^{(1)})\)-modules.  It sends the cuspidal module attached
to a vertex \((i,p)\in(\Upsilon_Q)_0\) to the corresponding
fundamental module:
\[
\frakF_Q\bigl(S_Q(i,p)\bigr)
\simeq
V(\varpi_i)_{(-q)^p}.
\]
For quantum affine modules, \(\Lambda\) and \(\de\) denote the
integer invariants defined from the renormalizing coefficient of the
universal \(R\)-matrix in \cite[Definition~2.10]{KKOP23P}.
Generalized quantum affine Schur--Weyl duality preserves these invariants:
\[
\Lambda(M,N)
=
\Lambda\bigl(\frakF_Q(M),\frakF_Q(N)\bigr),
\qquad
\de(M,N)
=
\de\bigl(\frakF_Q(M),\frakF_Q(N)\bigr)
\]
for simple modules \(M,N\in R\text{-gmod}\)
(\cite[Theorem~4.12]{KKOP23P}).   

In particular, the quiver Hecke KR-module associated with the
consecutive vertices
\[
(i,p),\ (i,p-2),\ldots,(i,p-2r)
\]
is sent to the KR-module over $U_q'(\g^{(1)})$
\[
\Vkm{i^{r+1}}_{p-r}
\]
(see \cite{OS25} for the centered notation).  Hence, whenever all the
vertices defining the corresponding KR-modules lie in \(\Upsilon_Q\),
the quiver Hecke and quantum affine \(\de\)-invariants agree.  

The quantum affine formulation also allows arbitrary centered
parameters, without requiring all the vertices defining the
corresponding KR-modules to lie in \(\Upsilon_Q\).
The denominator formulas for KR-modules and their relation with the
symmetrized inverse \(t\)-Cartan matrix are established in
\cite{FO21,OS25,CFM26}.  We first state this formula
for arbitrary centered quantum affine KR-modules of type \(ADE\).
Since \(\frakF_Q\) preserves \(\de\), the same formula holds
for quiver Hecke KR-modules whenever all their defining
vertices lie in \(\Upsilon_Q\).

Let \(\g\) be of finite type \(A\), \(D\), or \(E\), and let \(h\)
be its Coxeter number.   For \(i\in I\), \(k\in\Z_{>0}\), and \(p\in\Z\), write
$\Vkm{i^k}_p$ for the quantum affine KR-module over
$U_q'(\g^{(1)})$ with the centered spectral parameters
\[
\Vkm{i^k}_p
\seteq \head\bigl(
V(\varpi_i)_{(-q)^{p+k-1}}\otimes
V(\varpi_i)_{(-q)^{p+k-3}}\otimes\cdots\otimes
V(\varpi_i)_{(-q)^{p-k+1}}\bigr).
\] 
 
For finite-dimensional simple modules \(U\) and \(V\) over
$U_q'(\g^{(1)})$, let \(V_z\) be the spectral shift of \(V\) by
the indeterminate \(z\). Following \cite[(5.3)]{CFM26}, let
\(\frako(U,V)\) denote the pole order at \(z=1\) of the normalized
\(R\)-matrix from \(U\otimes V_z\) to \(V_z\otimes U\), and set 
\[
\de(U,V) \seteq \frako(U,V)+\frako(V,U).
\]
By~\cite[Proposition~3.16]{KKOP20}, it is expressed in terms of the
orders of zeros of the normalized \(R\)-matrix denominators.

Finally, set
\[
\tC_{i,j}^{[h]}(t)
\seteq
\sum_{v=1}^{h}\tc_{i,j}(v)t^v,
\qquad
[m]_t
\seteq
t^{m-1}+t^{m-3}+\cdots+t^{1-m}.
\]

\begin{theorem}[{\cite{OS25} and \cite[Theorem~6.27]{CFM26}}]
\label{thm: ADE quantum affine KR d-values}
Assume that \(\g\) is of simply-laced finite type \(A\), \(D\), or
\(E\), and let \(h\) be its Coxeter number.
Let \(i,j\in I\), \(k,\ell\in\Z_{>0}\), and \(p,s\in\Z\).  Set
\[
K\seteq\max(k,\ell), \qquad L\seteq\min(k,\ell), \qquad u\seteq|p-s|.
\]
Then we have
\begin{equation}\label{eq: ADE quantum affine KR d-value}
\de\bigl(\Vkm{i^k}_p, \Vkm{j^\ell}_s \bigr)  =[t^u]\! \left( t^K[L]_t\tC_{i,j}^{[h]}(t) \right). 
\end{equation}
\end{theorem}

For \(v\in\Z\), set 
\[
\tc_{i,j}^{[h]}(v) \seteq
\begin{cases}
\tc_{i,j}(v),&1\leq v\leq h,\\
0,&\text{otherwise}.
\end{cases}
\]

\begin{corollary}\label{cor: ADE KR d-value coefficient sum}
With the notation of
Theorem~\ref{thm: ADE quantum affine KR d-values}, one has
\begin{equation}\label{eq: ADE KR d-value coefficient sum}
\de\bigl(\Vkm{i^k}_p,\Vkm{j^\ell}_s\bigr)
= \sum_{a=0}^{L-1} \tc_{i,j}^{[h]} \bigl( u-|k-\ell|-2a-1 \bigr). 
\end{equation}
\end{corollary}

\begin{proof}
Since
$t^K[L]_t = \displaystyle\sum_{r=0}^{L-1}t^{K+L-1-2r}$,
the coefficient of \(t^u\) in the right-hand side of
\eqref{eq: ADE quantum affine KR d-value} is
\[
\sum_{r=0}^{L-1} \tc_{i,j}^{[h]} \bigl( u-K-L+1+2r \bigr). 
\]
Replacing \(r\) by \(L-1-a\) and using $K-L=|k-\ell|$
gives the desired formula.
\end{proof}

We now make
Corollary~\ref{cor: ADE KR d-value coefficient sum} explicit in
type \(A\).  Assume that \(\g\) is of type \(A_n\), and set
\[
\brho^{(n)}_{i,j} \seteq \min(i,j,n+1-i,n+1-j).
\]
The coefficients of \(\tC_{i,j}^{[h]}(t)\) are 
\begin{equation}\label{eq: type A inverse t-Cartan coefficients}
\tc_{i,j}^{[h]}(v)= 
\begin{cases}
1, & v=|i-j|+2r-1 \text{ for some } 1\leq r\leq\brho^{(n)}_{i,j}, \\
0,&\text{otherwise}
\end{cases}
\end{equation}
(see \cite[Theorem~4.12(1)]{KO23} for more details).

\begin{corollary}\label{cor: type A closed d-value formula}
Let \(\g\) be of type \(A_n\).  For
\(i,j\in I\), \(k,\ell\in\Z_{>0}\), and \(p,s\in\Z\), set
\[
L\seteq\min(k,\ell), \qquad u\seteq|p-s|, \qquad \nu \seteq \frac{u-|i-j|-|k-\ell|}{2}. 
\]
Then we have
\begin{equation}\label{eq: type A closed d-value formula}
\de\bigl( \Vkm{i^k}_p, \Vkm{j^\ell}_s \bigr)
=
\begin{cases}
\displaystyle  \bigl [ \min\bigl( \nu,\brho^{(n)}_{i,j},L,\brho^{(n)}_{i,j}+L-\nu \bigr) \bigr]_+, &\nu\in\Z, \\[2mm]
0,&\nu\notin\Z. 
\end{cases}
\end{equation}
\end{corollary}

\begin{proof}
By \eqref{eq: type A inverse t-Cartan coefficients} and
Corollary~\ref{cor: ADE KR d-value coefficient sum}, the desired
\(\de\)-value is the number of pairs
\[
(r,a)\in\Z_{\ge1} \times \Z_{\ge0}
\]
satisfying
\[
1\leq r\leq\brho^{(n)}_{i,j}, \qquad 0\leq a\leq L-1, \qquad r+a=\nu.
\]
It vanishes unless \(\nu\in\Z\).  If \(\nu\in\Z\), the possible
values of \(r\) form the interval
\[
\max(1,\nu-L+1) \leq r\leq \min(\brho^{(n)}_{i,j},\nu).
\]
Its cardinality is
\[
\bigl[\min(\brho^{(n)}_{i,j},\nu) - \max(1,\nu-L+1)+1]_+,
\]
which is equal to the right-hand side of
\eqref{eq: type A closed d-value formula}.
\end{proof}

\begin{example}\label{ex: type A closed d-value}
Let \(\g\) be of type \(A_5\), and consider $\Vkm{2^2}_p$ and $\Vkm{3^4}_s$.
Here we have
$\brho^{(5)}_{2,3}=2$ and $L=2$. 
Consequently,
\[
\de\bigl( \Vkm{2^2}_p, \Vkm{3^4}_s \bigr)
=
\begin{cases}
1,&|p-s|=5,\text{ or } 9,\\
2,&|p-s|=7,\\
0,&\text{otherwise}.
\end{cases}
\]
\end{example}

We next give the corresponding closed formulas in type \(D\).
We use the standard numbering in which \(n-1\) and \(n\) are the
two spin nodes of \(D_n\).

For \(r,m\in\Z_{>0}\), \(v\in\Q\), and \(d\in\{1,2\}\), define
\begin{equation}\label{eq: closed T}
\sfT^{(d)}_{r,m}(v)
\seteq
\bigl| \bigl\{ (a,b)\in\Z_{\ge1} \times \Z_{\ge0} 
\mid 1\leq a\leq r,\  0\leq b\leq m-1,\  da+b=v \bigr\} \bigr|.
\end{equation}
Equivalently,
\[
\sfT^{(1)}_{r,m}(v)
\seteq
\begin{cases}
\displaystyle \big[ \min(v,r,m,r+m-v) \big]_+, &v\in\Z, \\[1mm]
0,&v\notin\Z,
\end{cases}
\]
and
\[
\sfT^{(2)}_{r,m}(v)
\seteq
\begin{cases}
\displaystyle \left[ \min\left(r,\left\lfloor\frac{v}{2}\right\rfloor\right) -
\max\left(1,\left\lceil\frac{v-m+1}{2}\right\rceil\right) +1 \right]_+, &v\in\Z, \\[3mm]
0,&v\notin\Z.
\end{cases}
\]

\begin{corollary}\label{cor: type D closed d-value formulas}
Let \(\g\) be of type \(D_n\).  For
\(i,j\in I\), \(k,\ell\in\Z_{>0}\), and \(p,s\in\Z\), set
\[
L\seteq\min(k,\ell), \qquad u\seteq|p-s|, \qquad \kappa\seteq|k-\ell|.
\]
Then the following formulas hold.
\bnum
\item Suppose that \(1\leq i,j\leq n-2\).  Set
\[
r\seteq\min(i,j), \qquad \nu_1 \seteq \frac{u-\kappa-|i-j|}{2}, \qquad \nu_2 \seteq \frac{u-\kappa-(2n-i-j-2)}{2}.
\]
Then we have
\[
\de\bigl( \Vkm{i^k}_p, \Vkm{j^\ell}_s \bigr)  = \sfT^{(1)}_{r,L}(\nu_1) + \sfT^{(1)}_{r,L}(\nu_2).
\]

\item Suppose that exactly one of \(i,j\) is a spin node.  Let
\[
c\seteq\min(i,j,n-2),
\]
so that \(c\) is the non-spin node, and set
\[
\nu \seteq \frac{u-\kappa-(n-c-1)}{2}.
\]
Then we have
\[
\de\bigl( \Vkm{i^k}_p, \Vkm{j^\ell}_s \bigr) = \sfT^{(1)}_{c,L}(\nu). 
\]

\item Suppose that \(i=j\in\{n-1,n\}\).  Set
\[
r_{\mathrm{sp}} \seteq \left\lfloor\frac{n}{2}\right\rfloor \qtq \nu \seteq \frac{u-\kappa+2}{2}.
\]
Then we have
\[
\de\bigl( \Vkm{i^k}_p, \Vkm{i^\ell}_s \bigr) = \sfT^{(2)}_{r_{\mathrm{sp}},L}(\nu).
\]

\item Suppose that \(\{i,j\}=\{n-1,n\}\).  Set
\[
r_{\mathrm{op}} \seteq \left\lfloor\frac{n-1}{2}\right\rfloor
\qtq
\nu \seteq \frac{u-\kappa}{2}.
\]
Then we have
\[
\de\bigl( \Vkm{i^k}_p, \Vkm{j^\ell}_s \bigr) = \sfT^{(2)}_{r_{\mathrm{op}},L}(\nu).
\]
\ee
\end{corollary}

\begin{proof}
By \cite[Theorem~4.12(2)]{KO23}, the coefficients of \(\tC_{i,j}^{[h]}(t)\) in type \(D_n\)
are given as follows. Set
\[
\bar i\seteq\min(i,n-1)
\qtq
\bar j\seteq\min(j,n-1).
\]
Then
\[
\sum_v \tc_{i,j}^{[h]}(v)t^v
=
\begin{cases}
\displaystyle
\sum_{a=1}^{\min(\bar i,\bar j)}
\left(
t^{|\bar i-\bar j|+2a-1}
+
\delta\bigl(\max(\bar i,\bar j)<n-1\bigr)
t^{2n-\bar i-\bar j+2a-3}
\right),
& \min(\bar i,\bar j)<n-1,
\\[4mm]
\displaystyle
\sum_{a=1}^{\left\lfloor
  (n-1+\delta_{i,j})/2
\right\rfloor}
t^{4a-1-2\delta_{i,j}},
& \min(\bar i,\bar j)\geq n-1.
\end{cases}
\]
Substituting these coefficients into
\eqref{eq: ADE KR d-value coefficient sum}
and using~\eqref{eq: closed T} gives the four formulas.
\end{proof}

 \begin{remark}\label{rem: exceptional KR d-values}
Theorem~\ref{thm: ADE quantum affine KR d-values} also gives
explicit finite formulas in types \(E_6,E_7\), and \(E_8\).
We do not list all exceptional cases, since they are obtained
directly from the finitely many coefficients  
\[
\tc_{i,j}^{[h]}(v)
\qquad
(i,j\in I,\ 1\leq v\leq h).
\]

The higher root multiplicities in the exceptional denominator
formulas are understood according to the formulas in~\cite{OhSc19,OS19c}.  For example, for the normalized \(R\)-matrix denominator
\(d_{4,4}(z)\) between two color-\(4\) fundamental quantum affine
modules of type \(E_6\), with the numbering used in \cite[Proposition~2.1]{OS19c}, one has
\[
d_{4,4}(z)
=
(z-q^2)(z-q^4)^2(z-q^6)^3
(z-q^8)^3(z-q^{10})^2(z-q^{12}).
\]
In this example, the multiplicity of each root \(q^a\) in
\(d_{4,4}(z)\) is \(\tc_{4,4}(a-1)\).
Since the Coxeter number of \(E_6\) is \(12\), the corresponding
truncation of \(\tC_{4,4}(t)\) through degree \(12\) is  
\[
\tC_{4,4}^{[12]}(t)
=
t+2t^3+3t^5+3t^7+2t^9+t^{11}.
\]
Since \(\tc_{4,4}^{[12]}(5)=3\),
\eqref{eq: ADE quantum affine KR d-value} gives
\[
\de\bigl(\Vkm{4}_p, \Vkm{4}_s \bigr) = 3 \qquad\text{if}\qquad  |p-s|=6.
\]
In particular, the multiplicity \(3\) of the root \(q^6\) agrees
with the symmetrized inverse \(t\)-Cartan coefficient prescribed by
Theorem~\ref{thm: ADE quantum affine KR d-values}.

For arbitrary \(i,j\in I\) and \(k,\ell\in\Z_{>0}\), the exceptional
\(\de\)-values are computed uniformly by
\eqref{eq: ADE quantum affine KR d-value}, with
\(K=\max(k,\ell)\), \(L=\min(k,\ell)\), and \(u=|p-s|\).
\end{remark}

We now return to KR-modules over quiver Hecke algebras of non-simply-laced finite
type.  Keep the finite-type folding datum
\[
(\hg,\Gamma)\longrightarrow\g
\]
of Section~\ref{sec: folding cluster algebras}, where \(\hg\) is of
simply-laced finite type and \(\g\) is of finite type \(B\), \(C\), \(F\), or
\(G\) (\eqref{eq: finite folding}).  Let \(\hQ\) be a \(\Gamma\)-foldable simply-laced Dynkin quiver 
and set
\[
Q\seteq\hQ/\Gamma.
\]
Keep the notation of
Lemma~\ref{lem: folding of adapted repetition quivers}.  In
particular, the resulting vertex fibers induce the chain-dependent
box fibers used in
Theorem~\ref{thm: categorical folding of d-invariants} and
Corollary~\ref{cor: folding of d-values between box modules}.

The Coxeter numbers of \(\hg\) and \(\g\) have the same value \(h\).  We set
\[
\tC_{i,j}^{\g,[h]}(t)
\seteq
\sum_{v=1}^{h}\tc_{i,j}^{\g}(v)t^v.
\]
For \(\im,\jm\in\hI\), we likewise set
\[
\tC_{\im,\jm}^{\hg,[h]}(t)
\seteq
\sum_{v=1}^{h}\tc_{\im,\jm}^{\hg}(v)t^v.
\]
For these standard finite-type foldings one has \(m=|\Gamma|\), so the
group-sum form of
Theorem~\ref{thm: folding of inverse t-Cartan matrices} has no
additional scalar factor.
\begin{theorem}\label{thm: BCFG quiver Hecke KR d-values}
Assume that \(\g\) is of finite type, and let \(h\) be its
Coxeter number.
Let \(i,j\in I\), \(k,\ell\in\Z_{>0}\), and \(p,s\in\Z\).
Suppose that all vertices
\[
(i,p+k-1-2a)\quad(0\leq a<k),
\qquad
(j,s+\ell-1-2b)\quad(0\leq b<\ell)
\]
belong to \(\Upsilon_Q\).
Set
\[
K\seteq\max(k,\ell),\qquad
L\seteq\min(k,\ell),\qquad
u\seteq|p-s|.
\]
Then
\begin{equation}
\label{eq: BCFG quiver Hecke KR d-value}
\de\bigl(\Vkm{i^k}_p,\Vkm{j^\ell}_s\bigr)
=[t^u]\bigl(t^K[L]_t\tC_{i,j}^{\g,[h]}(t)\bigr)
=\sum_{a=0}^{L-1}
\tc_{i,j}^{\g,[h]}(u-|k-\ell|-2a-1),
\end{equation}
where
\[
\tc_{i,j}^{\g,[h]}(v)
\seteq
\begin{cases}
\tc_{i,j}^{\g}(v),&1\leq v\leq h,\\
0,&\text{otherwise}.
\end{cases}
\]
\end{theorem}

\begin{proof}
In simply-laced type, the assertion follows from
Theorem~\ref{thm: ADE quantum affine KR d-values}
and generalized quantum affine Schur--Weyl duality.
For the remaining types, use the standard finite-type folding
\(\hg\to\g\), for which the Coxeter numbers agree.

Let \(\frakc,\frakc'\) be the boxes corresponding to
\(\Vkm{i^k}_p,\Vkm{j^\ell}_s\), respectively.
Choose admissible chains \(\frakC,\frakD\) containing these boxes,
and choose
\(\hfrakc'\in\Fib_{\frakD}(\frakc')\) of color \(\jm\in I_j\).
By Corollary~\ref{cor: folding of d-values between box modules},
\[
\de\bigl(M^\ii[\frakc],M^\ii[\frakc']\bigr)
=d_i\sum_{\hfrake\in\Fib_{\frakC}(\frakc)}
\de\bigl(M^\hii[\hfrake],M^\hii[\hfrakc']\bigr).
\]
By Lemma~\ref{lem: folding preserves KR box parameters},
 the modules \(M^\hii[\hfrake]\) with
\(\hfrake\in\Fib_{\frakC}(\frakc)\) have parameters \(k,p\), and
\(M^\hii[\hfrakc']\) has parameters \(\ell,s\).
All their defining vertices lie in \(\Upsilon_{\hQ}^{\hii}\).
Their colors run through \(I_i\) and the fixed color \(\jm\),
respectively.
Applying the simply-laced formula and
Theorem~\ref{thm: folding of inverse t-Cartan matrices}, we obtain
\[
\begin{aligned}
\de\bigl(\Vkm{i^k}_p,\Vkm{j^\ell}_s\bigr)
=[t^u]\left(
t^K[L]_t\,d_i\sum_{\im\in I_i}
\tC_{\im,\jm}^{\hg,[h]}(t)
\right) =[t^u]\bigl(t^K[L]_t\tC_{i,j}^{\g,[h]}(t)\bigr).
\end{aligned}
\]
Finally,
\[
t^K[L]_t
=\sum_{a=0}^{L-1}t^{|k-\ell|+2a+1},
\]
which gives the coefficient sum.
\end{proof}

\begin{example}\label{ex: A5 C3 KR coefficient comparison}
Return to the folding \(A_5\to C_3\) of
Examples~\ref{ex: A5 C3 folding of F-invariants}
and~\ref{ex: A5 C3 t-Cartan folding}.  For the reduced sequence
\(\ii\) in the first example, the modules
\[
M=\Vkm{2^2}_3=M^\ii[2,5],\qquad
N=\Vkm{3^2}_6=M^\ii[7,9]
\]
have defining vertices \((2,2),(2,4)\) and \((3,5),(3,7)\) in the
same \(\Upsilon_Q\); see also
Example~\ref{ex: one-sided repetition quivers}~(ii).
Here \(K=L=2\), \(u=3\), \(I_2=\{2,4\}\), \(I_3=\{3\}\), and
\(d_2=1\).  Theorem~\ref{thm: folding of inverse t-Cartan matrices}
gives
\[
\tc_{2,3}^{C_3}(2)
=\tc_{2,3}^{A_5}(2)+\tc_{4,3}^{A_5}(2)=2.
\]
Since \(h=6\), we have
\[
\tc_{i,j}^{[6]}(v)=\tc_{i,j}(v)
\qquad(v=0,2),
\]
where both coefficients at \(v=0\) vanish.  Thus
Theorem~\ref{thm: BCFG quiver Hecke KR d-values} gives
\[
\de(M,N)
=\tc_{2,3}^{C_3,[6]}(2)+\tc_{2,3}^{C_3,[6]}(0)=2,
\]
in agreement with the categorical calculation following
Corollary~\ref{cor: folding of d-values between box modules}.
\end{example}

Recall the counting functions \(\sfT^{(d)}_{r,L}(\nu)\) from
\eqref{eq: closed T}.

\begin{corollary}\label{cor: type C quiver Hecke KR closed formula}
Let \(\g\) be of type \(C_n\), with \(n\) the long node.  Keep the
notation
\[
L\seteq\min(k,\ell), \qquad u\seteq|p-s|, \qquad \kappa\seteq|k-\ell|.
\]
Assume that all the vertices defining the two quiver Hecke
KR-modules in question lie in \(\Upsilon_Q\).
\bnum
\item If \(1\leq i,j<n\), set
\[
\nu_1\seteq\frac{u-\kappa-|i-j|}{2}
\qtq
\nu_2\seteq\frac{u-\kappa-(2n-i-j)}{2}.
\]
Then we have
\[
\de\bigl( \Vkm{i^k}_p, \Vkm{j^\ell}_s \bigr)
=
\sfT^{(1)}_{\brho^{(2n-1)}_{i,j},L}(\nu_1)
+
\sfT^{(1)}_{\brho^{(2n-1)}_{2n-i,j},L}(\nu_2).
\]
\item If \(i=n\) or \(j=n\), set
\[
\nu\seteq \frac{u-\kappa-|i-j|}{2}.
\]
Then we have
\[
\de\bigl(\Vkm{i^k}_p,\Vkm{j^\ell}_s\bigr)
=
2\sfT^{(1)}_{\brho^{(2n-1)}_{i,j},L}(\nu). 
\]
\ee 
\end{corollary}

\begin{proof}
Consider the standard folding
\[
A_{2n-1}\longrightarrow C_n, \qquad I_i=\{i,2n-i\}\quad(1\leq i<n), \qquad I_n=\{n\}.
\]
Let \(M\) and \(N\) denote the two downstairs KR-modules in the
statement.  For \(\im\in I_i\) and \(\jm\in I_j\), denote their
corresponding lifted KR-modules by
\(\hM_{\im}\) and \(\hN_{\jm}\), respectively.
By Lemma~\ref{lem: folding preserves KR box parameters},
\(\hM_\im\simeq\Vkm{\im^k}_p\) and
\(\hN_\jm\simeq\Vkm{\jm^\ell}_s\). 

Suppose first that \(i,j<n\).  Identify the upstairs vertices with
the standard numeric labels of \(A_{2n-1}\), and choose the lifts
\(\im=i\in I_i\) and \(\jm=j\in I_j\).  Thus the other lift of
\(i\) is \(2n-i\).  Since
\(d_i=1\), Corollary~
\ref{cor: folding of d-values between box modules} gives
\[
\de(M,N) = \de\bigl(\hM_\im,\hN_\jm\bigr) + \de\bigl(\hM_{2n-i},\hN_\jm\bigr).
\]
For the two pairs of upstairs colors, we have
\[
|\im-\jm|=|i-j|
\qtq
|(2n-i)-\jm|=2n-i-j,
\]
and
\[
\brho^{(2n-1)}_{i,j}
=
\brho^{(2n-1)}_{2n-i,j}
=
\min(i,j).
\]
Applying Corollary~
\ref{cor: type A closed d-value formula} to the two summands yields
\[
\de(M,N)
=
\sfT^{(1)}_{\brho^{(2n-1)}_{i,j},L}
\left( \frac{u-|i-j|-|k-\ell|}{2} \right)
+
\sfT^{(1)}_{\brho^{(2n-1)}_{2n-i,j},L}
\left( \frac{u-(2n-i-j)-|k-\ell|}{2} \right),
\]
which is the first formula.

Suppose now that \(i=n\).  Since \(I_n=\{n\}\) and \(d_n=2\),
Corollary~\ref{cor: folding of d-values between box modules} gives
\[
\de(M,N) = 2\de\bigl(\hM_n,\hN_{\jm}\bigr) \qquad(\jm\in I_j).
\]
The type \(A_{2n-1}\) closed formula therefore gives
\[
\de(M,N)
=
2\sfT^{(1)}_{\brho^{(2n-1)}_{i,j},L}
\left( \frac{u-|i-j|-|k-\ell|}{2} \right).
\]

Finally, suppose that \(j=n\) and \(i<n\).  Then
\[
\de(M,N) = \de\bigl(\hM_i,\hN_n\bigr) + \de\bigl(\hM_{2n-i},\hN_n\bigr).
\]
The two summands are equal, since
\[
|i-n|=|(2n-i)-n|
\qtq
\brho^{(2n-1)}_{i,n}
=
\brho^{(2n-1)}_{2n-i,n}
=i.
\]
Consequently,
\[
\de(M,N)
=
2\sfT^{(1)}_{\brho^{(2n-1)}_{i,j},L}
\left( \frac{u-|i-j|-|k-\ell|}{2} \right)
=
2\sfT^{(1)}_{\brho^{(2n-1)}_{i,j},L}(\nu).
\]
This proves the result.
\end{proof}

Example~\ref{ex: A5 C3 folding of F-invariants} is the case
\(n=3\), \(i=2\), \(j=3\), \(k=\ell=2\), and \(u=3\) of the
second formula:
\[
2\sfT^{(1)}_{2,2}(1)=2.
\]

\begin{example}\label{ex: type C5 folded KR d-value}
Let \(\g\) be of type \(C_5\), with node \(5\) long, and suppose
that all the vertices defining the modules below lie in
\(\Upsilon_Q\).  Consider
\[
\Vkm{5^2}_p \qtq \Vkm{5^3}_s \qquad\text{with}\qquad |p-s|=5.
\]
Here $L=2$, $\kappa=1$ and $\nu=\dfrac{5-1}{2}=2$. 

Hence
\[
\de\bigl( \Vkm{5^2}_p, \Vkm{5^3}_s \bigr)
= 2\sfT^{(1)}_{5,2}(2) = 4.
\]

This value may also be read directly from the folding
\[
A_9\longrightarrow C_5.
\]
The middle vertex \(5\) of \(A_9\) is fixed by the diagram
involution and lies above the long vertex \(5\) of \(C_5\).  Thus
its fiber has cardinality \(1\), while the corresponding
symmetrizer is
\[
d_5=2.
\]
The type \(A_9\) formula gives
\[
\de\bigl(
\Vkm{5^2}_p,
\Vkm{5^3}_s
\bigr)
=
\sfT^{(1)}_{5,2}(2)
=
2,
\]
and categorical folding therefore gives the downstairs type \(C_5\)
value
\[
\de\bigl(
\Vkm{5^2}_p,
\Vkm{5^3}_s
\bigr)
=
2\cdot2
=
4.
\]
\end{example}

The type \(B_n\) formula is obtained in the same way from the
standard folding \(D_{n+1}\to B_n\), using
Corollary~\ref{cor: folding of d-values between box modules} and the
type \(D\) closed formula.  We therefore omit its proof.

\begin{corollary}\label{cor: type B quiver Hecke KR closed formula}
Let \(\g\) be of type \(B_n\), with \(n\) the short node.  Keep the
notation
\[
L\seteq\min(k,\ell), \qquad u\seteq|p-s|, \qquad \kappa\seteq|k-\ell|.
\]
Assume that all the vertices defining the two quiver Hecke KR-modules in question lie in \(\Upsilon_Q\).
\bnum
\item If \(1\leq i,j<n\), set
\[
r\seteq\min(i,j), \qquad \nu_1\seteq\frac{u-\kappa-|i-j|}{2}, \qquad \nu_2\seteq\frac{u-\kappa-(2n-i-j)}{2}.
\]
Then we have
\[
\de\bigl( \Vkm{i^k}_p, \Vkm{j^\ell}_s \bigr) = 2\sfT^{(1)}_{r,L}(\nu_1) + 2\sfT^{(1)}_{r,L}(\nu_2).
\]
\item  If exactly one of \(i,j\) is \(n\), let
\[
c\seteq\min(i,j) \qtq \nu\seteq \frac{u-\kappa-(n-c)}{2}. 
\]
Then we have
\[
\de\bigl( \Vkm{i^k}_p, \Vkm{j^\ell}_s \bigr) = 2\sfT^{(1)}_{c,L}(\nu). 
\]
\item  Finally, if \(i=j=n\), set
\[
r_1\seteq \left\lfloor\frac{n+1}{2}\right\rfloor, \quad r_2\seteq \left\lfloor\frac{n}{2}\right\rfloor,
\quad
\nu_1\seteq\frac{u-\kappa+2}{2}, \quad \nu_2\seteq\frac{u-\kappa}{2}.
\]
Then
\[
\de\bigl( \Vkm{n^k}_p, \Vkm{n^\ell}_s \bigr) = \sfT^{(2)}_{r_1,L}(\nu_1) + \sfT^{(2)}_{r_2,L}(\nu_2).
\]
\ee 
\end{corollary}

\begin{example}\label{ex: type B5 folded KR d-value}
Let \(\g\) be of type \(B_5\), with node \(5\) short, and suppose
that all the vertices defining the modules below lie in
\(\Upsilon_Q\).  Consider
\[
\Vkm{2^2}_p \qtq \Vkm{3^3}_s \qquad\text{with}\qquad |p-s|=6.
\]
In this case,
\[
L=2,\quad \kappa=1, \quad
\nu_1 = \frac{6-1-|2-3|}{2} = 2,
\quad
\nu_2 = \frac{6-1-(10-2-3)}{2} = 0.
\]
Therefore
\[
\begin{aligned}
\de\bigl( \Vkm{2^2}_p, \Vkm{3^3}_s \bigr) = 2\sfT^{(1)}_{2,2}(2) + 2\sfT^{(1)}_{2,2}(0)= 4.
\end{aligned}
\]

This value is obtained from the folding
\[
D_6\longrightarrow B_5.
\]
The vertices \(2\) and \(3\) are fixed by the involution exchanging
the two spin vertices of \(D_6\), and their corresponding
symmetrizer is \(2\).  The type \(D_6\) formula gives
\[
\de\bigl(
\Vkm{2^2}_p,
\Vkm{3^3}_s
\bigr)
=
2.
\]
Consequently, the categorical fiber formula yields the downstairs
type \(B_5\) value
\[
\de\bigl(
\Vkm{2^2}_p,
\Vkm{3^3}_s
\bigr)
=
2\cdot2
=
4.
\]
\end{example}

Theorem~\ref{thm: BCFG quiver Hecke KR d-values} shows that the
categorical folding of \(\de\)-invariants and the folding of the
symmetrized inverse \(t\)-Cartan matrix coincide.   Both formulas contain the factor \(d_i\).  

\begin{example}\label{ex: type F4 folded KR d-value}
Consider type \(F_4\) with the convention
\[
C=
{\scriptsize\begin{pmatrix}
2&-1&0&0\\
-1&2&-1&0\\
0&-2&2&-1\\
0&0&-1&2
\end{pmatrix}},
\qquad
D=\diag(2,2,1,1).
\]
Truncating through degree \(12\), we obtain 
\[
\tC_{2,3}^{[12]}(t)
=
2t^2+4t^4+4t^6+4t^8+2t^{10}.
\]

Suppose that
\[
\Vkm{2^2}_p
\qtq
\Vkm{3^3}_s
\]
are defined by \(\bfi\)-boxes whose defining vertices all lie in
\(\Upsilon_Q\) and that
\[
|p-s|=6.
\] 
Here \(k=2\) and \(\ell=3\), so \(L=2\) and \(|k-\ell|=1\).
Theorem~\ref{thm: BCFG quiver Hecke KR d-values} gives
\[
\begin{aligned}
\de\bigl( \Vkm{2^2}_p, \Vkm{3^3}_s \bigr) = \tc_{2,3}^{[12]}(4) + \tc_{2,3}^{[12]}(2) = 4+2 = 6.
\end{aligned}
\]

Under the folding
$E_6\longrightarrow F_4$,
the fibers are
$I_1=\{2\},I_2=\{4\}, 
I_3=\{3,5\}, I_4=\{1,6\}$.
Thus this equality is obtained from the singleton vertex fiber over
\(2\), with folding factor \(d_2=2\), by summing the corresponding
type \(E_6\) KR-module \(\de\)-values.
Thus it is the higher KR analogue of the cuspidal folding formula.
\end{example}

\begin{example}\label{ex: type G2 folded KR d-value}
Consider type \(G_2\) with the convention
\[
C=
\begin{pmatrix}
2&-3\\
-1&2
\end{pmatrix},
\qquad
D=\diag(1,3).
\]
The truncation through degree \(6\) of the off-diagonal entry
of the symmetrized inverse \(t\)-Cartan matrix is 
\[
\tC_{1,2}^{[6]}(t) = 3t^2+3t^4.
\]

Suppose that
\[
\Vkm{1^2}_p \qtq \Vkm{2^2}_s
\]
are defined by \(\bfi\)-boxes whose defining vertices all lie in
\(\Upsilon_Q\) and that
$|p-s|=3$. 

Since $L=2$ and $|2-2|=0$, 
Theorem~\ref{thm: BCFG quiver Hecke KR d-values} yields
\[
\de\bigl( \Vkm{1^2}_p, \Vkm{2^2}_s \bigr) = \tc_{1,2}^{[6]}(2) + \tc_{1,2}^{[6]}(0) = 3.
\]

This value is obtained from the folding
\[
D_4\longrightarrow G_2
\] 
by summing the corresponding type \(D_4\) KR-module
\(\de\)-values over the three outer vertices above color \(1\),
while keeping the central vertex above color \(2\) fixed.
Here \(d_1=1\). 
\end{example}

\begin{remark} \label{rmk: Lambda computable}
For two quiver Hecke KR-modules \(M,N\) defined by boxes
in the same \(Q\)-adapted reduced sequence,
Proposition~\ref{prop: KKOP directional box Lambda formula}
computes at least one of \(\Lambda(M,N)\) and
\(\Lambda(N,M)\), as explained in
Remark~\ref{rmk: one Lambda is computable}. 
\eqref{eq: BCFG quiver Hecke KR d-value}
then determines the other through
\[
\Lambda(M,N)+\Lambda(N,M)=2\de(M,N).
\]
\end{remark}

\begin{example}\label{ex: overlapping boxes directional Lambda}
In type \(C_5\), take \(\xi_i=i-1\) and the
\(Q\)-adapted reduced sequence \(\ii=(1,2,3,4,5)^5\).
We use the standard realization of \(C_5\) in
\(\mathbb R^5\) with orthonormal basis \(\bfe_1,\ldots,\bfe_5\).
The positions \(4,9,14\) have color \(4\), with
coordinates \(3,5,7\) in the repetition quiver and root labels
\[
\be_4^\ii=\bfe_1-\bfe_5,\qquad
\be_9^\ii=\bfe_1+\bfe_2,\qquad
\be_{14}^\ii=\bfe_2+\bfe_3.
\]
Set
\[
M=\Vkm{4^2}_{6}=M^\ii[9,14],
\qquad
N=\Vkm{4^2}_{4}=M^\ii[4,9].
\]
Since \(9>4^-=-\infty\),
Proposition~\ref{prop: KKOP directional box Lambda formula}
applies to \((M,N)\).
The two boxes share position \(9\), and
\(\lambda^\ii(\ve_9,\ve_9)=0\).  Hence
\[
\begin{aligned}
\lambda^\ii(\calE_{[9,14]},\calE_{[4,9]})
&=(\be_9^\ii,\be_4^\ii)_\g
  +(\be_{14}^\ii,\be_4^\ii)_\g
  +(\be_{14}^\ii,\be_9^\ii)_\g\\
&=1+0+1=2.
\end{aligned}
\]
Thus \(\Lambda(M,N)=-2\).
Since \(k=\ell=2\) and \(u=2\),
Theorem~\ref{thm: BCFG quiver Hecke KR d-values} gives
\[
\de(M,N)=\tc_{4,4}^{[10]}(1)
+\tc_{4,4}^{[10]}(-1)=d_4=1.
\]
Therefore,
\[
\Lambda(N,M)=2\cdot1-(-2)=4.
\]
\end{example}

\subsection{Beyond finite type}
\label{subsec: beyond finite type}

We conclude by returning to arbitrary symmetrizable Cartan data.
Although the uniform formula of
Theorem~\ref{thm: BCFG quiver Hecke KR d-values} is proved here only
in finite type, its coefficient expression continues to detect three
elementary local categorical interactions without any finite-type
assumption: the interaction across the boundary of a box, the exchange
interaction in a \(T\)-system, and the vanishing associated with
commuting boxes.

 Let \(\ii=(i_1,\ldots,i_l)\) be a reduced sequence, and let
\(\frakc=[x,y]\) be an \(\bfi\)-box of color \(i\) associated
with \(\ii\).   By
\cite[Lemma~5.8]{KK24}, we have
\begin{align}\label{eq: boundary}
\de\bigl(M^\ii[x,y],S^\ii_{y^+}\bigr)=d_i
\quad \text{if }y^+\leq l \qtq  \de\bigl(S^\ii_{x^-},M^\ii[x,y]\bigr)=d_i
\quad \text{if }x^-\geq1.
\end{align}
We next record a coordinate estimate for commuting
\(\bfi\)-boxes. 

\begin{lemma}
\label{lem: commuting boxes coordinate inequality}
Let \(\ii\) be a finite \(Q\)-adapted sequence, and let
\(\frakc=[x,y]\) and \(\frakc'=[x',y']\) be commuting
\(\bfi\)-boxes of colors \(i\) and \(j\), respectively.  Let \(k\) and \(\ell\) be the numbers of occurrences of colors
\(i\) and \(j\) in \([x,y]\) and \([x',y']\), respectively, and set
\(p=(p_x^\ii+p_y^\ii)/2\) and
\(s=(p_{x'}^\ii+p_{y'}^\ii)/2\).
Then
\[
|p-s|\leq |k-\ell|+d(i,j).
\]
When \(\ii\) is reduced, the corresponding KR-modules are 
\[
\Vkm{i^k}_p
\qtq
\Vkm{j^\ell}_s.
\]
\end{lemma}

\begin{proof}
 Write \(P_a=p_a^\ii\).
If \(a<b\), let \(\eta\) be the height function immediately
before position \(b\) is read.  Since position \(a\) has
already been read, we have
 \begin{equation}\label{eq: position height inequality}
P_a+2\leq\eta_{i_a}
\leq\eta_{i_b}+d(i_a,i_b)
=P_b+d(i_a,i_b).
\end{equation}
Here the second inequality follows from
\(|\eta_i-\eta_j|=1\) for adjacent vertices \(i,j\).

Suppose first that
\[
x^-<x'\leq y'<y^+.
\]
We have
\[
P_x=p-k+1,\quad P_y=p+k-1
\qtq{}
P_{x'}=s-\ell+1,\quad P_{y'}=s+\ell-1.
\]
 If \(x^-\) exists, then \(P_{x^-}=p-k-1\), and
\eqref{eq: position height inequality} gives
\[
p-s\leq k-\ell+d(i,j).
\] 
If \(x^-=-\infty\), then \(P_x=\xi_i\) and
\(P_{x'}\geq\xi_j\), so the same inequality follows by applying
the preceding shortest-path argument to the initial height function.

 Similarly, if \(y^+\) exists, then \(P_{y^+}=p+k+1\), and
\eqref{eq: position height inequality} gives
\[
s-p\leq k-\ell+d(i,j).
\]
If \(y^+=+\infty\), then, for the height function obtained after
the finite sequence has been read, the heights at \(i\) and \(j\)
are \(P_y+2\) and at least \(P_{y'}+2\), respectively.  The same
shortest-path argument therefore gives the inequality.

Consequently,
\[
|p-s|\leq k-\ell+d(i,j).
\]
For the other commutation condition
\[
(x')^-<x\leq y<(y')^+,
\]
interchanging the two boxes gives
\[
|p-s|\leq \ell-k+d(i,j). \qedhere
\]
\end{proof}

We also recall the following consequence of the determinantial
module \(T\)-system.

\begin{proposition}[{\cite[Lemma~5.18]{KK24}}]
\label{prop: t-system exchange pair d-value}
Let \(\ii=(i_1,\ldots,i_l)\) be a reduced sequence, and let
\([x,y]\) be an \(\bfi\)-box of color \(i\) with \(x<y\).  Then the
two determinantial modules associated with the \(T\)-system
exchange pair \([x^+,y]\) and \([x,y^-]\) satisfy
\[
\de\bigl(M^\ii[x^+,y],M^\ii[x,y^-]\bigr)=d_i.
\]
\end{proposition}

\begin{proposition}
\label{prop: local coefficient formula beyond finite type}
Let the Cartan datum be symmetrizable, and let the indicated
\(\bfi\)-boxes occur in a \(Q\)-adapted reduced sequence \( \ii \).
\bnum
\item \label{it: adj singleton} If an \(\bfi\)-box of color \(i\) corresponds to 
\(\Vkm{i^k}_p\) and  the singleton box at the next occurrence of color \(i\)
corresponds to \(\Vkm{i}_{p+k+1}\), then
\[
\de\bigl(\Vkm{i^k}_p,\Vkm{i}_{p+k+1}\bigr)  = [t^{k+1}] \left(t^k[1]_t\tC_{i,i}(t) \right) = d_i.
\]

\item \label{it: T-system exchange pair}
If \(\Vkm{i^k}_p\) and \(\Vkm{i^k}_s\) correspond to
a \(T\)-system exchange pair of \(\bfi\)-boxes, then
\(|p-s|=2\) and
\[
\de\bigl(\Vkm{i^k}_p,\Vkm{i^k}_s\bigr)
=
[t^2]\left(t^k[k]_t\tC_{i,i}(t)\right)
=d_i.
\]

\item \label{it: commuting}
If two commuting \(\bfi\)-boxes of colors \(i\) and \(j\) correspond,
respectively, to
\(\Vkm{i^k}_p\) and \(\Vkm{j^\ell}_s\), then
\[
\de\bigl(\Vkm{i^k}_p,\Vkm{j^\ell}_s\bigr) = [t^{|p-s|}] \left( t^K[L]_t\tC_{i,j}(t) \right) = 0,
\]
where
\[
K=\max(k,\ell), \qquad L=\min(k,\ell). 
\]
\ee
\end{proposition}

\begin{proof}
For~\eqref{it: adj singleton}, the categorical equality is given in~\eqref{eq: boundary}.  Since \([1]_t=1\), we have
\[
[t^{k+1}]
\left(
t^k[1]_t\tC_{i,i}(t)
\right)
=
[t^1]\tC_{i,i}(t)
=
\tc_{i,i}(1)
=
d_i.
\]

For~\eqref{it: T-system exchange pair}, the categorical equality
follows from
Proposition~\ref{prop: t-system exchange pair d-value}.  Indeed, the
two subboxes \([x^+,y]\) and \([x,y^-]\) each contain \(k\)
occurrences of color \(i\), 
and their centered coordinates differ by \(2\).
On the coefficient side,
\[
t^k[k]_t=t+t^3+\cdots+t^{2k-1}.
\]
Since \(\tC_{i,i}(t)\in t\Z[[t]]\) and
\([t]\tC_{i,i}(t)=\tc_{i,i}(1)=d_i\), it follows that
\[
[t^2]\left(t^k[k]_t\tC_{i,i}(t)\right)
=
[t]\tC_{i,i}(t)
=
\tc_{i,i}(1)
=d_i.
\]

For~\eqref{it: commuting}, Lemma~\ref{lem: commuting box modules} implies that the
corresponding quiver Hecke KR-modules strongly commute.  Hence
\[
\de\bigl(\Vkm{i^k}_p,\Vkm{j^\ell}_s\bigr)=0.
\]

Set
\[
u=|p-s|, \qquad K=\max(k,\ell), \qquad L=\min(k,\ell).
\]
By Lemma~\ref{lem: commuting boxes coordinate inequality}, we have
$u\leq |k-\ell|+d(i,j)$.

Since $t^K[L]_t = \sum_{a=0}^{L-1}t^{K-L+1+2a}$, 
we obtain
\[
[t^u]
\left( t^K[L]_t\tC_{i,j}(t) \right) = \sum_{a=0}^{L-1} \tc_{i,j}\bigl(u-|k-\ell|-2a-1\bigr).
\]
For every \(0\leq a\leq L-1\), we have
\[
u-|k-\ell|-2a-1
\leq
d(i,j)-2a-1
<
d(i,j).
\]
Therefore, Proposition~\ref{prop: coefficients of inverse t-Cartan}~\eqref{it: 0-range}
shows that every summand on the right-hand side vanishes:
\[
\tc_{i,j}\bigl(u-|k-\ell|-2a-1\bigr)=0.
\]
Consequently,
\[
[t^u] \left( t^K[L]_t\tC_{i,j}(t) \right) =0. \qedhere
\]
\end{proof}

Thus the \(\de\)-invariant in each of these three cases is given by
the coefficient formula with the full \(t\)-adic entry
\(\tC_{i,j}(t)\), without a finite-type assumption. 

\begin{example}
\label{ex: affine B3 commuting boxes of different colors}
Consider type \(B_3^{(1)}\) with \(D=\diag(2,2,2,1)\).
Choose the height function
\( (\xi_0,\xi_1,\xi_2,\xi_3)=(1,1,2,1) \).
The sequence
\[
\ii=(0,1,3,2,0,1,3,2,0,1,3,2,0)
\]
is \(Q\)-adapted and reduced.  Consider the boxes
$[1,5]$ and $[4,8]$
of colors \(0\) and \(2\), respectively. 
Their corresponding quiver Hecke KR-modules are
\[
M^\ii[1,5]\simeq\Vkm{0^2}_2 \qtq M^\ii[4,8]\simeq\Vkm{2^2}_3.
\]

The boxes are not nested, since \(1<4\leq5<8\).
Moreover, \(1^-<4\leq8<5^+=9\), so they commute
by the definition of commuting \(\bfi\)-boxes. Thus Lemma~\ref{lem: commuting box modules} gives
\[
\de\bigl(\Vkm{0^2}_2,\Vkm{2^2}_3\bigr)=0.
\]

The same value is obtained directly from the full \(t\)-adic
symmetrized inverse \(t\)-Cartan matrix.  Indeed, \(K=L=2\) and
\(|2-3|=1\), and therefore
\[
[t^1]\left(t^2[2]_t\tC_{0,2}^{B_3^{(1)}}(t)\right)
=\tc_{0,2}^{B_3^{(1)}}(-2)+\tc_{0,2}^{B_3^{(1)}}(0)=0.
\]
Here the last equality follows from
Proposition~\ref{prop: coefficients of inverse t-Cartan}, since both
coefficient indices are at most \(d(0,2)=1\).
This illustrates
Proposition~\ref{prop: local coefficient formula beyond finite type}
in affine type, using two non-nested boxes of different colors.
\end{example}

\bibliographystyle{amsplain}
\bibliography{ref}{}
\end{document}